\documentclass[12pt,vi,twoside]{mitthesis}
\usepackage{lgrind}
\usepackage{amssymb,amsmath}
\usepackage[dvips]{graphicx}

\newcommand{\hR}{\hat{R}}
\newcommand{\cxp}{{\rm cxp}\,}
\newcommand{\Linf}{L_{\infty}}

\newcommand{\bul}{{\bullet}}
\newcommand{\al}{{\alpha}}
\newcommand{\la}{{\lambda}}
\newcommand{\h}{{\hbar}}
\newcommand{\pia}{{\pi_{{\rm add}}}}

\newcommand{\mh}{{\mathfrak{h}}}
\newcommand{\mv}{{\mathfrak{v}}}
\newcommand{\mP}{{\mathfrak{P}}}
\newcommand{\ma}{{\mathfrak{a}}}
\newcommand{\mb}{{\mathfrak{b}}}

\newcommand{\md}{{\mathfrak{d}}}
\newcommand{\mG}{{\mathfrak{G}}}
\newcommand{\mC}{{\mathfrak{C}}}
\newcommand{\mU}{{\mathfrak{U}}}

\newcommand{\mH}{{\mathfrak{H}}}

\newcommand{\bb}{{\bf b}}

\newcommand{\hmb}{\hat{\mb}}
\newcommand{\dia}{\diamond}
\newcommand{\club}{\clubsuit}

\newcommand{\Om}{{\Omega}}

\newcommand{\si}{{\sigma}}
\newcommand{\ga}{{\gamma}}
\newcommand{\vf}{{\varphi}}
\newcommand{\ve}{{\varepsilon}}
\newcommand{\ka}{{\kappa}}
\newcommand{\vr}{{\varrho}}
\newcommand{\G}{{\Gamma}}
\newcommand{\pa}{{\partial}}

\newcommand{\M}{{\cal M}}
\newcommand{\N}{{\cal N}}
\newcommand{\cO}{{\cal O}}
\newcommand{\cF}{{\cal F}}
\newcommand{\cD}{{\cal D}}
\newcommand{\cT}{{\cal T}}

\newcommand{\cA}{{\cal A}}
\newcommand{\cG}{{\cal G}}
\newcommand{\cH}{{\cal H}}

\newcommand{\cS}{{\cal S}}
\newcommand{\cE}{{\cal E}}
\newcommand{\cB}{{\cal B}}
\newcommand{\cL}{{\cal L}}
\newcommand{\cC}{{\cal C}}
\newcommand{\cK}{{\cal K}}
\newcommand{\cU}{{\cal U}}
\newcommand{\cR}{{\cal R}}
\newcommand{\cV}{{\cal V}}

\newcommand{\bbK}{{\Bbb K}}
\newcommand{\bbA}{{\Bbb A}}

\newcommand{\bbC}{{\Bbb C}}
\newcommand{\bbR}{{\Bbb R}}
\newcommand{\bbZ}{{\Bbb Z}}
\newcommand{\bbRf}{{{\Bbb R}^d_{formal}}}

\newcommand{\n}{{\nabla}}

\newcommand{\de}{{\delta}}
\newcommand{\D}{{\Delta}}

\newcommand{\tPi}{{\widetilde{\Pi}}}
\newcommand{\tn}{{\widetilde{\nabla}}}

\newcommand{\tcL}{{\widetilde{\cL}}}
\newcommand{\tcU}{{\widetilde{\cU}}}
\newcommand{\tH}{{\widetilde{H}}}
\newcommand{\tM}{{\widetilde{M}}}

\newcommand{\tF}{{\widetilde{F}}}

\newcommand{\tG}{{\widetilde{\Gamma}}}

\newcommand{\erarrow}{\stackrel{\sim}{\rightarrow}}

\newcommand{\lrarrow}{\,\longrightarrow \,}
\newcommand{\llarrow}{\,\longleftarrow \,}
\newcommand{\brarrow}{\succ\rightarrow}
\newcommand{\blarrow}{\leftarrow\prec}
\newcommand{\bbrarrow}{\succ\succ\rightarrow}
\newcommand{\bblarrow}{\leftarrow\prec\prec}

\newcommand{\cVf}{{\cal V}^{fib}}
\newcommand{\OM}{C^{\infty}(M)}
\newcommand{\SM}{{\cal S}M}
\newcommand{\AM}{{\cal A}^{\bul}(M)}
\newcommand{\Tp}{T_{poly}^{\bul}(M)}
\newcommand{\Dp}{D_{poly}^{\bul}(M)}
\newcommand{\Cp}{C^{poly}_{\bul}(M)}
\newcommand{\JM}{J_{\bul}(M)}

\newcommand{\cTp}{\cT_{poly}}
\newcommand{\cDp}{\cD_{poly}}
\newcommand{\cCp}{\cC^{poly}}

\newcommand{\Omb}{\Om^{\bul}}
\newcommand{\OmS}{\Om^{\bul}(M,\SM)}
\newcommand{\OmT}{\Om^{\bul}(M,\cT_{poly})}
\newcommand{\OmD}{\Om^{\bul}(M,\cD_{poly})}
\newcommand{\OmC}{\Om^{\bul}(M,\cC^{poly})}
\newcommand{\OmE}{\Om^{\bul}(M,\cE)}
\newcommand{\FT}{ker\, \de \cap \G(M, \cT_{poly})}
\newcommand{\FD}{ker\, \de \cap \G(M, \cD_{poly})}

\newcommand{\inu}{{\nu^{-1}}}
\date{}
\newtheorem{defi}{Definition}
\newtheorem{pred}{Proposition}
\newtheorem{lem}{Lemma}
\newtheorem{teo}{Theorem}
\newtheorem{cor}{Corollary}
\newtheorem{cond}{Condition}
\newtheorem{conj}{Conjecture}

\begin{document}

%
% THIS IS COVER
%
% $Log: cover.tex,v $
% Revision 1.7  2001/02/08 18:53:16  boojum
% changed some \newpages to \cleardoublepages
%
% Revision 1.6  1999/10/21 14:49:31  boojum
% changed comment referring to documentstyle
%
% Revision 1.5  1999/10/21 14:39:04  boojum
% *** empty log message ***
%
% Revision 1.4  1997/04/18  17:54:10  othomas
% added page numbers on abstract and cover, and made 1 abstract
% page the default rather than 2.  (anne hunter tells me this
% is the new institute standard.)
%
% Revision 1.4  1997/04/18  17:54:10  othomas
% added page numbers on abstract and cover, and made 1 abstract
% page the default rather than 2.  (anne hunter tells me this
% is the new institute standard.)
%
% Revision 1.3  93/05/17  17:06:29  starflt
% Added acknowledgements section (suggested by tompalka)
%
% Revision 1.2  92/04/22  13:13:13  epeisach
% Fixes for 1991 course 6 requirements
% Phrase "and to grant others the right to do so" has been added to
% permission clause
% Second copy of abstract is not counted as separate pages so numbering works
% out
%
% Revision 1.1  92/04/22  13:08:20  epeisach
\title{A Proof of Tsygan's Formality Conjecture
for an Arbitrary Smooth Manifold}

\author{Vasiliy A. Dolgushev}
\department{Department of Mathematics}
% If the thesis is for two degrees simultaneously, list them both
% separated by \and like this:
% \degree{Doctor of Philosophy \and Master of Science}
\prevdegrees{Master of Science, Tomsk State University, 2001}
\degree{Doctor of Philosophy}
\degreemonth{June}
\degreeyear{2005}
\thesisdate{April 12, 2005}

%% By default, the thesis will be copyrighted to MIT.  If you need to copyright
%% the thesis to yourself, just specify the `vi' documentclass option.  If for
%% some reason you want to exactly specify the copyright notice text, you can
%% use the \copyrightnoticetext command.
%\copyrightnoticetext{\copyright IBM, 1990.  Do not open till Xmas.}

% If there is more than one supervisor, use the \supervisor command
% once for each.
\supervisor{Pavel Etingof}{Associate Professor of Mathematics, MIT}
\supervisor{Dmitry Tamarkin}{Assistant Professor of Mathematics,
Northwestern University}

% This is the department committee chairman, not the thesis committee
% chairman.  You should replace this with your Department's Committee
% Chairman.
\chairman{Pavel Etingof}{Chairman, Department Committee on Graduate Students}

% Make the titlepage based on the above information.  If you need
% something special and can't use the standard form, you can specify
% the exact text of the titlepage yourself.  Put it in a titlepage
% environment and leave blank lines where you want vertical space.
% The spaces will be adjusted to fill the entire page.  The dotted
% lines for the signatures are made with the \signature command.
\maketitle

% The abstractpage environment sets up everything on the page except
% the text itself.  The title and other header material are put at the
% top of the page, and the supervisors are listed at the bottom.  A
% new page is begun both before and after.  Of course, an abstract may
% be more than one page itself.  If you need more control over the
% format of the page, you can use the abstract environment, which puts
% the word "Abstract" at the beginning and single spaces its text.

%% You can either \input (*not* \include) your abstract file, or you can put
%% the text of the abstract directly between the \begin{abstractpage} and
%% \end{abstractpage} commands.

% First copy: start a new page, and save the page number.
\cleardoublepage
% Uncomment the next line if you do NOT want a page number on your
% abstract and acknowledgments pages.
% \pagestyle{empty}
\setcounter{savepage}{\thepage}
\begin{abstractpage}
Proofs of Tsygan's formality conjectures for chains
would unlock important algebraic tools which might
lead to new generalizations of the Atiyah-Patodi-Singer index
theorem and the Riemann-Roch-Hirzebruch theorem. Despite this
pivotal role in the traditional investigations and the efforts of
various people the most general version of Tsygan's formality
conjecture has not yet been proven. In my thesis I propose Fedosov
re\-so\-lu\-tions for the Hochschild cohomological and homological
complexes of the algebra of functions on an arbitrary smooth
manifold. Using these re\-so\-lu\-tions together with Kontsevich's
formality quasi-iso\-mor\-phism for Hochschild cochains
of ${\Bbb R}[[y^1, \dots, y^d]]$ and Shoikhet's formality
quasi-iso\-mor\-phism for Hochschild chains of ${\Bbb
R}[[y^1, \dots, y^d]]$ I prove Tsygan's formality conjecture for
Hochschild chains of the algebra of functions on an arbitrary
smooth manifold. The construction of the formality
quasi-isomorphism for Hochschild chains is manifestly functorial for isomorphisms
of the pairs $(M,\nabla)$, where $M$ is the manifold and $\nabla$ is
an affine connection on the tangent bundle. In my thesis I apply these results
to equivariant quantization, computation of Hochschild homology of quantum algebras and
description of traces in deformation quantization.
\end{abstractpage}

% Additional copy: start a new page, and reset the page number.  This way,
% the second copy of the abstract is not counted as separate pages.
% Uncomment the next 6 lines if you need two copies of the abstract
% page.
% \setcounter{page}{\thesavepage}
% \begin{abstractpage}
% \input{abstract}
% \end{abstractpage}

\cleardoublepage

\section*{Acknowledgments}
I would like to express my sincere thanks to
my thesis advisors Pavel Etingof and Dmitry Tamarkin
for their overall support, encouragement, and numerous
stimulating discussions.

I would like to thank V. Angeltveit,
R. Anno, A. Braverman, D. Calaque, A. Cattaneo, A. Chervov,
M. Ching, C. Douglas,
B. Enriquez, M. Fedorchuk, B. Feigin, G. Felder, J. Francis,
A. Gerasimov, A. Goncharov, V. Guillemin, G. Halbout,
R. Heluani, L. Hesselholt, M. Hill,
V. Ka$\rm\check{c}$, A. Khoroshkin,
M. Markl, R. Melrose,
A. Retakh, F. Rochon, Y. Rubinstein, L. Rybnikov,
P. Severa, A.A. Sharapov,
B. Shoikhet, D. Sullivan, D. Testa,
B. Tsygan, B. Vallette, D. Vera, A. Voronov,
Z. Zhang, and A.V. Zotov for helpful conversations.

I am also thankful to D. Calaque and G. Halbout for
the collaboration and to J. Francis for his valuable
criticisms concerning a preliminary version
of my thesis.

I acknowledge the hospitality of Institut de
Recherche Math\'ematique Avanc\'ee in Strasbourg and
Northwestern University's
Mathematics Department where part of my thesis
was done.

I would like to thank S. Kleiman, H. Miller, and V. Ostrik
for giving me a hard time during my studies.

I would like to thank our graduate administrator
L. Okun, and her assistant M. Gallarelli for their
help and patience.

I would like to thank H. Rogers for the unique
opportunity to tutor a research of undergraduate
students at MIT in the Summer Program for Undergraduate
Research. I had very interesting time working with
the students PoNing Chen, Teal Guidici, and Olga
Stroilova.

An essential role in my PhD studies has been played
by a warm and friendly atmosphere of
the Department of Mathematics of MIT.
In this respect I would like to mention
my friends I. Elson, J. Francis, R. Heluani,
Y. Rubinstein, and B. Santoro. To me
there is something in our relations
which will last for ever.

I would like to thank my friend Alexej Abyzov for
his moral support.

\vfill
\parbox[t]{\textwidth}{This research is supported by
the NSF grant DMS-9988796, the Grant
for Support of Scientific Schools NSh-1999.2003.2,
the grant INTAS 00-561 and the grant CRDF
RM1-2545-MO-03.}

%%%%%%%%%%%%%%%%%%%%%%%%%%%%%%%%%%%%%%%%%%%%%%%%%%%%%%%%%%%%%%%%%%%%%%
% END OF COVER
%

\newpage
~\\[7cm]
\begin{center}
{\Large \it To my brother Nickolay Dolgushev}
\end{center}

\pagestyle{plain}
\tableofcontents
\newpage
\listoffigures

%%%
%%
%% CHAPTER 1 (INTRODUCTION)
%% 1 MILE
%%
%%
%% 1/2 MILE
%%
%%

\chapter{Introduction}
Proofs of Tsygan's formality conjectures for chains
\cite{TT, TT1, Tsygan} would unlock important algebraic tools
which might lead to new generalizations of
the Atiyah-Patodi-Singer index theorem and
the Riemann-Roch-Hirzebruch theorem \cite{AS, BNT, Fedosov1, Losev, H, NT, NT1, TT}.
Despite this pivotal role in traditional
investigations and the efforts of various people
\cite{F-Sh, Sh, Sh1, TT, TT1} the most general
version of Tsygan's formality conjecture \cite{TT}
has resisted proof.

In my thesis I prove Tsygan's conjecture for Hochschild chains of the
algebra of functions on an arbitrary smooth manifold $M$ using the
globalization technique proposed in \cite{CFT} and \cite{CEFT} and
the formality quasi-isomorphism for Hochschild chains of
$\bbR[[y^1, \dots y^d]]$ constructed by Shoikhet \cite{Sh}.
This result allows me to prove Tsygan's conjecture \cite{Tsygan}
about Hochschild homology of the quantum algebra of functions on
an arbitrary Poisson manifold and to describe traces on this
algebra.

The most general version of Tsygan's formality
conjecture for chains says
that a pair of spaces of Hochschild cochains and Hochschild chains
of any associative algebra is endowed with the so-called
$Calc_{\infty}$-structure and if the algebra in question
is the algebra of functions on a smooth manifold then
the corresponding $Calc_{\infty}$-structure is
formal.
This statement was announced in \cite{Tsygan1} and
\cite{TT} but the proof has not yet been formulated.

In this context I would like to mention paper \cite{F-Sh}, in
which the authors prove a statement closely related to the cyclic
formality theorem. In particular, this assertion allows them to
prove a generalization of Connes-Flato-Sternheimer conjecture
\cite{CFS} in the Poisson framework.

The structure of my thesis is as follows. In the next chapter I
recall basic notions related to $\Linf$- or the so-called homotopy
Lie algebras. I introduce a notion of partial homotopy
between $\Linf$-morphisms and describe a useful technical tool that allows
me to utilize Maurer-Cartan elements of differential graded Lie
algebras (DGLA). In the third chapter I recall algebraic
structures on Hochschild complexes of associative algebra and
introduce the respective versions of these complexes for the
algebra of functions on a smooth manifold. In this section I
formulate the main result of my thesis (see theorem
\ref{thm-chain} on page \pageref{thm-chain})
and recall Kontsevich's and Shoikhet's formality
theorems for $\bbR^d$\,. The main part of this work concerns
the construction of Fedosov resolutions of the algebras of
polydifferential operators and polyvector fields, as well as the
modules of Hochschild chains and exterior forms. These resolutions
are constructed in chapter $4$. Using Fedosov's
resolutions in chapter $5$, I prove theorem \ref{thm-chain}.
In this chapter I also show that the Fedosov resolutions
provide me with a simple functorial construction
of Kontsevich's quasi-isomorphism from the DGLA
of polyvector fields to the DGLA of polydifferential
operators (see theorem \ref{Konets} on
page \pageref{Konets}).
At the end of chapter $5$ I apply
theorems \ref{thm-chain} and \ref{Konets}
to equivariant quantization, computation
of Hochschild homology of quantum algebras and
description of traces in deformation quantization.
In the concluding chapter I discuss recent works related to
generalizations and applications of the formality theorems
for Hochschild (co)chains.

My thesis is based on papers \cite{FTC, CEFT}.

~\\
{\bf Notation.} Throughout this work I assume the summation over repeated
indices. $M$ is a smooth real ma\-ni\-fold of dimension $d$\,.
The definition of antisymmetrization goes without
any auxiliary factors. Thus,
$$
v_1 \wedge v_2 = v_1 \otimes v_2 -
(-)^{|v_1||v_2|} v_2 \otimes v_1.
$$
I assume the Koszul rule of signs which says that
a transposition of any two vectors $v_1$ and
$v_2$ of degrees $k_1$ and $k_2$, respectively,
yields the sign
$$
(-1)^{k_1 k_2}\,.
$$
``DGLA'' always means a differential graded Lie algebra,
while ``DGA'' means a differential graded associative
algebra.
The arrow $\brarrow$ denotes an $\Linf$-morphism
of $\Linf$-algebras, the arrow $\bbrarrow$ denotes a morphism
of $\Linf$-modules, and the notation
$$
\begin{array}{c}
\cL\\[0.3cm]
\downarrow_{\,mod}\\[0.3cm]
\M
\end{array}
$$
means that $\M$ is an $\Linf$-module
over the $\Linf$-algebra $\cL$\,.
$S_n$ denotes the symmetric group of permutations
of $n$ elements and for natural numbers $k_1, \dots, k_q$,
$k_1 + \dots + k_q = n$
$Sh(k_1, \dots, k_q)\subset S_n$ is
the subset of $(k_1, \dots, k_q)$-shuffles.
Namely,
$$
Sh(k_1, \dots, k_q) =
$$
$$
\{\ve \in S_{n}\,|\,
\ve(1)<\ve(2)< \dots < \ve(k_1)\,,
\dots, \ve(n-k_q+1)<\ve(n-k_q+2)< \dots <
\ve(n)\}\,.
$$
I omit the symbol
$\wedge$ referring to a local basis of exterior forms, as if one
thought of $dx^i$'s as anti-commuting variables. The symbol
$\circ$ always stands for a composition of morphisms.
I denote by $\cxp(x)$ the following function
$$
\cxp(x)= e^x-1\,.
$$
Finally, I
denote by $\G(M, \cG)$  the vector space of
smooth sections of the bundle
$\cG$ and by  $\Omb(M, \cG)$ the vector space
of exterior forms with values in $\cG$\,.

%%
%% CHAPTER 1 ENDS
%%
%% CHAPTER 2
%%
%% 1 MILE
%%

\chapter{$L_{\infty}$-structures}
In this chapter I recall the notions of
$L_{\infty}$-algebras, $L_{\infty}$-morphisms,
$\Linf$-modules and morphisms
between $\Linf$-modules.
I introduce a notion of partial homotopy
between $\Linf$-morphisms and
describe an important technical tool, which
allows me to modify $\Linf$-structures
with the help of a Maurer-Cartan element.
A more detailed discussion of
this theory and its applications can be found in
papers \cite{Fuk,HS,LS}.

In this chapter all the
vector spaces, $L_{\infty}$-algebras, and
$\Linf$-modules are considered over a field
of characteristic zero.

\section{$\Linf$-algebras and $\Linf$-morphisms}
Let $\cL$ be a $\bbZ$-graded vector space
\begin{equation}
\label{mh}
\cL=\bigoplus_{k\in \bbZ} \cL^{k}\,.
\end{equation}
I assume that the direct sum in the right hand
side of (\ref{mh}) is bounded below.
To the space $\cL$ I assign a coassociative
cocommutative coalgebra (without counit) $C(\cL)$ cofreely
cogenerated by $\cL$ with a shifted parity.

The vector space of $C(\cL)$ is the exterior algebra of $\cL$
\begin{equation}
\label{C(L)}
C(\cL)= \bigwedge \cL\,,
\end{equation}
where the antisymmetrization is graded, that is
for any $\ga_1\in \cL^{k_1}$ and $\ga_2\in \cL^{k_2}$
$$
\ga_1 \wedge \ga_2 = -(-)^{k_1k_2} \ga_2 \wedge \ga_1\,.
$$

The comultiplication
\begin{equation}
\label{copro-C}
\D\, :\, C(\cL)\mapsto
C(\cL) \bigwedge C(\cL)
\end{equation}
is defined by the formulas  $(n>1)$
$$
\D(\ga_1)=0\,,
$$
\begin{equation}
\label{copro-eq}
\D (\ga_1\wedge \dots \wedge \ga_n)  =
\sum_{k=1}^{n-1} \sum_{\ve\in Sh(k,n-k)}
\pm \ga_{\ve(1)} \wedge \dots \wedge \ga_{\ve(k)}
\bigotimes
\ga_{\ve(k+1)} \wedge \dots \wedge \ga_{\ve(n)}\,,
\end{equation}
where $\ga_1, \dots,\, \ga_n$ are homogeneous elements
of $\cL$\,.

~\\
{\bf Remark.} I would like to mention that although I use
the Koszul rule the definition of the signs in (\ref{copro-eq})
is delicate. In fact one has to define $C(\cL)$ as the cofree
coalgebra of the suspended cooperad of cocommutative
coalgebras in the category of graded vector spaces.
To determine the correct signs in (\ref{copro-eq})
it is also helpful to use the fact that
the functor $\cL \mapsto C(\cL)$ should give
a cotriple. In the setting of commutative algebras
the reader can see the remark of E. Getzler
on p. 217 in \cite{GK}.

I now give the definition of $L_{\infty}$-algebra.
\begin{defi}
A graded vector space $\cL$
is said to be endowed with a structure of
an $L_{\infty}$-algebra if the cocommutative
coassociative coalgebra $C(\cL)$ cofreely cogenerated by
the vector space $\cL$ with a shifted parity is equipped with
a $2$-nilpotent coderivation $Q$ of degree $1$\,.
\end{defi}

To unfold this definition I first mention that
the kernel of $\D$ coincides with the subspace $\cL\subset
C(\cL)$.
\begin{equation}
\label{ker-D}
ker\D= \cL\,.
\end{equation}
Next, I recall that a map $Q$ is a
coderivation of $C(\cL)$ if and only if for
any $X\in C(\cL)$
\begin{equation}
\label{QD=DQ}
\D Q X = - (Q\otimes I \pm I \otimes Q) \D X\,.
\end{equation}
Substituting $X=\ga_1 \wedge \dots \wedge \ga_n$ in (\ref{QD=DQ}),
using (\ref{ker-D}), and performing the induction on $n$ I get that equation
(\ref{QD=DQ}) has the following general solution
$$
Q\, \ga_1\wedge \dots \wedge \ga_n  = Q_n (\ga_1, \dots, \ga_n)+
$$
\begin{equation}
\label{Qstruc}
\sum_{k=1}^{n-1} \sum_{\ve\in Sh(k, n-k)} \pm
Q_k(\ga_{\ve(1)},\dots, \ga_{\ve(k)})\wedge \ga_{\ve(k+1)}\wedge \dots\wedge
\ga_{\ve(n)}\,,
\end{equation}
where $\ga_1 \dots \ga_n$ are homogeneous elements of $\cL$
and $Q_n$ for $n\ge 1$ are arbitrary
polylinear antisymmetric graded maps
\begin{equation}
\label{stru-m-Q}
Q_n \,:\, \wedge^n \cL \mapsto \cL[2-n]\,, \qquad n\ge 1\,.
\end{equation}

It is not hard to see that $Q$ can be expressed inductively
in terms of the structure maps (\ref{stru-m-Q}) and
vice-versa.

Similarly, one can show that the nilpotency condition $Q^2=0$
is equivalent to a semi-infinite collection of quadratic
relations on (\ref{stru-m-Q}). The lowest of these relations
are
\begin{equation}
\label{Q-1is-diff}
(Q_1)^2 \ga=0\,, \qquad \forall~\ga\in \cL\,,
\end{equation}
\begin{equation}
\label{Q-1-Q-2}
Q_1(Q_2(\ga_1, \ga_2)) - Q_2(Q_1(\ga_1), \ga_2)
-(-)^{k_1} Q_2 (\ga_1, Q_1 (\ga_2))=0\,,
\end{equation}
and
$$
(-)^{k_1k_3}Q_2(Q_2(\ga_1, \ga_2), \ga_3)+ \,c.p.(1,2,3) =
$$
\begin{equation}
\label{Jacobi}
=Q_1 Q_3(\ga_1, \ga_2, \ga_3) + Q_3(Q_1\ga_1, \ga_2, \ga_3)
+(-)^{k_1} Q_3(\ga_1, Q_1\ga_2, \ga_3)
\end{equation}
$$
+(-)^{k_1+k_2} Q_3(\ga_1, \ga_2, Q_1\ga_3)\,,
$$
where $\ga_i\in \cL^{k_i}$\,.

Thus (\ref{Q-1is-diff}) says that $Q_1$ is a differential
in $\cL$\,, (\ref{Q-1-Q-2}) says that $Q_2$ satisfies Leibniz rule
with respect to $Q_1$\,, and (\ref{Jacobi}) implies that $Q_2$ satisfies
Jacobi identity up to $Q_1$-cohomologically trivial terms.

~\\
{\bf Example.} Any differential graded Lie algebra (DGLA)
$(\cL, \md, [,])$ is
an example of an $L_{\infty}$-algebra with the only two nonvanishing
structure maps
$$Q_1=\md\,, \qquad Q_2=[\,,\,]\,,$$
$$
Q_3 = Q_4 = Q_5 = \dots = 0\,.
$$

\begin{defi}
An $L_{\infty}$-morphism $F$ from the $L_{\infty}$-algebra
$(\cL, Q)$ to the $L_{\infty}$-algebra $(\cL^{\dia}, Q^{\dia})$ is
a homomorphism of the cocommutative coassociative coalgebras
$$
F\,:\,C(\cL) \mapsto C(\cL^{\dia})\,,
$$
\begin{equation}
\label{U=homo}
\D F(X)= F\otimes F (\D X)\,, \qquad X\in C(\cL)
\end{equation}
compatible with the coderivations $Q$ and
$Q^{\dia}$
\begin{equation}
\label{QF=FQ}
Q^{\dia} F (X)= F(Q X)\,, \qquad \forall~X\in C(\cL)\,.
\end{equation}
\end{defi}
In what follows the notation
$$
F\,:\, (\cL, Q) \brarrow
(\cL^{\dia}, Q^{\dia})
$$
means that $F$ is an $\Linf$-morphism from the
$\Linf$-algebra $(\cL, Q)$ to the $\Linf$-algebra
$(\cL^{\dia}, Q^{\dia})$\,.

The compatibility of the map (\ref{U=homo}) with coproducts
in $C(\cL)$ and  $C(\cL^{\dia})$ means that
$F$ is uniquely determined by the semi-infinite collection
of polylinear graded maps
\begin{equation}
\label{struct}
F_n : \wedge^n \cL \mapsto \cL^{\dia} [1-n], \qquad n\ge 1
\end{equation}
via the equations $(n\ge 1)$
\begin{equation}
\label{eq-forFn}
F(\ga_1\wedge \dots \wedge \ga_n)= F_n(\ga_1, \dots, \, \ga_n)+
\end{equation}
$$
\sum_{p> 1} \sum_{k_1,\dots , k_p\ge 1}^{k_1+\dots + k_p=n}
\sum_{\ve\in Sh(k_1, \dots, k_p) }
\pm F_{k_1}(\ga_{\ve (1)}, \dots, \, \ga_{\ve (k_1)})
\wedge \dots
$$
$$
\dots \wedge F_{k_p}(\ga_{\ve (n-k_p+1)},\dots,\,  \ga_{\ve (n)})\,,
$$
where $\ga_1, \dots, \, \ga_n$ are homogeneous elements
of $\cL$\,.

The compatibility of $F$ with coderivations (\ref{QF=FQ})
is a rather complicated condition for general $\Linf$-algebras.
However it is not hard to see that if (\ref{QF=FQ}) holds
then
$$
F_1(Q_1\ga)=Q^{\dia}_1 F_1(\ga)\,, \qquad\forall~ \ga\in \cL\,,
$$
that is the first structure map $F_1$ is always a
morphism of complexes $(\cL, Q_1)$ and $(\cL^{\dia}, Q^{\dia}_1)$\,.
This observation motivates the following
natural definition:
\begin{defi}
\label{qis}
A quasi-isomorphism $F$ from the $L_{\infty}$-algebra
$(\cL, Q)$ to the $L_{\infty}$-algebra $(\cL^{\dia}, Q^{\dia})$ is
an $L_{\infty}$-morphism from $\cL$ to $\cL^{\dia}$, the
first structure map $F_1$ of which induces a quasi-isomorphism
of complexes
\begin{equation}
\label{F1=qis}
F_1\,:\,(\cL, Q_1) \mapsto (\cL^{\dia}, Q^{\dia}_1)\,.
\end{equation}
\end{defi}

Let us suppose that
our $\Linf$-algebras $(\cL, Q)$ and $(\cL^{\dia}, Q^{\dia})$ are
just DGLAs
$(\cL, \md, [,])$ and $(\cL^{\dia}, \md^{\dia}, [ ,]^{\dia})$\,.
Then if $F$ is an $\Linf$-morphism from $\cL$
to $\cL^{\dia}$ the compatibility of $F$ with the
respective coderivations $Q$ and
$Q^{\dia}$ is equivalent to the following semi-infinite
collection of equations $(n\ge 1)$
$$
\md^{\dia} F_n(\ga_1, \ga_2, \ldots, \ga_n)-
\sum_{i=1}^n (-)^{k_1+\ldots+k_{i-1}+1-n}
F_n(\ga_1, \ldots, \md \ga_i, \ldots, \ga_n)=
$$
\begin{equation}
=\frac12 \sum_{k,l\ge 1}^{k+l=n}
 \sum_{\ve\in Sh(k,l)}
\pm [ F_k (\ga_{\ve_1}, \ldots, \ga_{\ve_k}), F_l (\ga_{\ve_{k+1}}, \ldots,
\ga_{\ve_{k+l}})]^{\dia}-
\label{q-iso}
\end{equation}
$$
-\sum_{i\neq j}
\pm F_{n-1}([\ga_i,\ga_j], \ga_1, \ldots, \hat{\ga_i}, \ldots, \hat{\ga_j}, \ldots \ga_n),
\qquad \ga_i \in \cL^{k_i}\,,
$$
where $\hat{\ga_i}$ means that the polyvector $\ga_i$ is missing.

~\\
{\bf Example.} An important example of a quasi-iso\-mor\-phism
from a DGLA $\cL$ to a DGLA $\cL^{\dia}$ is provided by a
DGLA-homomorphism
$$
\la \,:\,\cL \mapsto \cL^{\dia}\,,
$$
which induces an iso\-mor\-phism on the spaces
of cohomology
$H^{\bul}(\cL,\md)$ and $H^{\bul}(\cL^{\dia},\md^{\dia})$.
In this case the quasi-iso\-mor\-phism has the only
nonvanishing structure map $F_1$
$$
F_1=\la\,, \qquad F_2 = F_3 = \dots = 0\,.
$$

\section{$L_{\infty}$-modules and their morphisms}
Another important object of the ``$\Linf$-world'' I am going
to deal with is an $\Linf$-module over an $\Linf$-algebra.
Namely,
\begin{defi}
Let $\cL$ be an $L_{\infty}$-algebra. Then a graded
vector space $\M$ is endowed with a structure of an
$L_{\infty}$-module over $\cL$ if the cofreely cogenerated
comodule $C(\cL)\otimes \M$ over the coalgebra
$C(\cL)$ is endowed with a
$2$-nilpotent coderivation $\vf$ of degree $1$\,.
\end{defi}
To unfold the definition I first mention that the
total space of the comodule $C(\cL)\otimes \M$ is
\begin{equation}
\label{comodule}
C(\cL)\otimes \M=\bigwedge(\cL)\otimes \M\,,
\end{equation}
and the coaction
$$
\ma\,:\,C(\cL)\otimes \M \mapsto C(\cL)\bigotimes
(C(\cL)\otimes \M)
$$
is defined on homogeneous elements as follows
$$
\ma (\ga_1\wedge \dots  \wedge \ga_n \otimes v)=
$$
\begin{equation}
\label{coact}
\sum_{k=1}^{n-1} \sum_ {\ve\in Sh(k,n-k)}
 \pm\, \ga_{\ve (1)} \wedge \dots
\wedge \ga_{\ve (k)} \bigotimes \ga_{\ve(k+1)} \wedge  \dots \wedge
\ga_{\ve(n)}\otimes v
\end{equation}
$$
+\ga_1\wedge \dots  \wedge \ga_n \bigotimes v\,,
$$
where $\ga_1, \dots \ga_n$ are homogeneous elements of $\cL$, $v\in \M$\,,
$S_{n}$ is the group of permutations of $n$ elements and
the signs are determined with the help of the Koszul rule.
For example,
$$
\ma (v)=0\,, \qquad \forall~v\in \M\,,
$$
$$
\ma (\ga\otimes v)=\ga\bigotimes v\,,
\qquad \forall~v\in \M\,,~\ga\in \cL\,,
$$
and
$$
\ma (\ga_1\wedge \ga_2\otimes v)= \ga_1\wedge \ga_2\bigotimes v
+\ga_1 \bigotimes (\ga_2\otimes v) -(-)^{k_1 k_2} \ga_2 \bigotimes (\ga_1\otimes
v)
$$
for any $v\in \M$ and for any pair
$\ga_1\in \cL^{k_1}$\,, $\ga_2\in \cL^{k_2}$\,.

A direct computation shows that the coaction (\ref{coact})
satisfies the required axiom
$$
(I\otimes\ma)\ma (X) = (\D\otimes I)\ma (X)\,, \qquad \forall~X\in
C(\cL)\otimes \M\,,
$$
where $\D$ is the comultiplication (\ref{copro-C}) in
the coalgebra $C(\cL)$\,. It is also easy to see that
\begin{equation}
\label{ker-ma}
\ker \ma = \M\subset C(\cL)\otimes \M \,.
\end{equation}

By definition $\vf$ is a coderivation
of $C(\cL)\otimes \M$\,. This means that
for any $X\in C(\cL)\otimes \M$
\begin{equation}
\label{coder-mod}
\ma\,\vf X=  - Q\otimes I\, (\ma X) \pm I\otimes \vf \, (\ma X)\,,
\end{equation}
where $Q$ is the $\Linf$-algebra structure on $\cL$ (that
is a $2$-nilpotent coderivation of $C(\cL)$\,).

Substituting $X=\ga_1 \wedge \dots \wedge \ga_n$ in (\ref{coder-mod}),
using (\ref{ker-ma}), and performing the induction on $n$ I get that equation
(\ref{coder-mod}) has the following general solution
$$
\vf (\ga_1\wedge\dots \wedge\ga_n\otimes v)=
\vf_n (\ga_1, \dots, \ga_n, v)+
$$
\begin{equation}
\label{coder-mod-str}
\sum_{k=1}^{n-1}  \sum_{\ve\in Sh(k, n-k)} \pm\,
\ga_{\ve(1)}\wedge \dots \wedge \ga_{\ve(k)}\wedge
\vf_{n-k}(\ga_{\ve(k+1)}, \dots,\, \ga_{\ve(n)},v)
\end{equation}
$$
+(-)^{k_1 +\dots +k_n} \ga_1\wedge \dots \wedge \ga_n
\otimes\vf_0(v)+
$$
$$
\sum_{k=1}^{n} \sum_{\ve\in Sh(k,n-k)} \pm
Q_k(\ga_{\ve(1)},\dots, \ga_{\ve(k)})\otimes\ga_{\ve(k+1)}\wedge \dots\wedge
\ga_{\ve(n)}\otimes v\,,
$$
where $\ga_i\in \cL^{k_i}$, $v\in \M$, $Q_k$'s
represent the $L_{\infty}$-algebra structure on
$\cL$ and $\{\vf_n\}$ for $n\ge 0$ are arbitrary
polylinear antisymmetric graded maps
\begin{equation}
\label{coll-vf}
\vf_n \,:\, \wedge^n\cL\otimes \M \mapsto \M[1-n]\,.
\end{equation}
Equation (\ref{coder-mod-str}) allows us to express $\vf$ inductively
in terms of its structure maps (\ref{coll-vf}) and vice-versa.

Similarly, one can show that the nilpotency condition $\vf^2=0$
is equivalent to the following semi-infinite collection of quadratic
relations in $\vf_k$ and $Q_l$
$(n\ge 0)$
$$
\vf_0 (\vf_n(\ga_1,\dots,\, \ga_n,v))
-(-)^{1-n}\vf_n(Q_1(\ga_1),\dots,\, \ga_n,v)- \dots
$$
$$
\ldots -(-)^{k_1 + \dots + k_{n-1} + 1-n} \vf_n(\ga_1,\dots,\, Q_1(\ga_n),v)-
$$
$$
(-)^{k_1 + \dots + k_{n}+1-n} \vf_n(\ga_1,\dots,\,\ga_n,\vf_0(v))=
$$
\begin{equation}
\label{nil-fi}
\frac12\sum_{k=1}^{n-1} \sum_{\ve\in Sh(k, n-k)}
\pm \vf_k(\ga_{\ve(1)},\dots,\, \ga_{\ve(k)},
\vf_{n-k}(\ga_{\ve(k+1)}, \dots,\, \ga_{\ve(n)},v))+
\end{equation}
$$
\frac12\sum_{k=1}^{n-1}\sum_{\ve\in Sh(k+1,n-k-1)}
\pm \vf_{n-k}(Q_{k+1}(\ga_{\ve(1)},\dots,\, \ga_{\ve(k+1)}),
\ga_{\ve(k+2)}, \dots,\, \ga_{\ve(n)},v)\,,
$$
$$
\ga_i \in \cL^{k_i},\qquad v\in\M\,.
$$
The signs in the above equations are defined
with the help of the Koszul rule.

For $n=0$ equation (\ref{nil-fi}) says that $\vf_0$
is a differential on $\M$
$$
(\vf_0)^2=0\,,
$$
and for $n=1$ it
says that $\vf_1$ is closed with
respect to the natural differential acting on
the vector space $Hom(\cL\otimes \M, \M)$
$$
\vf_0\vf_1 (\ga, v)-\vf_1(Q_1\ga,v) -
(-)^{k}\vf_1(\ga,\vf_0(v))=0\,,
$$
$$
\forall~ \ga\in \cL^k\,,\qquad v\in\M\,.
$$
For an $\Linf$-module structure I reserve the following
notation
$$
\begin{array}{c}
\cL\\[0.3cm]
\phantom{aa}\downarrow_{\,mod}^{\vf}\\[0.3cm]
(\M, \vf_0)
\end{array}
$$
where $\cL$ stands for the $\Linf$-algebra and
$\M$ stands for the respective graded vector space.

~\\
{\bf Example.} The simplest example of an $\Linf$-module
is a DG module $(\M,\bb)$ over a DGLA $(\cL, \md,[\,,\,])$.
In this case the only nonvanishing structure maps of $\vf$
are
$$
\vf_0(v)=\bb(v)\,, \qquad v\in \M\,,
$$
and
$$
\vf_1(\ga,v)=\rho(\ga)\, v\,, \qquad \ga\in\cL\,,
~v\in \M\,,
$$
where $\rho$ is the action of $\cL$ on $\M$.
The axioms of DGLA module
$$
\bb^2=0\,,
$$
$$
\bb (\rho(\ga)\,v)=\rho(\md \ga)\,v+(-)^{k}\rho(\ga)\,\bb(v)\,,
\qquad \ga\in \cL^{k}\,,
$$
$$
\rho(\ga_1)\rho(\ga_2)\,v- (-)^{k_1k_2} \rho(\ga_2)\rho(\ga_1)\,v=
\rho([\ga_1,\ga_2])v\,,
$$
$$
\ga_1\in \cL^{k_1}\,, \qquad \ga_2\in \cL^{k_2}
$$
are exactly the axioms of $\Linf$-module.

\begin{defi}
Let $\cL$ be an $\Linf$-algebra and $(\M, \vf^{\M})$\,,
$(\N,\vf^{\N})$ be
$\Linf$-modules over $\cL$. Then a morphism $\ka$ from the
comodule $C(\cL)\otimes \M$ to the comodule
$C(\cL)\otimes \N$ compatible with the coderivations
$\vf^{\M}$ and $\vf^{\N}$
\begin{equation}
\label{Kvf=vfK}
\ka(\vf^{\M} X)=\vf^{\N}(\ka X)\,, \qquad
\forall~X\in C(\cL)\otimes \M
\end{equation}
is called an morphism between $\Linf$-modules
$(\M, \vf^{\M})$ and $(\N,\vf^{\N})$\,.
\end{defi}
Unfolding this definition one can easily show that
the morphism $\ka$ is uniquely determined by its
structure maps
\begin{equation}
\label{struct-ka}
\ka_n\,:\,\wedge^n \cL \otimes \M \mapsto \N[-n]\,, \qquad n\ge 0
\end{equation}
via the following equations
$$
\ka(\ga_1\wedge\dots \wedge\ga_n\otimes v)=
\ka_n(\ga_1, \dots,\,\ga_n,v)+
$$
\begin{equation}
\label{struct-ka-ka}
\sum_{k=1}^{n-1}
\sum_{\ve \in Sh(k,n-k)} \pm\, \ga_{\ve(1)}\wedge \dots \wedge \ga_{\ve(k)}
\otimes \ka_{n-k}(\ga_{\ve(k+1)},\dots,\,\ga_{\ve(n)},v)
\end{equation}
$$
+\ga_1\wedge \dots \wedge \ga_n \otimes \ka_0(v)\,.
$$

Relation (\ref{Kvf=vfK}) is equivalent to the
following semi-infinite collection of equations
$(n\ge 0)$
$$
\vf^{\N}_0 \ka_n(\ga_1,\dots,\, \ga_n,v) -
(-)^n\ka_n(Q_1\ga_1,\ga_2, \dots,\, \ga_n)- \dots
$$
$$
\ldots -(-)^{k_1+\dots + k_n+n} \ka_n(\ga_1, \dots,\, \ga_n, \vf^{\M}_0 v)=
$$
$$
\sum_{p=0}^{n-1}
\sum_{\ve \in Sh(p,n-p)} \pm
\ka_p(\ga_{\ve(1)}, \dots, \, \ga_{\ve(p)},
\vf^{\M}_{n-p}(\ga_{\ve(p+1)}, \dots,\,
\ga_{\ve(n)},v))
$$
\begin{equation}
\label{ka-vf=vf-ka}
-\sum_{p=1}^n
\sum_{\ve \in Sh(p,n-p)} \pm
\vf^{\N}_p(\ga_{\ve(1)}, \dots, \, \ga_{\ve(p)},\ka_{n-p}(\ga_{\ve(p+1)},
 \dots,\, \ga_{\ve(n)}),v)
\end{equation}
$$
+\sum_{p=2}^{n}
\sum_{\ve \in Sh(p,n-p)} \pm
\ka_{n-p+1}(Q_p(\ga_{\ve(1)}, \dots,\,\ga_{\ve(p)}),\ga_{\ve(p+1)},
\dots, \,\ga_{\ve(n)},v)\,,
$$
$$
\ga_i\in\cL^{k_i}\,,\qquad v\in \M\,.
$$
It is not hard to check that an ordinary morphism of
DG modules over an ordinary DGLA provides us with the
simplest example of the morphism between $\Linf$-modules.

For $n=0$ equation (\ref{ka-vf=vf-ka}) reduces to
$$
\ka_0(\vf^{\M}_0 v)=\vf^{\N}_0 \ka_0(v)\,, \qquad
v\in \M
$$
and hence the zero-th structure map of $\ka$ is
always a morphism of complexes $(\M, \vf^{\M}_0)$ and
$(\N, \vf^{\N}_0)$\,. This motivates the following
\begin{defi}
A quasi-isomorphism $\ka$ of $\Linf$-modules
$(\M, \vf^{\M})$ and $(\N, \vf^{\N})$ is a morphism
between these $\Linf$-modules with the zero-th structure
map $\ka_0$ being a quasi-isomorphism of complexes
 $(\M, \vf^{\M}_0)$ and
$(\N, \vf^{\N}_0)$\,.
\end{defi}

In what follows the notation
$$
\M \stackrel{\ka}{\bbrarrow} \N
$$
means that $\ka$ is a morphism from the
$\Linf$-module $\M$ to the $\Linf$-module $\N$\,.

To this end I mention that there is
another definition of an $\Linf$-module
over an $\Linf$-algebra which is known \cite{Tsygan} to be
equivalent to the definition I gave above.
\begin{defi}
Let $\cL$ be an $\Linf$-algebra. Then a complex $(\M, \bb)$ is
called an $\Linf$-module over $\cL$ if there is an
$\Linf$-morphism $\eta$ from $\cL$ to $Hom(\M,\M)$\,, where
$Hom(\M,\M)$ is naturally viewed as a DGLA with the
differential induced by $\bb$\,.
\end{defi}
The structure maps $\vf_n$ of the respective
coderivation of the comodule $C(\cL)\otimes \M$
are related to $\bb$ and the structure maps of
the $\Linf$-morphism $\eta$ in the following simple way
\begin{equation}
\label{vf-b-chi}
\bb=\vf_0\,,\qquad
\eta_n(\ga_1, \dots,\, \ga_n)(v)=
\vf_n(\ga_1,\dots,\,\ga_n,v)\quad (n\ge 1)\,,
\end{equation}
$$
\ga_i\in \cL,~ v\in \M\,.
$$

%
%%
%%  Section on homotopies in 1/2 mile.
%%
%%
%

\section{Partial homotopies between $\Linf$-morphisms}
\label{homotopy}
In this section I introduce a notion of
partial homotopy between two
$\Linf$-morphisms. I will
use this notion in
section \ref{revisited}.

Let $(\cL, Q)$ and $(\cL^{\dia}, Q^{\dia})$ be
two $\Linf$-algebras. As above,
I denote by $Q$ and $Q^{\dia}$ the
corresponding codifferentials of
the cocommutative coassociative
coalgebras $C(\cL)$ and $C(\cL^{\dia})$\,.
Let
$$
F : C(\cL) \mapsto C(\cL^{\dia})
$$
be an $\Linf$-morphism from $\cL$ to
$\cL^{\dia}$.

One can observe that if a map
$$
H : C(\cL) \mapsto C(\cL^{\dia})
$$
is of degree $-1$ then the map
$$
\tF = F+ Q^{\dia}\, H + H\, Q
\,:\, C(\cL) \mapsto C(\cL^{\dia})
$$
is of degree zero, and moreover
it is compatible with the coderivations
$Q$ and $Q^{\dia}$
\begin{equation}
\label{Q-tF}
Q^{\dia} \tF = \tF Q\,.
\end{equation}

A compatibility of $\tF$ with the coproducts (\ref{copro-eq})
in $C(\cL)$ and $C(\cL^{\dia})$ is equivalent
to a rather complicated equation for the map $H$
$$
\D \, H \, Q  - (Q^{\dia}\otimes I \pm
I\otimes Q^{\dia}) \D\, H =
$$
\begin{equation}
\label{copro-H}
F \otimes (Q^{\dia} H + H Q) +(Q^{\dia} H + H Q) \otimes F
+ (Q^{\dia} H + H Q)\otimes (Q^{\dia} H + H Q)\,.
\end{equation}
However, if $H$ satisfies the following
equation
\begin{equation}
\label{copro-H1}
\D H = - \Big (
F \otimes H + H \otimes F  +
\frac1{2}(H \otimes Q^{\dia} H  + Q^{\dia} H \otimes H)
+
\frac1{2}(H Q \otimes H  + H \otimes H Q)
\Big )
\end{equation}
then due to (\ref{QD=DQ}) and
(\ref{U=homo}) $H$ satisfies
(\ref{copro-H}) as well.

Using (\ref{ker-D}) it is not
hard to get the most general solution
of equation (\ref{copro-H1}).
Namely, any solution $H$ of (\ref{copro-H1})
is uniquely determined by a semi-infinite
collection of polylinear graded maps
\begin{equation}
\label{H-struct}
H_n : \wedge^{n}\cL \mapsto \cL^{\dia}[-n]\,,
\qquad
H_n = p r \circ H\, \Big|_{\wedge^{n}\cL}\,,
\end{equation}
where $p r$ is the canonical projection
\begin{equation}
\label{pr-cL1}
p r : \bigwedge(\cL^{\dia}) \to \cL^{\dia}\,.
\end{equation}
In order to restore the map $H$ from
the collection (\ref{H-struct}) one solves
(\ref{copro-H1}) iteratively from
$\wedge^{< n}\cL$ to $\wedge^n\cL$
starting with
\begin{equation}
\label{H-1}
H(\ga) = H_1 (\ga )\,, \qquad \forall~ \ga \in \cL\,.
\end{equation}
I refer to (\ref{H-struct}) as
structure maps of $H$\,.

It is immediate from (\ref{eq-forFn}) and
(\ref{H-1}) that for any $\ga \in \cL$
\begin{equation}
\label{F1-H1}
\tF_1(\ga)= F_1(\ga) + Q^{\dia}_1 \, H_1(\ga)
+ H_1\, Q_1(\ga)\,,
\end{equation}
where $F_1$ and $\tF_1$ are the first
structure maps of $F$ and $\tF$, respectively.
This observation motivates the following definition
\begin{defi}
A map
$$
H : C(\cL) \mapsto C(\cL^{\dia})[-1]
$$
is called a partial homotopy
between $\Linf$-morphisms
$$
F,\, \tF : (C(\cL),Q) \mapsto (C(\cL^{\dia}),Q^{\dia})
$$
if it satisfies (\ref{copro-H1}) and
$$
\tF = F+ Q^{\dia}\, H + H\, Q\,.
$$
Two $\Linf$-morphisms are called
partially homotopic if they are
connected by a finite chain
of partial homotopies.
\end{defi}
{\bf Remark 1.} It is easy to see that
equation (\ref{copro-H1}) still holds
if I replace $H$ by $-H$ and $F$ by
$F+ Q^{\dia} H + H Q$. However a composition
of two partial homotopies is not
in general a partial homotopy\footnote{I am thankful
to G. Felder for this observation.}.
That is why I extend the relation of
partial homotopy to an equivalence relation
by transitivity.

~\\
{\bf Remark 2.} The correct notion
of homotopy between $\Linf$-morphisms
is based on the structure of
the closed model category on the category
of $\Linf$-algebras \cite{Hinich}, \cite{Quillen}.
Unfortunately, I do not know how to
relate the above notion of the partial
homotopy to the correct notion of homotopy
based on the closed model
category structure. For my purposes
the above {\it ad hoc} notion will be
sufficient.

Let me prove the following auxiliary
statement:
\begin{lem}
\label{styag}
Let
$$
F : C(\cL) \mapsto C(\cL^{\dia})
$$
be a quasi-isomorphism from
an $\Linf$-algebra $(\cL, Q)$ to
an $\Linf$-algebra $(\cL^{\dia}, Q^{\dia})$.
For $n\ge 1$ and any map
$$
\tH : \wedge^{n}\cL \mapsto \cL^{\dia}[-n]
$$
one can construct
a quasi-isomorphism
$$
\tF : C(\cL) \mapsto C(\cL^{\dia})
$$
such that for any $m<n$
\begin{equation}
\label{kak}
\tF_m = F_m \, :\, \wedge^m \cL\mapsto \cL^{\dia}
\end{equation}
and
$$
\tF_n (\ga_1, \dots, \ga_n)= F_n (\ga_1, \dots, \ga_n) +
$$
\begin{equation}
\label{kak1}
Q^{\dia}_1 \tH(\ga_1, \dots, \ga_n) -
(-)^{n} \tH (Q_1(\ga_1), \ga_2, \dots, \ga_n) - \dots
\end{equation}
$$
\dots -(-)^{n+k_1+ \dots + k_{n-1}}
\tH(\ga_1, \dots, \ga_{n-1}, Q_1(\ga_n))\,,
$$
where $\ga_i \in \cL^{k_i}$\,.
\end{lem}
{\bf Proof.} It is obvious that if
a partial homotopy $H$ has
the following structure maps:
$$
H_m =
\begin{cases}
\begin{array}{cc}
\tH& {\rm if}\quad  m=n\,,\\
0 & {\rm otherwise}
\end{array}
\end{cases}
$$
then $\tF= F + Q^{\dia} H + H Q $
satisfies the desired properties
(\ref{kak}), (\ref{kak1}).
Since $F$ is a quasi-isomorphism
so is $\tF$. $\Box$

~\\
{\bf REMARK.} From now on {\bf all} $\Linf$-algebras
are DGLAs. ``Weird'' things I still borrow
from the ``$\Linf$-world'' are $\Linf$-morphisms,
$\Linf$-modules, and morphisms between such modules.

%%
%%
%% MAURER-CARTAN IN 1/2 MILE!!!
%%
%%

\section{Maurer-Cartan elements and twisting procedures}
\label{section-MC}
Motivated by deformation theory I consider
DGLAs $\cL$ equipped with a
complete descending
filtration
\begin{equation}
\label{Filt-cL}
\cL = \cF^0 \cL \supset \cF^1 \cL \supset \dots\,,
\qquad
\cL = \lim_{n} \cL/\cF^n\cL\,.
\end{equation}
In this section I assume that
all DGLAs and $\Linf$-modules are
equipped with complete descending filtrations
and all $\Linf$-morphisms as well as
morphisms of $\Linf$-modules
are compatible with
these filtrations. Furthermore,
I require that all quasi-isomorphisms
of the corresponding complexes
are {\it strongly compatible} with
the filtrations. Namely,
\begin{cond}
\label{condition}
Let $\la$ be a quasi-isomorphism
$$
\la : \cL^{\bul} \mapsto \tcL^{\bul}
$$
of filtered complexes
$\cL^{\bul}$, $\tcL^{\bul}$\,.
I say that $\la$ is
compatible with the
filtrations if for any
filtration subcomplex
$\cF^k \cL^{\bul}\subset \cL^{\bul}$
$$
\la\,  \Big|_{\cF^k \cL^{\bul}} \, : \,
 \cF^k \cL^{\bul} \mapsto \cF^k \tcL^{\bul}
$$
is a quasi-isomorphism.
\end{cond}
I will assume this compatibility condition
throughout my thesis.

If $\cL$ is such a filtered DGLA then $\cF^1\cL^0$ is
a projective limit of nilpotent Lie algebras. Therefore,
$\cF^1\cL^0$ can be ``integrated'' to a prounipotent group.
I denote this group by $\mG(\cL)$\,.

Let me recall the following definition:
\begin{defi}
\label{MCdefi}
Let $(\cL, \md, [,])$ be a filtered DGLA.
Then $\pi\in \cF^1\cL^1$
is called a Maurer-Cartan element if
\begin{equation}
\label{MC}
\md \pi+ \frac{1}{2}[\pi,\pi]=0\,.
\end{equation}
\end{defi}

The Lie algebra $\cF^1\cL^0$ acts naturally on
the cone (\ref{MC}) of Maurer-Cartan elements
\begin{equation}
\label{act-onMC}
\rho(\xi) \pi = \md \xi + [\pi, \xi]\,,
\qquad  \xi \in \cF^1\cL^0\,,
\end{equation}
and the action (\ref{act-onMC}) obviously
lifts to the action of the corresponding prounipotent
group $\mG(\cL)$. The quotient space of the cone (\ref{MC})
with respect to the $\mG(\cL)$-action is called the
{\it moduli space} of the DGLA $\cL$.

It turns out that a
quasi-iso\-mor\-phism (see definition \ref{qis}) between DGLAs
gives a bijective correspondence between their moduli spaces.
A weaker version of this statement is proved in proposition
\ref{twist-morph} (see claim $4$). This version says that
if $F$ is an $\Linf$-morphism from a DGLA $(\cL, \md,[,])$
to a DGLA $(\cL^{\dia}, \md^{\dia},[,]^{\dia})$ and
$\pi\in \cF^1\cL^1$ is a Maurer-Cartan element of
$\cL$ then
\begin{equation}
\label{MCelem}
S =\sum_{n\ge 1}\frac{1}{n!} F_n(\pi, \ldots, \pi)
\end{equation}
is a Maurer-Cartan element of $\cL^{\dia}$.

Notice that the infinite sum in (\ref{MCelem}) is well-defined
because $\cL^{\dia}$ is assumed to be complete with
the respect to the corresponding filtration.
All elements of this sum are of degree $1$ since
for any $n$ $F_n$ shifts the degree by $1-n$ (see
(\ref{struct}))\,.

Using a Maurer-Cartan element $\pi\in \cF^1\cL^1 $ one can
naturally modify the structure of the DGLA on $\cL$
by adding the inner derivation $[\pi,\cdot \,]$ to the
initial differential $\md$. Thanks to Maurer-Cartan
equation (\ref{MC}) this new derivation $\md+[\pi,\cdot \,]$
is $2$-nilpotent. This modification can
be described in terms of the respective $\Linf$-structure .
Namely, the coderivation $Q^{\pi}$ on the coassociative cocommutative coalgebra
$C(\cL)$ corresponding to the new DGLA structure
$(\cL, \md+[\pi,\cdot \,], [\,,\,])$ is
related to the initial coderivation $Q$ by the equation
\begin{equation}
\label{twistQ}
Q^{\pi}(X)=\exp((-\pi)\wedge) Q(\exp(\pi\wedge) X)\,, \qquad
X \in C(\cL)\,,
\end{equation}
where the sum
$$
\exp(\pi\wedge) \underbrace{~} =  \underbrace{~} +
\pi\wedge  \underbrace{~} +\frac{1}{2!} \pi\wedge \pi\wedge  \underbrace{~}
+ \dots
$$
is well-defined since $\pi\in \cF^1\cL^1$\,.

I call this procedure of changing the initial DGLA structure
on $ \cL$ {\it twisting} of the DGLA $\cL$ by the Maurer-Cartan
element $\pi$\,.  This terminology is borrowed from
Quillen's paper \cite{Q} (see App. B $5.3$). This twisting
procedure is also extensively used in paper \cite{Ye}
by A. Yekutieli on deformation quantization
in algebraic geometry setting.

Similar twisting procedures by a Maurer-Cartan element
can be defined for an $\Linf$-morphism,
for an $\Linf$-module, and for a morphism
of $\Linf$-modules. In the following propositions
I describe these procedures.
\begin{pred}[See also theorem 0.1 in  \cite{Ye1}]
\label{twist-morph}
If $F$ is an $\Linf$-morphism
$$
F\,:\, (\cL, Q) \brarrow (\cL^{\dia}, Q^{\dia})
$$
of DGLAs, $\pi\in \cF^1 \cL^1 $
and an element $S\in \cF^1(\cL^{\dia})^1$ is
given by equation (\ref{MCelem}) then
\begin{enumerate}
\item For any homogeneous element $X\in C(\cL)$
\begin{equation}
\label{Del-MC}
\D (\exp(\pi\wedge) X)=
\exp(\pi\wedge)\bigotimes \exp(\pi\wedge)(\D X)+
\cxp(\pi)\bigotimes \exp(\pi\wedge) X -
\end{equation}
$$
(-)^{|X|} \exp(\pi\wedge) X \bigotimes \cxp(\pi)\,,
$$
where
\begin{equation}
\label{cxp}
\cxp(\pi)=\sum_{k=1}^{\infty}\frac1{k!}\underbrace{\pi \wedge\dots \wedge
\pi}_{k}\,.
\end{equation}
\item Equation (\ref{MC}) is equivalent to
\begin{equation}
\label{Q-MC}
Q(\cxp(\pi))=0\,.
\end{equation}
\item
\begin{equation}
\label{F-MC}
F(\cxp(\pi))= \cxp(S)\,.
\end{equation}

\item If $\pi$ a Maurer-Cartan element then
so is $S$ and the map
\begin{equation}
\label{twist-F}
F^{\pi}=\exp(-S\wedge) F \exp(\pi \wedge)\,:\,
C(\cL)\mapsto C(\cL^{\dia})
\end{equation}
defines an $\Linf$-morphism between the
DGLAs $\cL^{\pi}$ and $\cL^{\dia\,S}$,
obtained via twisting by the Maurer-Cartan elements
$\pi$ and $S$, respectively.

\item Let $\pi$ be a Maurer-Cartan element.
If $F$ is a quasi-isomorphism satisfying condition
\ref{condition} (on page \pageref{condition} )
then so is $F^{\pi}$\,.

\end{enumerate}

\end{pred}
In what follows I refer to $F^{\pi}$ in (\ref{twist-F})
as an $\Linf$-mor\-phism (or a quasi-iso\-mor\-phism)
{\it twisted by the Maurer-Cartan element} $\pi$\,.
It is not hard to see that the structure maps
of the twisted $\Linf$-morphism $F^{\pi}$ are given by
\begin{equation}
\label{twist-F-str}
F^{\pi}_n(\ga_1, \dots,\,\ga_n)=
\sum_{k=0}^{\infty} \frac1{k!}
F_{n+k} (\pi,\dots,\, \pi, \ga_1, \dots, \, \ga_n)\,,
\qquad \ga_i\in \cL\,.
\end{equation}
~\\[0.3cm]
{\bf Proof.} In order to prove claim $1$ I
introduce an auxiliary variable $t$ and analyze
a slightly stronger statement
\begin{equation}
\label{Del-MC1}
\D (\exp(t \pi\wedge) X) \stackrel{?}{=}
\exp(t\pi\wedge)\bigotimes \exp(t\pi\wedge)(\D X)+
\cxp(t\pi)\bigotimes \exp(t\pi\wedge) X -
\end{equation}
$$
(-)^{|X|} \exp(t\pi\wedge) X \bigotimes \cxp(t\pi)\,.
$$
It is clear that (\ref{Del-MC1}) holds for $t=0$.
On the other hand a direct computation shows
that both sides of (\ref{Del-MC1}) satisfies
the following differential equation:
$$
\frac{d}{dt} W(t) =
\Big(\pi \wedge \bigotimes 1  + 1 \bigotimes\pi \wedge \Big) W(t)
+ \pi \bigotimes \exp(t\pi\wedge)X - (-)^{|X|}
\exp(t\pi\wedge)X \bigotimes \pi\,.
$$
Thus equation (\ref{Del-MC1}) holds and claim $1$
follows.

It is obvious that (\ref{Q-MC}) implies
(\ref{MC}). Let me prove the converse statement.
First, I observe that if $(\cL, \md, [,])$ is a DGLA
then the collection
$$
Q^t_1 = t\md, \quad
Q^t_2= [\,,\,],\quad
Q^t_3 = Q^t_4 = \dots =0
$$
defines a DGLA  on $\cL[t]$\,.
Second, $t\pi$ is a Maurer-Cartan
element in $(\cL[t], t \md ,[,])$ and
the equation
\begin{equation}
\label{Q-MC1}
Q^t(\cxp(t\pi))\stackrel{?}{=}0
\end{equation}
obviously holds for $t=0$.
Using the Maurer-Cartan equation
(\ref{MC}) it is not hard to prove that
the left hand side $Z(t)= Q^t(\cxp(t\pi))$
of (\ref{Q-MC1}) satisfies the following
differential equation:
$$
\frac{d}{d t} Z(t) = Z(t)\wedge \pi\,.
$$
Since this equation is homogeneous
(\ref{Q-MC1}) holds for any $t$ and
claim $2$ follows.

To prove that the element
$$
Y= F(\cxp(\pi)) - \cxp (S) \in C(\cL^{\dia})
$$
is vanishing I observe that for
any $\pi \in \cL^1$
\begin{equation}
\label{D-cxp}
\D (\cxp (\pi)) = \cxp(\pi) \bigotimes \cxp(\pi)\,.
\end{equation}
Furthermore, due to (\ref{MCelem})
$Y$ lies in the kernel of
the natural projection
\begin{equation}
\label{pr-cL}
p r : \bigwedge(\cL^{\dia}) \to \cL^{\dia}\,.
\end{equation}

Let us prove by induction that
\begin{equation}
\label{Fm-all}
Y\in \cF^m (C(\cL^{\dia}))
\end{equation}
for all $m$.

By definition of the Maurer-Cartan
element $\pi \in \cF^1\cL$. Therefore
the element $S$ (\ref{MCelem}) belongs
to $\cF^1\cL^{\dia}$ and hence
$$
Y\in \cF^1 (C(\cL^{\dia}))\,.
$$
Let me take it as base of the induction
and suppose that (\ref{Fm-all}) holds
for some $m$.

Equation (\ref{D-cxp}) and the compatibility
of the map $F$ with the coproducts (\ref{copro-eq})
in $C(\cL)$ and $C(\cL^{\dia})$ implies that
$$
\D Y \in \cF^{m+1}(\wedge^2 C(\cL^{\dia}))\,.
$$
Therefore due to (\ref{ker-D}) the image of $Y$ in
$\cF^{m}(C(\cL^{\dia}))/\cF^{m+1}(C(\cL^{\dia}))$
belongs to
$$\cF^m \cL^{\dia} / \cF^{m+1}\cL^{\dia}\,.$$
But the image of $Y$ vanishes under the
projection $p r$ (\ref{pr-cL}). Hence,
$$
Y\in \cF^{m+1} (C(\cL^{\dia}))\,,
$$
and therefore (\ref{Fm-all}) holds
for all $m$. Since $\cL^{\dia}$ is complete
with respect to the filtration equation
(\ref{F-MC}) is proved.

Let me now turn to claim $4$.
While the compatibility of $F^{\pi}$
with the coderivations $Q^{\pi}$ and $Q^{\dia\,S}$ follows directly from
the definitions the compatibility
of $F^{\pi}$ with the coproducts in $C(\cL)$ and
$C(\cL^{\dia})$ requires some work.
Using claim $1$ and $3$
I get that for any homogeneous $X\in C(\cL)$
$$
\D \exp(-S\wedge) F \exp(\pi \wedge)X =
\exp(-S\wedge)\bigotimes \exp(-S\wedge)(F\bigotimes F)
\D\exp(\pi \wedge)X +
$$
$$
\cxp(-S)\bigotimes \exp(-S\wedge)  F \exp(\pi \wedge) X
-(-)^{|X|}
\exp(-S\wedge) F\exp(\pi \wedge)X \bigotimes \cxp(-S)=
$$
$$
\cxp(-S)\bigotimes F^{\pi}X +
F^{\pi}X \bigotimes \cxp(-S)+
$$
$$
(F^{\pi}\bigotimes F^{\pi}) (\D X)+
$$
$$
\exp(-S\wedge)\bigotimes \exp(-S\wedge)(F\bigotimes F)
(\cxp(\pi)\bigotimes \exp(\pi \wedge)  X)-
$$
$$
(-)^{|X|}\exp(-S\wedge)\bigotimes \exp(-S\wedge)(F\bigotimes F)
(\exp(\pi \wedge) X \bigotimes \cxp(\pi))\,.
$$
The first and the second terms in the latter expression
cancel with the forth and the
fifth terms, respectively, due to
claim $3$ and the following obvious identity
between Taylor series
\begin{equation}
\label{Taylor}
e^{-S}\cxp (S)= -\cxp(-S)\,.
\end{equation}
Thus, I get the desired relation
$$
\D F^{\pi}(X)= (F^{\pi}\bigotimes F^{\pi}) (\D X)\,.
$$

Claim $5$ is proved by the standard argument of the
spectral sequence. We have to prove that
the first structure map $F^{\pi}_1$ of the
twisted $\Linf$-morphism (\ref{twist-F}) is a quasi-isomorphism
from the complex $(\cL, \md + [\pi,\cdot ])$ to
the complex $(\cL^{\dia}, \md^{\dia} + [S,\cdot ]^{\dia} )$.
These complexes are filtered and $F^{\pi}_1$
is compatible with the filtration. Since $F_1$ is
a quasi-isomorphism between the complexes
$(\cL, \md)$ and $(\cL^{\dia}, \md^{\dia})$
and $F_1$ satisfies condition \ref{condition}
(on page \pageref{condition})
the map $F^{\pi}_1$ induces a quasi-isomorphism
on the zeroth level of the corresponding
spectral sequences. Therefore $F^{\pi}_1$ gives
a quasi-isomorphism on the terminal $E_{\infty}$-level.
Hence, due to the standard snake-lemma argument of
homological algebra $F^{\pi}_1$ is also a
quasi-isomorphism.

Proposition \ref{twist-morph} is proved.  $\Box$

\begin{pred}
\label{twist-mod}
If $(\cL, \md, [,])$ is a DGLA, $(\M, \vf)$ is an $\Linf$-module over
$\cL$ and $\pi\in \cF^1\cL^1$ is a Maurer-Cartan element then
\begin{enumerate}

\item For any\footnote{if $X=v\in \M $ I set
``$\pi\wedge X=\pi \otimes X$''}
$X\in C(\cL)\otimes \M$
\begin{equation}
\label{Coact-MC}
\ma (\exp(\pi\wedge) X)=
\exp(\pi\wedge)\bigotimes \exp(\pi\wedge)(\ma X)+
\cxp(\pi)\bigotimes  \exp(\pi\wedge) X\,,
\end{equation}
where $\ma$ is the coaction (\ref{coact}) and
$\cxp(\pi)$ is defined in the previous proposition.

\item The following map
\begin{equation}
\label{twist-vf}
\vf^{\pi}=\exp(-\pi\wedge)\vf \exp(\pi\wedge)\,:\,
C(\cL )\otimes \M  \mapsto
C(\cL )\otimes \M
\end{equation}
is a $2$-nilpotent coderivation of the comodule
$C(\cL )\otimes \M $\,.

\item  If
$\tilde{\vf}\,:\, \cL \brarrow (Hom(\M,\M),\vf_0)$
is the $\Linf$-morphism induced by the module
structure $\vf$ then the twisted $\Linf$-morphism
$\tilde{\vf}^{\pi}$ defines the $\Linf$-module
structure given in (\ref{twist-vf})\,.

\item If $\ka\,:\,\M \bbrarrow \N$ is a morphism
of $\Linf$-modules $(\M, \vf)$ and
$(\N, \psi)$ over $\cL$ then the map
\begin{equation}
\label{twist-kappa}
\ka^{\pi}=\exp(-\pi\wedge)\ka \exp(\pi\wedge)\,:\,
C(\cL)\otimes \M \mapsto
C(\cL)\otimes \N
\end{equation}
is a morphism between $\Linf$-modules
$(\M , \vf^{\pi})$ and
$(\N , \psi^{\pi})$ over
$(\cL , \md+[\pi,\cdot \,], [,])$

\item If $\ka$ is a quasi-isomorphism of
$\Linf$-modules $\M$ and $\N$ and $\ka$
satisfies condition \ref{condition}
(on page \pageref{condition})
then $\ka^{\pi}$ (\ref{twist-kappa}) is also
a quasi-isomorphism.

\end{enumerate}
\end{pred}
In what follows I refer to $\vf^{\pi}$ in (\ref{twist-vf})
and $\ka^{\pi}$ in (\ref{twist-kappa}), respectively,
as an $\Linf$-module structure and a morphism of
$\Linf$-modules twisted by the
Maurer-Cartan element $\pi$\,.
It is not hard to see that the structure maps
of the twisted coderivation $\vf^{\pi}$ and
the twisted morphism $\ka^{\pi}$ are given by
\begin{equation}
\label{twist-vf-str}
\vf^{\pi}_n(\ga_1, \dots,\,\ga_n,v)=
\sum_{m=0}^{\infty} \frac1{m!}
\vf_{n+m} (\pi,\dots,\, \pi, \ga_1, \dots, \, \ga_n,v)\,,
\end{equation}
\begin{equation}
\label{twist-ka-str}
\ka^{\pi}_n(\ga_1, \dots,\,\ga_n,v)=
\sum_{m=0}^{\infty} \frac1{m!}
\ka_{n+m} (\pi,\dots,\, \pi, \ga_1, \dots, \, \ga_n,v)\,,
\end{equation}
where
$$
\ga_i\in \cL\,,~v\in \M\,.
$$

~\\[0.3cm]
{\bf Proof.} Claim $1$ is proved with the help
of the similar scale trick ($\pi \to t \pi$) I used
in the proof of the previous proposition.
Claim $2$ follows from claim $1$ of
this proposition and claim $2$ of the previous
proposition. Claim $4$ essentially follows from
claim $1$ of this proposition and claim $3$ is proved by
comparing the corresponding structure maps.

Claim $5$ is proved by the standard argument of the
spectral sequence. We have to prove that
the zeroth structure map $\ka^{\pi}_0$ of the
twisted morphism (\ref{twist-kappa}) is a quasi-isomorphism
from the complex $(\M, \vf^{\pi}_0)$ to
the complex $(\N, \psi^{\pi}_0 )$.
These complexes are filtered and $\ka^{\pi}_0$
is compatible with the filtration. Since $\ka_0$ is
a quasi-isomorphism between the complexes
$(\M, \vf_0)$ and $(\N, \psi_0)$ and $\ka_0$
satisfies condition \ref{condition}
(on page \pageref{condition})
the map $\ka^{\pi}_0$ induces a quasi-isomorphism
on the zeroth level of the corresponding
spectral sequences. Therefore $\ka^{\pi}_0$ gives
a quasi-isomorphism on the terminal $E_{\infty}$-level.
Hence, due to the standard snake-lemma argument of
homological algebra $\ka^{\pi}_0$ is also a
quasi-isomorphism. $\Box$

From the definitions of the above twisting procedures,
it is not hard to see that these procedures are functorial.
Namely,
\begin{pred}
\label{functor}
If $F\,:\,\cL\brarrow\cL^{\dia}$ and
$F^{\dia}\,:\,\cL^{\dia}\brarrow\cL^{\club}$ are
$\Linf$-morphisms of DGLAs, $\pi$
is a Maurer-Cartan element of $\cL$
and $S$ is the corresponding Maurer-Cartan element
(\ref{MCelem}) of $\cL^{\dia}$ then
$$
(F^{\dia}\circ F)^{\pi}= F^{\dia\, S}\circ F^{\pi}\,,
$$
where $\circ$ stands for the composition of $\Linf$-morphisms.
Furthermore, the twisting procedure assigns to
any Maurer-Cartan element of a DGLA $\cL$
a functor from the category of $\Linf$-modules
to itself. $\Box$
\end{pred}

Let us turn to the moduli functor of
Maurer-Cartan elements and prove that this
functor provides us with a homotopy
invariant of a DGLA.
\begin{pred}[K. Fukaya, \cite{Fuk}, theorem 2.2.2]
\label{Fukaya}
Let $(\cL,\md, [,])$ and
$(\cL^{\dia}, \md^{\dia}, [,]^{\dia})$ be
two completely filtered DGLAs
and let $F$ be a quasi-isomorphism (see definition
\ref{qis}) from $\cL$ to
$\cL^{\dia}$ compatible with the filtrations
in the sense of condition \ref{condition}.
Then (\ref{MCelem}) gives a bijective correspondence
between the moduli spaces of $\cL$ and $\cL^{\dia}$.
\end{pred}
{\bf Remark.} The case of the ordinary
(not $\Linf$) quasi-isomorphism is treated
by Goldman and Millson \cite{Gold, Gold1}.
Its generalization to $\Linf$ setting has been
a folklore\footnote{I learnt this statement from B. Shoikhet.}
and was quoted by several authors
(without proofs). In principle, using the
``nonsense'' of the homotopy theory \cite{Hinich},
\cite{Quillen}
it is possible to reduce the statement of
the above proposition to the result of
Goldman and Millson \cite{Gold, Gold1}.
In \cite{Fuk} K. Fukaya gives a rigorous proof of this statement
both in $\Linf$ and $A_{\infty}$ settings. However,
since his proof is based on other results which appear
elsewhere, I decided to give my own proof.

~\\
{\bf Proof.} First I have to prove that (\ref{MCelem})
gives a well defined map from the moduli space
of $\cL$ to the moduli space of $\cL^{\dia}$.
Due to claim $4$ of proposition \ref{twist-morph}
it suffices to check that the map (\ref{MCelem})
of cones of Maurer-Cartan elements is compatible
with the action (\ref{act-onMC}) of
$\cF^1\cL^0$ and $\cF^1(\cL^{\dia})^0$, respectively.

If $\pi$ is a Maurer-Cartan
element of $\cL$ and $\xi\in \cF^1\cL^0$ then
$$
\rho(\xi)(\cxp(\pi)) =\exp (\pi \wedge) Q^{\pi}(\xi)\,,
$$
where $Q^{\pi}$ the DGLA structure on $\cL$ twisted
by the Maurer-Cartan element $\pi$.
Hence, due to claim $3$ of proposition \ref{twist-morph}
$$
\rho(\xi) (\cxp(S_{\pi})) = \rho(\xi)F \cxp(\pi)=
F\exp (\pi \wedge) Q^{\pi}(\xi)\,,
$$
where
$$
S_{\pi} = \sum_{k=1}^{\infty}
\frac1{k!} F_k (\pi, \dots, \pi)\,.
$$
Or equivalently,
$$
\rho(\xi) (\cxp(S_{\pi})) =
Q^{\dia}  F(\exp(\pi\wedge)\xi)=
\exp(S\wedge) (Q^{\dia})^{S_{\pi}}(F^{\pi}(\xi))
\,,
$$
where $F^{\pi}$ is the twisted quasi-isomorphism.
Thus,
$$
\rho(\xi) S_{\pi} = \rho(F^{\pi}(\xi)) S_{\pi}
$$
and hence (\ref{MCelem}) gives a well-defined map
\begin{equation}
\label{MC-corr}
F_{MC} : M C(\cL) \mapsto M C(\cL^{\dia})
\end{equation}
from the moduli space $M C(\cL)$ of
the DGLA $\cL$ to the moduli space $M C(\cL^{\dia})$ of
the DGLA $\cL^{\dia}$.

Let $S \in \cF^1 (\cL^{\dia})^1$ be
a Maurer-Cartan element of $\cL^{\dia}$. I denote by
$\mG$ the prounipotent group
corresponding to the Lie algebra $\cF^1\cL^0$ and
by $\mG[S]$ the $\mG$-orbit that passes through
$S$. To prove surjectivity of the map (\ref{MC-corr})
I show by induction that there exists a collection
of pairs $(S_m, \pi_m)$, $m\ge 1$ where $S_m$ are
Maurer-Cartan
elements of $\cL^{\dia}$ belonging to the orbit
$\mG[S] $,  $\pi_m\in \cF^1 \cL^1$,
\begin{equation}
\label{eq:conv}
S_{m+1} - S_m \in \cF^m \cL^{\dia}\,,
\qquad \pi_{m+1} - \pi_m \in  \cF^m \cL\,,
\end{equation}
\begin{equation}
\label{pi-MC-m}
\md \pi_m +\frac1{2} [\pi_m,\pi_m] \in \cF^m \cL\,,
\end{equation}
and
\begin{equation}
\label{pi-S-m}
S_m - \sum_{k=1}^{\infty}
\frac1{k!} F_k (\pi_m, \dots, \pi_m)
 \in \cF^m \cL^{\dia}\,.
\end{equation}

For $m=1$ I set $S_1=S$, $\pi_1=0$. Then equations
(\ref{pi-MC-m}) and (\ref{pi-S-m}) obviously hold.
Let me take it as a base of the induction and
assume that (\ref{eq:conv}), (\ref{pi-MC-m}),
and (\ref{pi-S-m}) hold up to $m$ but $S_{m+1}$ and
$\pi_{m+1}$ are not chosen. It suffices to prove
that there exists a pair $(S_{m+1}, \pi_{m+1})$ such that
$S_{m+1}\in \mG[S_m]$, (\ref{eq:conv}) is satisfied
and equations (\ref{pi-MC-m}), (\ref{pi-S-m}) hold
for $m$ replaced by $m+1$\,.

Due to assumption
(\ref{pi-S-m}) and the Maurer-Cartan equation
$\displaystyle \md^{\dia} S_m+ \frac1{2}[S_m,S_m]^{\dia}$
I get that
\begin{equation}
\label{S-pi}
\md^{\dia} (S_m-S_{\pi_m}) +
(\md^{\dia}S_{\pi_m} + \frac1{2}[S_{\pi_m}, S_{\pi_m}]^{\dia} )
\in  \cF^{m+1}\cL^{\dia}\,,
\end{equation}
where I denoted by $S_{\pi_m}$ the sum
$$
S_{\pi_m} = \sum_{k=1}^{\infty}
\frac1{k!} F_k (\pi_m, \dots, \pi_m)\,.
$$
Hence due to claim $2$ of proposition \ref{twist-morph}
$$
\md^{\dia} (S_m-S_{\pi_m}) + Q^{\dia} \cxp(S_{\pi_m})
\in  \cF^{m+1}(C(\cL^{\dia}))\,.
$$
Therefore using claim $3$ of proposition
\ref{twist-morph} one gets
$$
\md^{\dia} (S_m-S_{\pi_m}) +  F Q \cxp(\pi_m)
\in \cF^{m+1}(C(\cL^{\dia}))\,.
$$
Applying assumption (\ref{pi-MC-m}) and
claim $2$ of proposition \ref{twist-morph}
once again I get
$$
\md^{\dia} (S_m-S_{\pi_m}) +
F_1 (\md \pi_m + \frac1{2}[\pi_m , \pi_m])
\in \cF^{m+1}\cL^{\dia}\,.
$$
Since $F_1$ is a quasi-isomorphism
of complexes $(\cL, \md)$ and
$(\cL^{\dia}, \md^{\dia})$ compatible
with the filtrations in the sense of
condition \ref{condition}
(on page \pageref{condition})
there
exist an element $\pia \in \cF^m \cL^1$
and an element $\xi\in \cF^m (\cL^{\dia})^0$
such that
$$
\md \pia + \md \pi_m + \frac1{2}[\pi_m, \pi_m] \in \cF^{m+1} \cL\,,
$$
and
$$
S_m - S_{\pi_m} + F_1(\pia) + \md^{\dia} \xi \in \cF^{m+1} \cL^{\dia}\,.
$$
Thus, if I set $S_{m+1}= \exp (\rho(\xi)) S_m$ and
$\pi_{m+1} = \pi_m + \pia$ then
$S_{m+1}$ and $\pi_{m+1}$ satisfy condition
(\ref{eq:conv}) and, moreover, equations
(\ref{pi-MC-m}), (\ref{pi-S-m}) hold
with $m$ replaced by $m+1$.

Since the DGLAs $\cL$ and $\cL^{\dia}$ are
complete with respect to the
filtrations the surjectivity
of the map (\ref{MC-corr})
follows from the existence of the
desired collection $(S_m,\pi_m)$.

The injectivity is proved by
analyzing the differential of
the map (\ref{MC-corr}).

Indeed, let $\pi$ be a Maurer-Cartan
element of $\cL$. Then the tangent space
to the cone (\ref{MC}) is cut in
$\cL^1$ by the equation
\begin{equation}
\label{tang-pi}
\md \pi^t + [\pi, \pi^t]=0\,, \qquad \pi^t \in \cL^1\,.
\end{equation}

Therefore by definition of the action (\ref{act-onMC})
of $\cF^1 \cL^0$ on the cone (\ref{MC}) the tangent space of
the moduli space $M C(\cL)$ of Maurer-Cartan elements
to the orbit $[\pi]$ passing through $\pi$ is
the first cohomology group of the
complex $(\cL, \md + [\pi,\cdot ])$
$$
T_{[\pi]} (M C(\cL)) = H^1 (\cL, \md + [\pi,\cdot ])\,.
$$

By the assumption of the proposition
$F$ is a quasi-isomorphism between $\cL$ and
$\cL^{\dia}$. Hence, due to claim $5$ $F^{\pi}$
is a quasi-isomorphism of the twisted
DGLAs
$(\cL, \md + [\pi,\cdot\, ], [,])$ and
$(\cL^{\dia}, \md^{\dia} + [S_{\pi},\cdot\, ]^{\dia}, [,])$.
Therefore the differential of
the map (\ref{MC-corr}) is
an isomorphism. Thus, it
is injective and the proposition follows. $\Box$

%%
%% CHAPTER 2 ENDS
%%
%% CHAPTER 3
%%
%% 1 MILE
%%

\chapter{Mosaic}
In this chapter I recall the basic
algebraic structures on Hochschild (co)chains.
I formulate the main result of this thesis, the formality
theorem for Hochschild chains
of the algebra of functions on a smooth
manifold, and state Kontsevich's and Shoikhet's
formality theorems for Hochschild (co)chains
of the algebra of functions on $\bbR^d$.

\section{Algebraic structures on Hochschild (co)chains}
For a unital associative algebra $\bbA$ (over a field of characteristic zero)
I denote by $C^{\bul}(\bbA)$ the vector space of Hochschild cochains
with a shifted grading
\begin{equation}
\label{H-coch}
C^{n} (\bbA)=Hom(\bbA^{\otimes (n+1)}, \bbA)\,,~(n\ge 0)\,, \qquad
C^{-1}(\bbA)=\bbA\,.
\end{equation}
The space $C^{\bul}(\bbA)$ can be endowed with the
Gerstenhaber bracket \cite{G}, defined
between homogeneous elements $P_1\in C^{k_1}(\bbA)$ and
$P_2\in C^{k_2}(\bbA)$ as follows
\begin{equation}
\label{Gerst}
[P_1, P_2]_G = P_1 \bul P_2
-(-)^{k_1k_2} P_2 \bul P_1\,,
\end{equation}
where
$$ (P_1 \bul P_2)(a_0,\,\dots, a_{k_1+k_2})=$$
\begin{equation}
\label{bullet}
\sum_{i=0}^{k_1}(-)^{i k_2}
P_1(a_0,\,\dots , P_2 (a_i,\,\dots,a_{i+k_2}),\, \dots,
a_{k_1+k_2})\,.
\end{equation}
Direct computation shows that (\ref{Gerst}) is a Lie bracket and
therefore $C^{\bul}(\bbA)$ is a graded Lie algebra.

For the same unital algebra $\bbA$, I denote by $C_{\bul}(\bbA)$
the vector space of Hochschild chains
\begin{equation}
\label{H-chain}
C_{n} (\bbA)= \bbA\otimes \bbA^{\otimes n}\,,~(n\ge 1),
 \qquad
C_{0}(\bbA)=\bbA\,.
\end{equation}
The space $C_{\bul}(\bbA)$ can be endowed with the structure
of a graded module over the Lie algebra $C^{\bul}(\bbA)$ of Hochschild
cochains. For homogeneous elements the action of $C^{\bul}(\bbA)$ on
$C_{\bul}(\bbA)$ is defined as follows:
$$
R : C^k(\bbA) \otimes C_n(\bbA) \to C_{n-k}(\bbA)\,, \quad
P \otimes (a_0\otimes a_1\otimes \dots \otimes a_n)
\mapsto R_P(a_0\otimes a_1\otimes \dots \otimes a_n)
$$
\begin{equation}
\label{cochain-act}
R_{P}(a_0\otimes a_1\otimes \dots \otimes a_n)=
\sum_{i=0}^{n-k}(-)^{ki} a_0 \otimes \dots \otimes
P(a_i, \dots, a_{i+k})\otimes \dots \otimes a_n+
\end{equation}
$$
\sum_{j=n-k}^{n-1}(-)^{n(j+1)}P(a_{j+1}, \dots, a_n, a_0, \dots, a_{k+j-n})
\otimes a_{k+j+1-n}\otimes \dots \otimes a_j\,, \qquad a_i\in \bbA\,.
$$
The proof of the required axiom of the Lie algebra module
\begin{equation}
\label{R-OK}
R_{[P_1, P_2]_G}= R_{P_1} R_{P_2} -
(-)^{|P_1||P_2|} R_{P_2} R_{P_1}
\end{equation}
can be found in paper \cite{Ezra}, in which
it is discussed in a more general
$A_{\infty}$ setting (see lemma $2.3$ in \cite{Ezra}).

The multiplication $\mu_0$ in the algebra $\bbA$ can be
naturally viewed as an element of $C^1(\bbA)$ and
the associativity condition for $\mu_0$ can be rewritten
in terms of bracket (\ref{Gerst}) as
\begin{equation}
\label{assoc}
[\mu_0,\mu_0]_G=0\,.
\end{equation}

Thus, on the one hand $\mu_0$ defines a
$2$-nilpotent interior derivation
of the graded Lie algebra $C^{\bul}(\bbA)$
\begin{equation}
\label{pa}
\pa P = [\mu_0,P]_G : C^k(\bbA)\mapsto C^{k+1}(\bbA)\,, \qquad
\pa^2 =0\,,
\end{equation}
and on the
other hand $\mu_0$ endows the graded vector space $C_{\bul}(\bbA)$
with the differential
\begin{equation}
\label{b-chain}
\mb  = R_{\mu_0} : C_k(\bbA)\mapsto C_{k-1}(\bbA)\,, \qquad
\mb^2 =0\,.
\end{equation}
Equation (\ref{R-OK}) implies that
$$
R_{\pa P}= \mb R_{P}- (-)^{k} R_{P}\mb\,,
\qquad P\in C^k(\bbA)
$$
and therefore the vector spaces $C^{\bul}(\bbA)$ and $C_{\bul}(\bbA)$
become a pair of a DGLA and a DG module over
this DGLA.

~\\
{\bf Remark 1.}
Notice that the differentials (\ref{pa}) and (\ref{b-chain})
are exactly the Hochschild coboundary and boundary operators
on $C^{\bul}(\bbA)$ and $C_{\bul}(\bbA)$, respectively.
Thus, the Hochschild (co)homology groups $HH^{\bul}(\bbA)$
and $HH_{\bul}(\bbA)$ form a pair of graded Lie
algebra and a graded module of this Lie algebra.

~\\
{\bf Remark 2.} Notice that the action $R$
(\ref{cochain-act}) is not compatible with the grading on
$C_{\bul}(\bbA)$ and the differential
(\ref{b-chain}) have degree $-1$ (not $+1$).
In order to get the DGLA module in the
sense of the previous chapter one has to use
the converted grading on $C_{\bul}(\bbA)$.
However, I prefer to restrain the conventional grading on
the space of Hochschild chains keeping in mind
the above remark.

Let me also recall that the graded vector space
$C^{\bul-1}(\bbA)$ is endowed with the obvious
associative product
$$
\cup : C^{k_1-1}(\bbA)\otimes C^{k_2-1}(\bbA)
\mapsto C^{k_1+ k_2 -1}(\bbA)\,,
$$
\begin{equation}
\label{cup}
P_1 \cup P_2 (a_1, \dots, a_{k_1+k_2})=
P_1 (a_1, \dots, a_{k_1}) \cdot
P_2 (a_{k_1+1}, \dots, a_{k_1+k_2})\,,
\end{equation}
$$
P_i\in C^{k_i-1}(\bbA)\,, \qquad a_j\in \bbA\,,
$$
where $\cdot$ denotes the product in
the algebra $\bbA$\,.

The product (\ref{cup}) is compatible with
the Hochschild differential (\ref{pa})
in the sense of the following
equation
\begin{equation}
\label{cup-pa}
\pa (P_1 \cup P_2) =
P_1 \cup \pa(P_2) + (-)^{k_2} \pa (P_1)
\cup P_2\,, \qquad P_2\in C^{k_2-1}(\bbA)\,.
\end{equation}
Thus Hochschild chains also form a DGA.

I will refer to the product (\ref{cup})
as {\it the cup-product} and I will use it in
the proof of proposition \ref{vot-ono}.

\section{Formality theorems}
I will be mainly interested in the algebra
$\bbA_0=C^{\infty}(M)$
where $M$ is a smooth ma\-ni\-fold of dimension $d$.
A natural analogue of the complex of Hochschild cochains for
this algebra is the complex $\Dp$ of
polydifferential operators with the same
differential as in $C^{\bul}(C^{\infty}(M))$
\begin{equation}
\label{D-poly}
\Dp=\bigoplus_{k=-1}^{\infty} D_{poly}^k(M)\,,
\qquad D^{-1}_{poly}(M)= C^{\infty}(M)\,,
\end{equation}
where $D_{poly}^k(M)$ consists of po\-ly\-dif\-fe\-ren\-tial
operators of rank $k+1$
$$
P\,:\,C^{\infty}(M)^{\otimes (k+1)} \mapsto C^{\infty}(M)\,.
$$
Similarly, instead of the complex $C_{\bul}(C^{\infty}(M))$ I consider
three versions of the vector space $\Cp$ of Hochschild
chains for $C^{\infty}(M)$
\begin{enumerate}
\item

\begin{equation}
\label{H-chains1}
C_{function}^{poly}(M)= \bigoplus_{n\ge 0} C^{\infty}(M^{n+1})\,,
\end{equation}

\item

\begin{equation}
\label{H-chains2}
C_{germ}^{poly}(M)= \bigoplus_{n\ge 0}
germs_{\D(M^{n+1})}C^{\infty}(M^{n+1})\,,
\end{equation}

\item

\begin{equation}
\label{H-chains3}
C_{jet}^{poly}(M)= \bigoplus_{n\ge 0}
jets^{\infty}_{\D(M^{n+1})}C^{\infty}(M^{n+1})\,,
\end{equation}
\end{enumerate}
where $\D(M^{n+1})$ is the diagonal in $M^{n+1}$\,.

It is not hard to see that
the Gerstenhaber bracket (\ref{Gerst}), the action
(\ref{cochain-act}), the differentials (\ref{pa}),
(\ref{b-chain}), and the cup-product (\ref{cup})
still make sense if I replace $C^{\bul}(C^{\infty}(M))$ by
$\Dp$ and $C_{\bul}(C^{\infty}(M))$ by either of versions
(\ref{H-chains1}), (\ref{H-chains2}),
(\ref{H-chains3}) of $\Cp$.
Thus, $\Dp$ and $\Cp$
are DGLA and a DG module over this DGLA, respectively,
and, moreover, $D^{\bul-1}_{poly}(M)$ is a DGA.
I use the same notations for all the
operations $[,]_G$, $R_{P}$, $\pa$, $\mb$,
and $\cup$ when I speak of $\Dp$ and $\Cp$\,.

The cohomology of $\Dp$ and of $\Cp$ is
described by Hochschild-Kostant-Rosenberg type theorems.
The original version of the theorem
\cite{HKR} by Hoch\-schild, Kos\-tant, and Ro\-sen\-berg
says that the module of Hochschild
homology of a smooth affine algebra is isomorphic to the
module of exterior differential forms on the respective
affine algebraic variety. A dual version of this theorem
was proved in \cite{Ye1} (see corollary $4.12$).
In the $C^{\infty}$ setting we have
\begin{pred} [J. Vey, \cite{Vey}] \label{HKR-Vey}
Let
\begin{equation}
\label{T}
\Tp=\bigoplus_{k=-1}^{\infty} T^k_{poly}(M)\,, \qquad
T^k_{poly}(M)=\G(M, \wedge^{k+1} TM)
\end{equation}
be a vector space of the polyvector fields on
$M$ with shifted grading. If $\Tp$ is regarded as a
complex with a vanishing differential then
the natural map
\begin{equation}
\label{U-1}
\cV(\ga)(a_0, \dots, a_k)= i_{\ga}(d a_0\wedge \dots\wedge d a_k)
\,:\,T^k_{poly}(M)\mapsto D^k_{poly}(M)\,, \qquad k\ge -1
\end{equation}
defines a quasi-isomorphism of complexes
$(\Tp, 0)$ and $(\Dp, \pa)$\,. Here
$d$ stands for the De Rham differential and
$i_{\ga}$ denotes the contraction of the polyvector
field $\ga$ with an exterior form.
\end{pred}

The most general $C^{\infty}$-manifold version of the
Hoch\-schild-Kos\-tant-Ro\-sen\-berg theorem is
due to N. Teleman \cite{Tel}\footnote{See also \cite{Connes},
in which this statement was proven for any compact
smooth manifold.}
\begin{pred} [Teleman, \cite{Tel}] \label{HKR-Connes}
Let
\begin{equation}
\label{exter-alg}
\cA^{\bul}(M)=\bigoplus_{k\ge 0} \cA^{k}(M)\,, \qquad
\cA^k(M)=\G(M, \wedge^{k} T^{\ast}M)
\end{equation}
be a vector space of the exterior forms on $M$.
If $\cA^{\bul}(M)$ is regarded as a
complex with a vanishing differential then
the natural map
\begin{equation}
\label{C-1}
\mC(a_0\otimes \dots \otimes a_k)= a_0 d a_1\wedge \dots\wedge d a_k
\,:\,C^{poly}_{k}(M)\mapsto A^{k}(M)\,, \qquad k\ge 0
\end{equation}
defines a quasi-isomorphism of complexes
$(\Cp, \mb)$ and $(\cA^{\bul}(M), 0)$
for either of versions (\ref{H-chains1}),
(\ref{H-chains2}), (\ref{H-chains3}) of $\Cp$.
\end{pred}

One can easily check that the Lie algebra structure
induced on cohomology
$$H^{\bul}(\Dp, \pa)=\Tp$$
coincides with the one given by the so-called
Schouten-Nijenhuis bracket
$$
[,]_{SN}\,:\,\Tp\bigwedge \Tp\mapsto \Tp\,.
$$
This bracket
is defined as an ordinary Lie bracket between vector fields
and then extended by Leibniz rule
\begin{equation}
\label{eq:Leib}
[\ga_1,\ga_2\wedge \ga_3]_{SN}=[\ga_1,\ga_2]_{SN} \wedge
\ga_3 + (-1)^{|\ga_1|(|\ga_2|-1)}
\ga_2 \wedge [\ga_1,\ga_3]_{SN}\,, \quad
\ga_i \in T^{\bul}_{poly}(M)
\end{equation}
with respect to the $\wedge$-product
to an arbitrary pair of po\-ly\-vec\-tor
fields.

Furthermore, the DGLA $\Dp$-module structure on
$\Cp$ induces a $\Tp$-module
structure on the vector space $\AM$
which coincides with the one
defined by the action of a polyvector field on
exterior forms via the Lie derivative
\begin{equation}
\label{Tpoly-act}
L_\ga = d \, i_{\ga} + (-)^k i_{\ga} \, d\,,
\qquad \ga\in T^k_{poly}(M)\,,
\end{equation}
where as above $d$ stands for the De Rham differential
and $i_{\ga}$ denotes the contraction of the
polyvector field $\ga$ with an exterior form.

~\\
{\bf Remark.} In what follows
I will restrict myself to the third version
(\ref{H-chains3}) of $\Cp$ and
since all $\Dp$-modules
(\ref{H-chains1}), (\ref{H-chains2}),
(\ref{H-chains3}) are naturally
quasi-isomorphic the further results
will hold for versions (\ref{H-chains1}), (\ref{H-chains2})
as well.

For my purposes it will be very convenient
to represent the chains (\ref{H-chains3})
as $\OM$-linear homomorphisms from
$\Dp$ to $\OM$. Namely, one can equivalently
define ($k\ge 0$)
\begin{equation}
\label{jets}
C^{poly}_{k}(M) = Hom_{\OM} (D_{poly}^{k-1}(M), \OM)\,.
\end{equation}

To avoid the shift in the above formula
let me introduce the auxiliary graded
bundle of polyjets
placed in non-negative degrees
\begin{equation}
\label{polyjets}
J_{\bul} = \bigoplus_{k\ge 0} J_k\,,
\qquad
J_k  = Hom_{\cO_M} (D_{poly}^{k}, \cO_M)\,,
\end{equation}
where $\cO_M$ denotes the structure sheaf
of (smooth) functions on $M$ and
$D_{poly}^{\bul}$ is the sheaf of
polydifferential operators.

Note that although
$$
J_k(M) = C^{poly}_{k+1}(M), \qquad k\ge 0
$$
I would like to reserve special notation
for the bundle (\ref{polyjets}) and distinguish
$\JM$ and $\Cp$. Let me, from now on, refer to elements of
$\Cp$ as {\it Hochschild chains} and to elements
of $\JM$ as {\it polyjets}.

The bundle $J_{\bul}$ is endowed
with a canonical flat connection $\nabla^G$ which is
called the \emph{Grothendieck connection} and defined by
the formula
\begin{equation}
\label{eq:gro}
\nabla^G_u(j)(P):= u (j(P))-j(u \bul P)\,,
\end{equation}
where $j\in J_k(M)$,
$P\in D^k_{poly}(M)$, and $u$ is
a vector field which is viewed,
in the right hand side, as a differential operator.
The operation $\bul$ is
defined in (\ref{bullet}).

For this connection we have
the following remarkable proposition:
\begin{pred}
\label{pr:chi}
Let $\chi$ be a linear map ($k\ge 0$)
$$
\chi: J_k(M) \to C^{poly}_{k}(M)
$$
defined by the formula
\begin{equation}
\label{chi}
\chi(a)(P)=a(1\otimes P)\,, \quad
P \in D^{k-1}_{poly}(M)\,, \quad a\in J_{k}(M)\,.
\end{equation}
The restriction of the
map $\chi$
to the $\n^G$-flat polyjets
gives the ($\bbR$-linear) isomorphism
($k\ge 0$)
\begin{equation}
\label{chi1}
\chi: \ker \n^G \cap J_{k}(M)
\erarrow
C^{poly}_{k}(M)\,.
\end{equation}
\end{pred}
{\bf Proof.} To see that the map (\ref{chi1})
is surjective one has to notice that
for any Hoch\-schild chain
$b\in C^{poly}_{k}(M)$
the equations
$$
a(1\otimes P) = b(P)\,,  \qquad P\in D_{poly}^{k-1}(M)
$$
and
\begin{equation}
\label{reshaiu}
a (u\cdot Q\otimes P) = u\, a(Q\otimes P)
- a (Q\otimes ( u \bul P))\,,
\end{equation}
$$
Q\in D_{poly}^0(M)\,, \qquad u \in \G(M, TM)
$$
define a $\n^G$-flat polyjet $a$ of
degree $k$\,.

On the other hand, if $a$ is a $\n^G$-flat polyjet of
degree $k$ equation (\ref{reshaiu}) is
automatically satisfied. Thus $a$ is uniquely
determined by its image $\chi(a)$. $\Box$

Let $t$ be the cyclic permutation acting on
the sheaf $J_{\bul}$ of polyjets
\begin{equation}
\label{eq:t}
t(a)(P_0\otimes\cdots\otimes P_l):=
a (P_1\otimes\cdots\otimes P_l\otimes P_0)\,,
\end{equation}
$$
a\in J_l(M)\,, \qquad
P_i \in D_{poly}^0(M)\,.
$$

Using this operation and the bilinear
product (\ref{bullet}) I define the
map
$$
\hR : D^{k}_{poly} \otimes J_{l} \to J_{l-k}\,, \qquad
P \otimes a \mapsto \hR_P(a)\,,
$$
$$
\hR_P(a)(Q_0\otimes Q)=a((Q_0\otimes Q) \bullet P)+
$$
\begin{equation}
\label{ono}
\sum_{j=1}^k (-1)^{lj}t^j(a)
\big( (Q_0 \bul P) \otimes Q \big)\,,
\end{equation}
$$
P\in D^{k}_{poly}(M)\,, \qquad a\in J_{l}(M)\,,
\qquad Q\in D^{l-k-1}_{poly}(M)\,,
\qquad Q_0 \in D^0_{poly}(M)\,.
$$

Following the lines of the proof of
lemma $2.3$ in \cite{Ezra} one can show that
\begin{equation}
\label{hR-OK}
\hR_{[P_1, P_2]_G}= \hR_{P_1} \hR_{P_2} -
(-)^{|P_1||P_2|} \hR_{P_2} \hR_{P_1}\,,
\end{equation}
and hence, (\ref{polyjets}) is a sheaf of graded
modules over the sheaf of graded Lie algebras $D_{poly}^{\bul}$\,.
Furthermore, using the multiplication
$\mu_0\in D^1_{poly}(M)$ in $\OM$ one can
turn the $D_{poly}^{\bul}$-module (\ref{polyjets})
into a DGLA $D_{poly}^{\bul}$-module by
introducing the following differential
\begin{equation}
\label{hb-chain}
\hmb  = \hR_{\mu_0} : J_k \mapsto J_{k-1}\,.
\end{equation}

It follows from the construction that
both the action (\ref{ono}) and the differential
(\ref{hb-chain}) commute with the Grothendieck
connection (\ref{eq:gro}). Thus the $\n^G$-flat
polyjets $\ker \n^G \cap J_{\bul}(M)$ form
a DG module over the DGLA $\Dp$\,.

A direct but slightly tedious computation
shows that
\begin{pred}
\label{jet-chain}
The DG module structure on $\Cp$ over
the DGLA $\Dp$ induced from (\ref{ono})
and (\ref{hb-chain}) via the isomorphism
(\ref{chi1}) coincides with standard one
given by (\ref{cochain-act}) and
(\ref{b-chain}). $\Box$
\end{pred}
This proposition allows me to identify
$\n^G$-flat polyjets $\ker\n^G \cap \JM$
and Hochschild chains $\Cp$ as DGLA modules.
This identification will be very handy
for the construction of the Fedosov resolution
for $\Cp$\,.

Unfortunately, the maps (\ref{U-1}) and
(\ref{C-1}) are not compatible with the Lie brackets
on $\Tp$ and $\Dp$ and with
the respective actions (\ref{cochain-act}) and
(\ref{Tpoly-act}). In particular, the equation
\begin{equation}
\mC \circ R_{\cV(\ga)} \stackrel{?}{=}
L_{\ga} \mC
\end{equation}
does not hold in general. In \cite{Tsygan}
B. Tsygan suggested that this defect
could be cured by the following
statement:
\begin{conj}[B. Tsygan, \cite{Tsygan}]
For any smooth manifold $M$
the DGLA modules $(\Tp, \AM)$ and
$(\Dp, \Cp)$ are quasi-isomorphic.
\end{conj}
The following theorem gives a positive
answer to the question of B. Tsygan.
\begin{teo} \label{thm-chain}
For any smooth manifold $M$
there exists a commutative
diagram of DGLAs and DGLA modules
\begin{equation}
\begin{array}{ccccccc}
\Tp & \brarrow &\cL_1 &  \brarrow  & \cL_2 & \blarrow & \Dp \\[0.3cm]
\downarrow_{\,mod}  & ~  & \downarrow_{\,mod}& ~ &\downarrow_{\,mod}  & ~
& \downarrow_{\,mod} \\[0.3cm]
\AM & \bbrarrow &\M_1 &  \bblarrow  & \M_2 & \bblarrow & \Cp\,,
\end{array}
\label{diag-thm}
\end{equation}
in which the horizontal arrows in the upper
row are quasi-isomorphisms of
DGLAs and the horizontal arrows in the lower
row are quasi-isomorphisms of $\Linf$-modules.
The terms ($\cL_1$, $\cL_2$,
$\M_1$, $\M_2$) and  the quasi-isomorphisms
of diagram (\ref{diag-thm}) are functorial for
diffeomorphisms of pairs ``manifold $M$ $+$
a torsion free connection on $TM$''.
\end{teo}
The construction of the quasi-iso\-mor\-phisms
in diagram (\ref{diag-thm}) is explicit and in chapter $5$
I show how this result allows us to prove Tsygan's conjecture
(see the first part of corollary $4.0.3$ in \cite{Tsygan})
about Hochschild homology of the
quantum algebra of functions on an arbitrary Poisson manifold,
and in particular, to describe the space of
traces on this algebra.

The main part of the proof of theorem \ref{thm-chain}
concerns the construction of
Fedosov resolutions of the DGLA modules
$(\Tp, \cA^{\bul}(M))$ and $(\Dp,\Cp)$.
After completing this stage it will only remain to use
Kontsevich's \cite{K}
and Shoikhet's \cite{Sh} formality theorems for $\bbR^{d}_{formal}$
and apply the twisting procedures developed in the
previous chapter.

Let us now recall these formality theorems.
\begin{teo}[M. Kontsevich, \cite{K}] \label{aux}
There exists a quasi-iso\-mor\-phism $\cK$
\begin{equation}
\label{auxform}
\cK\,:\,T^{\bul}_{poly}(\bbR^d)\brarrow D^{\bul}_{poly}(\bbR^d)
\end{equation}
from the
DGLA $T^{\bul}_{poly}(\bbR^d)$ of po\-ly\-vec\-tor fields
to the DGLA $D^{\bul}_{poly}(\bbR^d)$ of po\-ly\-dif\-fe\-ren\-tial
operators on the space $\bbR^d$ such that
\begin{enumerate}
\item One can replace $\bbR^d$ in (\ref{auxform}) by its formal completion $\bbRf$
at the origin.

\item The quasi-iso\-mor\-phism $\cK$ is equivariant with respect to
linear transformations of the coordinates on $\bbRf$\,.

\item If $n>1$ then
\begin{equation}
\label{vanish}
\cK_{n}(u_1, u_2, \dots, u_n)=0
\end{equation}
for any set of vector fields
$u_1, u_2, \dots, u_n\in T^0_{poly}(\bbRf)$\,.

\item If $n\ge 2$ and $u \in T^0_{poly}(\bbRf)$ is linear in the
coordinates on $\bbRf$ then for any set of po\-ly\-vec\-tor fields
$\ga_2, \dots, \ga_n\in T^{\bul}_{poly}(\bbRf)$
\begin{equation}
\label{vanish1}
\cK_{n}(u,\ga_2, \dots, \ga_n)=0\,.
\end{equation}
\end{enumerate}
\end{teo}
Composing the quasi-isomorphism $\cK$ with the action
(\ref{cochain-act}) of $D^{\bul}_{poly}(\bbR^d)$ on
$C^{poly}_{\bul}(\bbR^d)$
I get an $\Linf$-module structure on $C^{poly}_{\bul}(\bbR^d)$ over
the DGLA $T^{\bul}_{poly}(\bbR^d)$. For this module structure
we have the following results:
\begin{teo}[B. Shoikhet, \cite{Sh}] \label{aux1}
There exists a quasi-iso\-mor\-phism $\cS$
\begin{equation}
\label{auxform-Sh}
\cS\,:\,C^{poly}_{\bul}(\bbR^d) \bbrarrow \cA^{\bul}(\bbR^d)
\end{equation}
of $\Linf$-modules over $T_{poly}(\bbR^d)$, the zeroth
structure map $\cS_0$ of which is the map (\ref{C-1}) of Connes
and such that
\begin{enumerate}
\item One can replace $\bbR^d$ in (\ref{auxform-Sh}) by its formal completion $\bbRf$
at the origin.

\item The quasi-iso\-mor\-phism $\cS$ is equivariant with respect to
linear transformations of the coordinates on $\bbRf$\,.

\end{enumerate}
\end{teo}

\begin{figure}
\includegraphics[width=1.6in, height=1.6in]{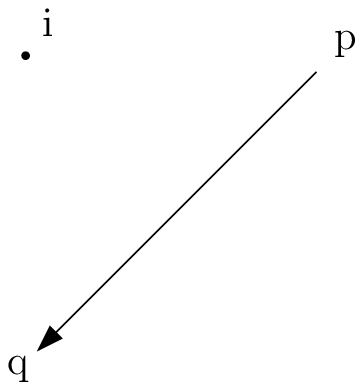}
\hfill
\includegraphics[width=1.6in, height=1.6in]{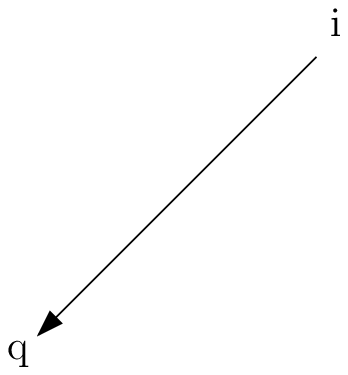}
\\
\parbox[t]{2.6in}{\caption{Edge of the first type} \label{fig1}}
\hfill
\parbox[t]{2.6in}{\caption{Edge of the second type} \label{fig2}}
\end{figure}

\begin{pred}
\label{John}
If $\cS$ be the quasi-isomorphism (\ref{auxform-Sh})
of B. Shoikhet, $n\ge 1$, and $u\in T^0_{poly}(\bbRf)$ is linear in the
coordinates on $\bbRf$ then for any set of po\-ly\-vec\-tor fields
$\ga_2, \dots, \ga_n\in T^{\bul}_{poly}(\bbRf)$ and any Hochschild
chain $a\in C^{poly}_{\bul}(\bbRf)$
\begin{equation}
\label{vanish-Sh}
\cS_{n}(u,\ga_2, \dots, \ga_n ; a)=0\,.
\end{equation}
\end{pred}
{\bf Proof.} The proof of (\ref{vanish-Sh}) reduces to
calculation of integrals entering the construction of
the structure maps $\cS_n$ (see section $2.2$ of \cite{Sh}).
To do this calculation I first
transform the unit disk $\{|\zeta |\le 1\}$
used in section $2.2$ of \cite{Sh}
into the upper half
plane $\cH^{+}=\{z,\, Im\,(z)\ge 0\}$ via the
standard fractional linear
transformation
\begin{equation}
\label{transf}
z= -i \frac{\zeta +1}{\zeta -1}\,.
\end{equation}
The origin of the unit disk goes to $z=i$ and
the point $\zeta =1$ goes to $z=\infty$\,.
The angle function corresponding to
an edge of the first type \cite{Sh} (see figure \ref{fig1})
connecting $p\neq i$ and $q\neq i$ looks as follows
\begin{equation}
\label{angle}
\al^{Sh}(p,q)=Arg(p-q)- Arg(\bar{p}-q)-
 Arg(p-i) + Arg(\bar{p}-i)\,.
\end{equation}
If I fix the rotation symmetry by placing
the first function of the Hochschild chain at
the point $z=\infty$ then the angle function
corresponding to an edge of the second type
(see figure \ref{fig2}) connecting $p=i$ and $q$
takes the form
\begin{equation}
\label{angle1}
\beta^{Sh}(q)=Arg(i-q)- Arg(-i-q)\,.
\end{equation}

\begin{figure}
\centering\includegraphics{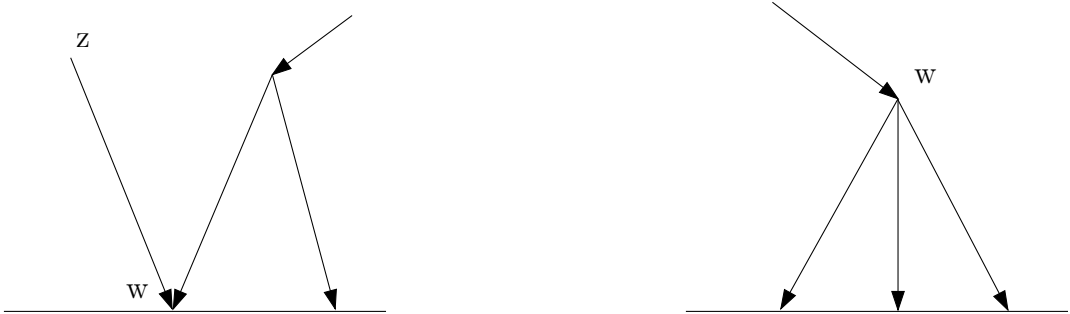}
\vspace{0.26in}
\caption{Diagrams of the first type} \label{fig3}
\end{figure}

\begin{figure}
\vspace{0.5in}
\centering\includegraphics{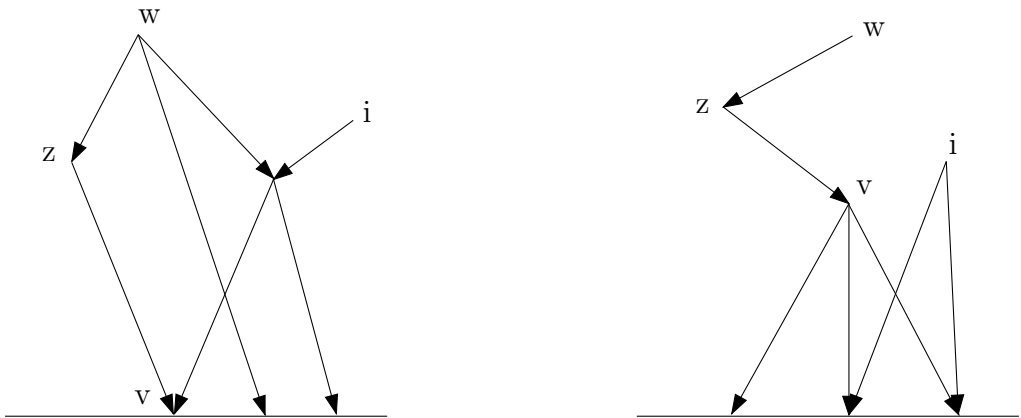}
\vspace{0.5in}
\caption{Diagrams of the second type} \label{fig4}
\end{figure}

\begin{figure}
\vspace{0.5in}
\centering\includegraphics{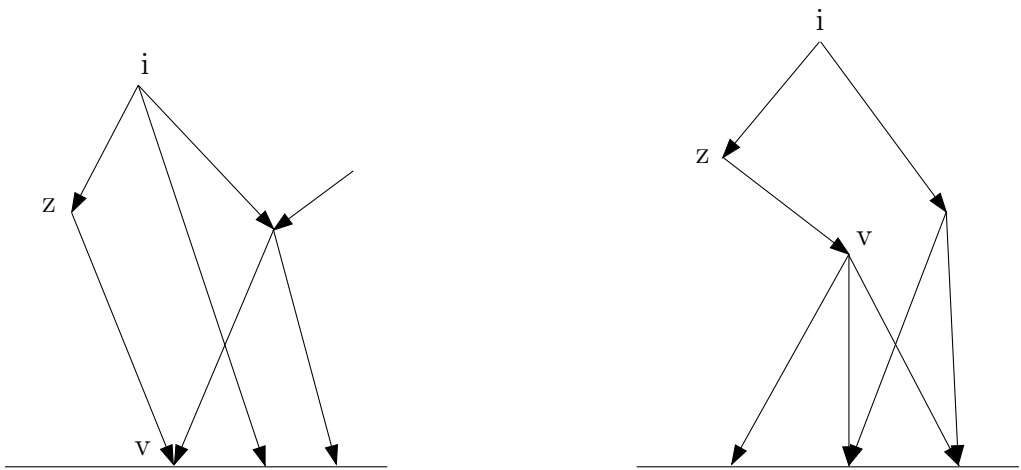}
\vspace{0.5in}
\caption{Diagrams of the third type} \label{fig5}
\end{figure}

Let us suppose that $u$ is a vector linear in
coordinates on $\bbRf$. Then there are three types of
the diagrams corresponding to $\cS_n(u, \dots)$ $n\ge 2$\,.
In the diagram of the first type (see figure \ref{fig3}) there are no
edges ending at the vertex $z$ corresponding to the
vector $u$. In the diagrams of the second type
(see figure \ref{fig4}) there is exactly one edge ending
at the vertex $z$ and this wedge does not start at
the vertex $i$. In the diagrams of the third type
(see figure \ref{fig5}) there is exactly one edge ending
at the vertex $z$ and this wedge starts at
the vertex $i$.

The coefficient corresponding to
a diagram of the first type vanishes because
the angle functions entering the
integrand form turn out
to be dependent. The coefficients  corresponding to
diagrams of the second and the third type vanish
since so do the following integrals
\begin{equation}
\label{integral}
\int_{z\in \cH^+\setminus\{w, v, i\}} d\al^{Sh}(w,z)
d\al^{Sh}(z,v)=0\,, \quad
\int_{z\in \cH^+\setminus\{v, i\}} d\beta^{Sh}(z)
d\al^{Sh}(z,v)=0\,.
\end{equation}
Equations (\ref{integral}) follow immediately from lemmas
$7.3$, $7.4$, and $7.5$ in \cite{K}. $\Box$

~\\
{\bf Remark.} Hopefully, alternative proofs of
theorems \ref{aux}, \ref{aux1}, and proposition \ref{John}
may be obtained along the lines of Tamarkin and
Tsygan \cite{Dima, TT, TT1}.

%%
%% CHAPTER 3 ENDS
%%
%% CHAPTER 4  1 mile
%%
%%

\chapter{Fedosov resolutions of the DGLA modules
$(\Tp$, $\cA^{\bul}(M))$ and $(\Dp$,  $\Cp)$}
In paper \cite{Fedosov} B. Fedosov proposed a simple geometric
construction for star-products on an arbitrary
symplectic manifold. The key idea of
Fedosov's construction has various
incarnations and it is referred to as the Gelfand-Fuchs
trick \cite{GF} or formal geometry \cite{G-Kazh}
in the sense of I.M. Gelfand and D.A. Kazhdan , or
mixed resolutions \cite{Ye11}. This idea can
be roughly formulated as the following slogan:
``In order to linearize a problem one has to
formulate it in terms of jets''.

If $M$ is smooth manifold the bundle of jets
$J_0$ (\ref{polyjets}) is non-canonically isomorphic
to the bundle $\SM$ of the formally
completed symmetric algebra of the cotangent bundle $T^*M$.
For this reason I start with the definition of
this bundle.
\begin{defi} The bundle $\SM$ of the formally completed symmetric algebra
 of the cotangent bundle $T^*M$
is defined as a bundle over the ma\-ni\-fold $M$  whose sections
are infinite collections of symmetric covariant tensors $a_{i_1\dots
i_p}(x)$\,, where $x^i$ are local coordinates, $p$
runs from $0$ to $\infty$\,, and the indices
$i_1,\dots, i_p$ run from $1$ to $d$\,.
\end{defi}
It is convenient to introduce auxiliary variables $y^i$\,, which
transform as contravariant vectors. These variables allow us to
rewrite any section $a\in \G(M, \SM)$ in the form
of the formal power series
\begin{equation}
\label{sect}
a=a(x,y)=\sum_{p=0}^{\infty} a_{i_1\dots
i_p}(x)y^{i_1}\dots y^{i_p}\,.
\end{equation}

It is easy to see that the vector space $\G(M, \SM)$ is
naturally endowed with the commutative product which is induced by a
fiberwise multiplication of formal power series in $y^i$\,. This product
makes $\G(M, \SM)$ into a commutative algebra with a unit.

Now I recall from \cite{CEFT} definitions of
formal fiberwise po\-ly\-vec\-tor fields and formal fiberwise
po\-ly\-dif\-fe\-ren\-tial operators on $\SM$\,.
\begin{defi}
A bundle $\cT^k_{poly}$ of formal fiberwise po\-ly\-vec\-tor
fields of degree $k$ is a bundle over $M$ whose sections
are $C^{\infty}(M)$-linear operators
$\mv : \wedge^{k+1} \G(M, \SM) \mapsto \G(M, \SM)$ of the form
\begin{equation}
\label{vect}
\mv =\sum_{p=0}^{\infty}\mv^{j_0\dots j_k}_{i_1\dots i_p}(x)y^{i_1}
\dots y^{i_p} \frac{\pa}{\pa y^{j_0}}\wedge  \dots \wedge \frac{\pa}{\pa
y^{j_k}}\,,
\end{equation}
where I assume that the infinite sum in $y$'s is formal and
$\mv^{j_0\dots j_k}_{i_1\dots i_p}(x)$ are tensors symmetric in
indices $i_1, \dots, i_p$ and antisymmetric in indices
$j_0, \dots, j_k$\,.
\end{defi}
Extending the definition of the formal fiberwise po\-ly\-vec\-tor
field by allowing the fields to be inhomogeneous
I define the total bundle $\cT_{poly}$ of
formal fiberwise po\-ly\-vec\-tor fields
\begin{equation}
\label{cal-T}
\cT_{poly} =\bigoplus_{k=-1}^{\infty} \cT_{poly}^k\,, \qquad
\cT_{poly}^{-1}=\SM\,.
\end{equation}
The fibers of the bundle $\cT_{poly}$ are endowed with
the DGLA structure $T_{poly}(\bbRf)$ of po\-ly\-vec\-tor fields on
the formal completion $\bbRf$ of $\bbR^d$ at the origin.
This turns $\cTp$ into a sheaf of DGLAs (with the vanishing
differential).

\begin{defi}
A bundle $\cD^k_{poly}$ of formal fiberwise po\-ly\-dif\-fe\-ren\-tial
operator of degree $k$ is a bundle over $M$ whose sections
are $C^{\infty}(M)$-polylinear
maps $\mP : \bigotimes^{k+1} \G(M, \SM) \mapsto \G(M, \SM)$ of the form
\begin{equation}
\label{operr}
\mP =\sum_{\al_0 \dots \al_k}\sum_{p=0}^{\infty}\mP^{\al_0\dots \al_k}_{i_1\dots i_p}(x)y^{i_1}
\dots y^{i_p} \frac{\pa}{\pa y^{\al_0}}\otimes  \dots \otimes \frac{\pa}{\pa
y^{\al_k}}\,,
\end{equation}
where $\al$'s are multi-indices $\al={j_1\dots j_l}$ and
$$
\frac{\pa}{\pa y^{\al}}=\frac{\pa}{\pa y^{j_1}} \dots \frac{\pa}{\pa
y^{j_l}}\,,
$$
the infinite sum in $y$'s is formal, and the
sum in the orders of derivatives $\pa/\pa y$
is finite.
\end{defi}
Notice that the tensors
$\mP^{\al_0\dots \al_k}_{i_1\dots i_p}(x)$ are
symmetric in covariant indices $i_1,\dots, i_p$\,.

As well as for po\-ly\-vec\-tor fields I define the total bundle
$\cD_{poly}$ of formal fiberwise po\-ly\-dif\-fe\-ren\-tial operators
as the direct sum
\begin{equation}
\label{cal-D}
\cD_{poly} =\bigoplus_{k=-1}^{\infty} \cD^k_{poly}\,, \qquad
\cD^{-1}_{poly}=\SM\,.
\end{equation}
The fibers of the bundle $\cD_{poly}$ are endowed with
the DGLA structure (and
DGA structure) $D_{poly}(\bbRf)$ of po\-ly\-dif\-fe\-ren\-tial operators
on $\bbRf$\,. This turns $\cDp$ into a sheaf of DGLAs and
a sheaf of DGAs.

\begin{defi}
A bundle $\cC^{poly}_{k}$ of formal fiberwise Hochschild chains
of degree $k$ $(k\ge 0)$ is a bundle over $M$ whose sections
are formal power series in $k+1$ collections of fiber coordinates
$y^i_0, \dots, y_k^i$ of the tangent bundle
\begin{equation}
\label{f-chain}
a(x,y_0, \dots, y_k)=\sum_{\al_0 \dots \al_k}
a_{\al_0\dots \al_k}(x)y_0^{\al_0} \dots y^{\al_k}_k\,,
\end{equation}
where $\al$'s are multi-indices $\al={j_1\dots j_l}$ and
$$
y^{\al}= y^{j_1}y^{j_2} \dots y^{j_l}\,.
$$
\end{defi}
The total bundle
$\cC^{poly}$ of formal fiberwise Hochschild chains
is the direct sum
\begin{equation}
\label{cal-C}
\cC^{poly} =\bigoplus_{k=0}^{\infty} \cC_{k}^{poly}\,, \qquad
\cC_{0}^{poly}=\SM\,.
\end{equation}
The operations $R$ (\ref{cochain-act}) and
$\mb$ (\ref{b-chain}) turn each fiber of $\cC^{poly}$
into a DGLA $D_{poly}(\bbRf)$-module.
Thus $\cCp$ is a sheaf of DG modules over
the sheaf of DGLAs $\cDp$.

As above, I denote by $\cA^{\bul}(M)$
the space of exterior forms
\begin{equation}
\label{cA}
\cA^{\bul}(M)=\bigoplus_{k=0}^{\infty} \cA^{k}(M)\,,
\qquad
\cA^{k}(M)=\{a=a_{i_1\dots i_k}(x) dx^{i_1}\dots dx^{i_k}\}\,.
\end{equation}
Furthermore,
\begin{defi}
The bundle $\cE$ of fiberwise exterior
forms is a bundle over $M$ whose sections
are exterior forms with values in $\SM$.
These sections are given
by the following formal power series
\begin{equation}
\label{cE}
a(x,y,d x)=\sum_{p,q\ge 0} a_{i_1\dots
i_p\, ;\, j_1 \dots j_q}(x) y^{i_1} \dots y^{i_p} dx^{j_1} \dots dx^{j_q}\,,
\end{equation}
where $a_{i_1\dots i_p \,;\, j_1 \dots j_q}(x)$ are components
of covariant tensors symmetric in indices
$i_1,\dots, i_p$ and antisymmetric in indices $j_1,\dots,
j_q$\,.
\end{defi}
The fiberwise analogue of the Lie derivative
(\ref{Tpoly-act}) allows me to speak of $\cE$ as
of a sheaf of modules over the sheaf of DGLAs $\cTp$\,.

For my purposes I will also need
``exterior forms with values in exterior forms''.
This forces me to introduce an additional copy
$\{dy^i \}$ of the local basis $\{ dx^i \}$
of exterior forms on $M$.
Having these two copies
I reserve the notation $\Om^{\bul}(M,\cB)$ for
the graded vector space of $d y$-exterior forms with
values in the bundle $\cB$. In particular, I
would like to distinguish the graded vector
spaces $\OmS$ and $\G(M, \cE)$\,. $\OmS$ consists
of {\it $d y$-forms} and $\G(M, \cE)$ consists of
{\it $d x$-forms}.

For the relations between $dx^i$ and $dy^j$ I accept
the following convention
$$
dx^i dy^j = - dy^j dx^i\,.
$$

Homogeneous elements of the graded vector spaces
$\OmT$ and $\OmD$ are the following formal series in $y$'s
\begin{equation}
\label{Om-vect}
\mv =\sum_{p\ge 0} dy^{l_1} \dots dy^{l_q}
\mv^{j_0\dots j_k}_{l_1\dots l_q \, ; \, i_1\dots i_p}(x)y^{i_1}
\dots y^{i_p} \frac{\pa}{\pa y^{j_0}}\wedge  \dots \wedge \frac{\pa}{\pa
y^{j_k}}\,,
\end{equation}
and
\begin{equation}
\label{Om-operr}
\mP =\sum_{\al_0 \dots \al_k}\sum_{p\ge 0}
dy^{l_1} \dots dy^{l_q} \mP^{\al_0\dots \al_k}_{l_1\dots l_q\, ; \, i_1\dots i_p}(x)y^{i_1}
\dots y^{i_p} \frac{\pa}{\pa y^{\al_0}}\otimes  \dots \otimes \frac{\pa}{\pa
y^{\al_k}}\,,
\end{equation}
where as above $\al$'s are multi-indices $\al={j_1\dots j_l}$ and
$$
\frac{\pa}{\pa y^{\al}}=\frac{\pa}{\pa y^{j_1}} \dots \frac{\pa}{\pa
y^{j_l}}\,.
$$
Similarly, homogeneous elements of $\OmE$ and $\OmC$ are
the formal series
\begin{equation}
\label{Om-cE}
a(x,d y,y, d x)=\sum_{p\ge 0}dy^{l_1} \dots dy^{l_q} a_{l_1\dots l_q \, ; \, i_1\dots
i_p j_1 \dots j_k}(x)y^{i_1}\dots y^{i_p}dx^{j_1} \dots dx^{j_k}\,,
\end{equation}
and
\begin{equation}
\label{Om-f-chain}
b(x,dy, y_0, \dots, y_k)=\sum_{\al_0 \dots \al_k}
dy^{l_1} \dots dy^{l_q} b_{l_1\dots l_q \, ; \, \al_0\dots \al_k}(x)y_0^{\al_0}
\dots y^{\al_k}_k\,,
\end{equation}
where as above $\al$'s are multi-indices $\al={j_1\dots j_l}$ and
$$
y^{\al}= y^{j_1}y^{j_2} \dots y^{j_l}\,.
$$
The symmetries of tensor indices in formulas
(\ref{Om-vect}), (\ref{Om-operr}), (\ref{Om-cE}), and
(\ref{Om-f-chain}) are obvious.

The space $\OmS$ is naturally endowed with the structure
of a $\bbZ$-graded commutative algebra and it is
also filtered with respect to the powers in  $y$'s.
The graded vector spaces $\OmT$ and $\OmD$ are, in turn,
endowed with fiberwise DGLA structures induced by those
on $T_{poly}(\bbRf)$ and $D_{poly}(\bbRf)$\,.
Similarly, $\OmE$ and $\OmC$ become fiberwise
DGLA modules over $\OmT$ and $\OmD$,
respectively\footnote{I regard $\OmT$ and $\OmE$ as
a DGLA and a DGLA module with vanishing differentials.}.
I denote the Lie bracket in $\OmD$
by $[,]_G$ and the Lie bracket in $\OmT$ by $[,]_{SN}$\,.
For fiberwise Lie derivative on $\OmE$ and for the
fiberwise action of $\OmD$ on $\OmC$ I also use
the same notation $L$ and $R$, respectively.
It is not hard to see that the formulas for
the fiberwise differentials on $\OmD$ and $\OmC$
can be written similarly to (\ref{pa}) and
(\ref{b-chain})
$$
\pa = [\mu, \cdot \,]\,, \qquad
\mb=R_{\mu}\,,
$$
where $\mu \in \G(M, \cD^1_{poly})$ is the (commutative) multiplication
in $\G(M, \SM)$\,. Notice that $\OmD$ is also
endowed with a fiberwise DGA structure induced by that
on $D_{poly}(\bbRf)$\,.

The parity of elements in the algebras $\OmT$, $\OmD$ and the
modules $\OmE$ and $\OmC$ is defined by
the sum of the exterior degree and the degree
in the respective fiberwise algebra or
the respective fiberwise module.

The following proposition shows that
I have a distinguished sheaf of graded
Lie algebras which acts on the sheaves
$\SM$, $\cTp$, $\cE$,
$\cDp$, and $\cCp\,.$
\begin{pred}
\label{ono1}
$\cTp^0$ is a sheaf of graded Lie
algebras. $\SM$, $\cE$, $\cTp$, $\cDp$, and
$\cCp$ are sheaves of modules over $\cTp^0$ and
the action of $\cTp^0$ is
compatible with the (DG) algebraic structures
on $\SM$, $\cE$, $\cTp$, $\cDp$, and $\cCp$\,.
\end{pred}
{\bf Proof.} Since the Schouten-Nijenhuis bracket (\ref{eq:Leib})
has degree zero $\cTp^0$ $\subset$ $\cTp$
$\subset$ $\cDp$ is a subsheaf of graded Lie algebras.
While the action of $\cTp^0$ on the
sections of $\SM$ is obvious, the action on
$\cE$ is given by the Lie derivative, the
action on $\cTp$ is the adjoint action corresponding
to the Schouten-Nijenhuis bracket, the action on
$\cDp$ is given by the Gerstenhaber bracket and
the action on $\cCp$ is induced by the
action of Hochschild cochains on
Hochschild chains (\ref{cochain-act}).
The compatibility of the action with the cup product
(\ref{cup}) in $\cDp$ essentially follows
from the fact that $\cTp^0$ acts by derivations
on the sheaf of algebras $\SM$. The compatibility
with the remaining DG algebraic structures follows
from the definitions. $\Box$

This proposition implies that the
canonical vector field $\displaystyle dy^i \frac{\pa}{\pa y^i}$
$\in$ $\Om^1(M,\cTp^0)$ defines the differential
\begin{equation}
\label{del}
\de= dy^i \frac{\pa}{\pa y^i} ~ \cdot  \,:\, \Om^{\bul}(M,\cB) \mapsto
\Om^{\bul+1}(M,\cB)\,,
\qquad \de^2=0\,,
\end{equation}
where $\cB$ is either of the bundles $\SM$,
$\cTp$, $\cDp$, $\cE$, or $\cCp$ and $\cdot$
denotes the corresponding action of $\cTp^0$.
Due to the above proposition the differential
$\de$ is compatible with the corresponding DG algebraic
structures.

The subspaces $\FT$ and $\FD$
will subsequently play an important role in
our construction. They can be described in the
following way. Elements of $\FT$ are
fiberwise po\-ly\-vec\-tor fields (\ref{vect})
$$
\mv =\sum_k \mv^{j_0\dots j_k}(x)
\frac{\pa}{\pa y^{j_0}}\wedge  \dots \wedge
\frac{\pa}{\pa y^{j_k}}\,
$$
whose components do not depend on $y$'s.
Similarly, elements of $\FD$ are fiberwise polydifferential
operators (\ref{operr})
$$
\mP =\sum_k\sum_{\al_0 \dots \al_k}
\mP^{\al_0\dots \al_k}(x)\frac{\pa}{\pa y^{\al_0}}\otimes
\dots \otimes \frac{\pa}{\pa y^{\al_k}}\,
$$
whose coefficients do not depend on $y$'s.

In the following proposition I
describe cohomology of the differential $\de$
in $\OmS$, $\OmT$, $\OmD$, and $\OmE$
\begin{pred}
\label{Cohom-del}
For $\cB$ be either of the bundles
$\SM$, $\cTp$, $\cDp$, or $\cE$
$$
H^{>0}(\Omb(M, \cB), \de) =0\,.
$$
Furthermore,
$$
H^0(\OmT,\de)=\FT\,,
$$
$$
H^0(\OmD,\de)=\FD\,,
$$
$$
H^0(\OmS,\de)= C^{\infty}(M)\,,
$$
$$
H^0(\OmE, \de)= \cA^{\bul}(M)\,.
$$
\end{pred}
{\bf Proof.}
The proposition will follow immediately
if I construct an operator
$$
\de^{-1} \,:\,\Om^{\bul}(M,\cB)\mapsto \Om^{\bul-1}(M,\cB)
$$
such that for any $a\in \Omb(M, \cB)$
\begin{equation}
a=\sigma(a) +\delta \delta ^{-1}a + \delta ^{-1}\delta a\,,
 \label{Hodge}
\end{equation}
where
\begin{equation}
\label{si}
\si a= a\Big|_{y^i=dy^i=0}\,.
\end{equation}

First, I define this operator on $\OmS$
\begin{equation}
\de^{-1}(a) =
\begin{cases}
\begin{array}{cc}
\displaystyle
y^k \frac {\vec{\partial}} {\partial (d y^k)}
\int\limits_0^1 a(x,t y,t dy)\frac{dt} t, & {\rm if}~
a \in \Om^{>0}(M, \SM)\,,\\
0, & {\rm otherwise}\,,
\end{array}
\end{cases}
\label{del-1}
\end{equation}
where the arrow over $\pa$ denotes
the left derivative with respect to
the anti-com\-muting variable $dy^k$.

Next, I extend $\de^{-1}$ to
the vector spaces $\OmE$, $\OmT$, $\OmD$
in the componentwise manner.
A direct computation shows that
equation (\ref{Hodge}) holds and the
proposition follows. $\Box$

It is worth noting that the operator $\de^{-1}$
is $2$-nilpotent for either of complexes
\begin{equation}
\label{needed}
(\de^{-1})^2=0\,.
\end{equation}

For our purposes I fix an affine torsion free connection
$\n$ on $M$. Since the bundles $\SM$,
$\cTp$, $\cDp$, $\cE$, or $\cCp$ are obtained from
the tangent bundle the connection $\n$ extends to
them in the natural way. I use the same
notation for all these connections
\begin{equation}
\label{nab}
\n \,:\, \Om^{\bul}(M,\cB) \mapsto  \Om^{\bul+1}(M,\cB)\,,
\end{equation}
where $\cB$ is either $\SM$, $\cTp$, $\cDp$, $\cE$,
or $\cCp$\,.

The following statement is
an easy exercise of differential geometry.
\begin{pred}
\label{nabla}
Let $\cB$ be either $\SM$, $\cTp$, $\cDp$, $\cE$,
or $\cCp$ and let $\cdot$ denote
the action of $\cTp^0$.
Then the connection $\n$ is given by
the following operator
\begin{equation}
\label{nab1}
\n= dy^i \frac{\pa}{\pa x^i} + \G \, \cdot
\,:\,  \Om^\bul (M, \cB) \mapsto  \Om^{\bul+1} (M, \cB)\,,
\end{equation}
where
\begin{equation}
\label{Christ}
\G= -dy^i \G^k_{ij}(x) y^j \frac{\pa}{\pa y^k}\,,
\end{equation}
and $\G^k_{ij}(x)$ are the corresponding
Christoffel symbols.
Furthermore,
\begin{equation}
\label{nab-sq}
\n^2 a = \cR \,\cdot a
\,:\, \Om^{\bul}(M,\cB) \mapsto \Om^{\bul+2}(M,\cB)\,,
\end{equation}
where
$$
\cR= -\frac12 dy^i dy^j (R_{ij})^k_l(x) y^l \frac{\pa}{\pa y^k}\,,
$$
and $(R_{ij})^k_l(x)$ is the standard Riemann curvature tensor
of the connection $\n$. $\Box$
\end{pred}
Notice that, due proposition \ref{ono1}
the operator (\ref{nab1}) is compatible with the
(DG) algebraic structures on $\OmS$, $\OmT$, $\OmE$,
$\OmD$, and $\OmC$\,. Moreover, since the connection
$\n$ is torsion free the derivations
(\ref{del}) and (\ref{nab1}) (anti)commute
\begin{equation}
\label{anticomm}
\de \n + \n \de = 0\,.
\end{equation}

I would like to combine
the operators (\ref{del}) and
(\ref{nab1}) into a $2$-nilpotent
derivation
\begin{equation}
\label{DDD}
D=\n - \de + A\, \cdot \,:\, \Om^{\bul}(M,\cB)
\mapsto
\Om^{\bul+1}(M,\cB)\,,
\end{equation}
where $\cB$ and $\cdot$ are as in
proposition \ref{nabla} and
$$
A=\sum_{p=2}^{\infty}dy^k A^j_{ki_1\dots i_p}(x) y^{i_1} \dots
y^{i_p}\frac{\pa}{\pa y^j} \in \Om^1(M,\cT^0_{poly})\,.
$$
is a $d y$- $1$-form with values in
the fiberwise vector fields $\cTp^0$.

Due to the following theorem it is
always possible to find the $1$-form
$A$ such that the derivation (\ref{DDD})
is $2$-nilpotent\,.
\begin{teo}
Iterating the equation
\begin{equation}
\label{iter_A}
A=\de^{-1} \cR + \de^{-1}(\n A +\frac12 [A,A]_{SN})
\end{equation}
in degrees in $y$ one constructs
$A\in \Om^1(M,\cT^0_{poly})$
such that $\de^{-1}A=0$ and the derivation $D$ (\ref{DDD})
is $2$-nilpotent
$$
D^2=0\,.
$$
\end{teo}
In what follows I refer to the differential $D$ (\ref{DDD})
as {\it the Fedosov differential}.

~\\
{\bf Proof.} First, I observe that the
recurrent procedure in (\ref{iter_A}) converges to
an element $A\in \Om^1(M,\cTp^0)$
since the operator $\de^{-1}$ raises the degree in $y$.
Moreover, due to equation (\ref{needed})
\begin{equation}
\label{de-1-A}
\de^{-1}A=0\,.
\end{equation}

Second, the equation $D^2=0$
is equivalent to
\begin{equation}
\label{flat}
\cR -\de A  + \n A +\frac12 [A,A]_{SN}=0\,.
\end{equation}

Denoting by $C\in \Om^2(M, \cTp^0)$
the left hand side of (\ref{flat})
$$
C= -\de A + \cR + \n A +\frac12 [A,A]_{SN}
$$
using (\ref{Hodge}), (\ref{iter_A}), and (\ref{de-1-A})
one gets that
\begin{equation}
\label{flat1}
\de^{-1} C=0\,.
\end{equation}

On the other hand $[D,D^2]=0$ and hence
\begin{equation}
\label{svo}
\n C - \de C + [A, C]_{SN}=0\,.
\end{equation}
Thus, applying (\ref{Hodge}) to $C$ and
using (\ref{flat1}) one gets
the equation
$$
C= \de^{-1} (\n C + [A,C]_{SN})\,.
$$
This equation has the unique vanishing
solution since the operator $\de^{-1}$ raises the degree in $y$\,.
The theorem is proved. $\Box$

In the next theorem I compute cohomology
of the Fedosov differential (\ref{DDD})
for $\OmS$, $\OmT$, $\OmE$, and $\OmD$.
\begin{teo}
\label{teo}
If $\cB$ is either
$\SM$, $\cE$, $\cTp$, or $\cDp$
then
\begin{equation}
\label{H-D>0}
H^{>0}(\Omb(M,\cB), D) =0\,.
\end{equation}
Furthermore,
\begin{equation}
\label{H-D0}
\begin{array}{c}
H^0(\Om(M,\SM), D) \cong \OM \,, \\[0.3cm]
H^0(\Om(M,\cE ), D) \cong \AM \\[0.3cm]
\end{array}
\end{equation}
as graded commutative algebras,
\begin{equation}
\label{H-D01}
H^0(\Om(M,\cTp), D) \cong  \FT
\end{equation}
as graded vector spaces, and
\begin{equation}
\label{H-D011}
H^0(\Om(M,\cDp), D) \cong   \FD
\end{equation}
as graded associative algebras.
\end{teo}
{\bf Proof.}
Although the first statement follows easily from
the spectral sequence argument I need a more
explicit proof.

To prove (\ref{H-D>0}) I construct
an $\bbR$-linear map
\begin{equation}
\label{Phi}
\Phi : \Om^{\bul}(M, \cB) \mapsto \Om^{\bul-1}(M,\cB)
\end{equation}
such that for any $a\in \Om^{>0}(M,\cB)$
\begin{equation}
\label{eq:Phi}
D \Phi (a) + \Phi D(a) = a\,.
\end{equation}

I define the map $\Phi$ with the help of
the following recurrent procedure
\begin{equation}
\label{nado-mne}
\Phi(a) = -\de^{-1} a+ \de^{-1}(\n \Phi(a)  + A\cdot \Phi(a))\,,
\end{equation}
where $\cdot$ denotes the action of
$\cTp^0$ (see proposition \ref{ono1}) and
the procedure (\ref{nado-mne})
converges since $\de^{-1}$ (\ref{del-1}) raises the
degree in the fiber coordinates $y^i$\,.

Due to equation (\ref{needed})
$\de^{-1}\Phi(a)=0$ and therefore
\begin{equation}
\label{nado-ochen}
\Phi^2 =0\,.
\end{equation}

Let me prove that for any element
$a\in \Om^{>0}(M,\cB)$ $\cap$ $\ker D$
\begin{equation}
\label{Exact}
a= D \Phi (a)\,.
\end{equation}
For this I denote by $h$ the element
$$
h=a- D \Phi(a) \in \Om^{>0}(M,\cB)
$$
and mention that $D h=0$ or equivalently
\begin{equation}
\label{Closed}
\de h = \n h + A \cdot h\,.
\end{equation}

Since $\de^{-1} \Phi(a)=0$ and
$\si(\Phi(a))=0$ equation (\ref{Hodge})
for $\Phi(a)$ boils down to
$$
\Phi(a) = \de^{-1} \de \Phi(a)\,.
$$
Thus, using (\ref{nado-mne}), I conclude that
$$
\de^{-1} h=0\,.
$$
Furthermore, since $h \in \Om^{>0}(M, \cB)$
$$
\si h=0\,.
$$
Hence applying (\ref{Hodge}) to $h$
and using (\ref{Closed}) I get
$$
h= \de^{-1} (\n h + A\cdot h)\,.
$$
The latter equation has the unique vanishing solution
since $\de^{-1}$ raises the degree in the
fiber coordinates $y^i$\,. Thus (\ref{Exact}) is
proved.

Using (\ref{Exact}) I conclude that
\begin{equation}
\label{Exact1}
D \circ \Phi \circ D = D\,.
\end{equation}

Let me now turn to our combination
$$
b= a- D \Phi(a) - \Phi D(a)\,,
$$
where $a\in \Om^{>0}(M,\cB)$\,.

Thanks to (\ref{Exact1}) and $D^2=0$
$$
D b=0\,.
$$
Hence, applying (\ref{Exact}) to $b$ I get
$$
b = D \Phi(b)\,.
$$
Using (\ref{nado-ochen}) and (\ref{Exact1})
once again I get that $b=0$, and therefore,
(\ref{eq:Phi}) holds.

Thus the first statement (\ref{H-D>0}) is proved.

Let $\cH$ denote either $\OM$, $\AM$,
$\FT$, or $\FD$ and $\cdot$ denote the
action of $\cTp^0$ (see proposition
\ref{ono1}) on $\SM$, $\cE$, $\cTp$,
and $\cDp$. I claim that iterating the
equation
\begin{equation}
\label{iter_a}
\tau(a)=a + \de^{-1}(\n \tau(a)+ A\cdot \tau(a))\,,
\qquad a \in \cH
\end{equation}
in degrees in $y$ one gets an isomorphism
\begin{equation}
\label{lift}
\tau : \cH \mapsto \ker D\cap  \G(M,\cB)\,.
\end{equation}
Indeed, let $a\in \cH$. Then, due
to formula (\ref{Hodge}) $\tau(a)$ satisfies the
following equation
\begin{equation}
\label{de-1-Y}
\de^{-1} (D (\tau(a))) =0\,.
\end{equation}
Let us denote $D \tau(a)$ by $Y$
$$
Y= D \tau(a)\,.
$$
The equation $D^2=0$ implies that
$D Y = 0$, or in other words
\begin{equation}
\label{del-Y}
\de Y = \n Y + A \cdot Y\,.
\end{equation}
Applying (\ref{Hodge}) to $Y$ and
using equations (\ref{de-1-Y}), (\ref{del-Y})
I get
$$
Y = \de^{-1} (\n Y + A\cdot Y)\,.
$$
The latter equation has the unique vanishing
solution since the operator $\de^{-1}$ (\ref{del-1})
raises the degree in the fiber coordinates $y^i$.

The map (\ref{lift}) is
injective since $\si$ (\ref{si})
is a section of (\ref{lift})
\begin{equation}
\label{si-tau}
\si \circ \tau = Id\,.
\end{equation}
To prove surjectivity of (\ref{lift})
it suffices to show that if
$b\in \G(M, \cB)\cap \ker D$
and
\begin{equation}
\label{si-b}
\si b = 0
\end{equation}
then $b$ vanishes.

The condition $b\in \ker D$ is
equivalent to the equation
$$
\de b = \n b + A \cdot b\,.
$$
Hence, applying (\ref{Hodge}) to
$a$ and using (\ref{si-b})
I get
$$
b = \de^{-1} (\n b + A\cdot b)\,.
$$
The latter equation has the unique vanishing
solution since the operator $\de^{-1}$ (\ref{del-1})
raises the degree in the fiber coordinates $y^i$.
Thus, the map (\ref{lift}) is bijective and
the map $\si$ (\ref{si}) provides me with
the inverse of (\ref{lift})
\begin{equation}
\label{tau-inv}
\tau \circ \si \Big|_ {\ker D\, \cap \, \G(M,\cB)} = Id\,.
\end{equation}

Since $\si$ respects the
multiplications in $\OmS$, $\OmE$, $\cDp^{\bul}(M)$,
$\OM$, $\AM$, and $\FD$ so does
the map $\tau$ and the theorem
follows. $\Box$

Notice that since the Fedosov differential (\ref{DDD})
is compatible with the DGLA structure on $\OmT$
and $\OmD$, the cohomology groups $H^{\bul}(\OmT,D)$ and
$H^{\bul}(\OmD,D)$ acquire structures of
a DGLA, and $H^{\bul}(\OmD,D)$ also becomes a DGA.
To analyze these structures
let me, first, observe that for any
function $a\in C^{\infty}(M)$ and for any integer $p\ge 0$
\begin{equation}
\label{deriv}
\frac{\pa}{\pa y^{i_1}} \dots \frac{\pa}{\pa y^{i_p}} \tau (a)
\Big|_{y=0}= \pa_{x^{i_1}} \dots \pa_{x^{i_p}} a(x)+
{\rm lower~order~ derivatives~ of}~a\,.
\end{equation}
Due to this observation the following map
$$
\nu\,:\, \FD \mapsto \Dp\,,
$$
\begin{equation}
\label{eq:nu}
\nu (\mP)(a_0, \dots, a_k) =
\Big( \mP (\tau(a_0), \dots, \tau(a_k)) \Big)\Big|_{y=0}\,,
\end{equation}
$$
\mP \in \ker \de \cap \G(M, \cDp^k)\,,\qquad
a_i \in \OM
$$
is an isomorphism of graded
associative algebras $\FD$ and  $\Dp$\,.

I claim that
\begin{pred}
\label{Pr-la-D}
The composition
\begin{equation}
\label{la-D}
\la_D =  \tau \circ \nu^{-1} : \Dp \mapsto \ker D \cap \G(M, \cDp^{\bul})
\end{equation}
induces an isomorphism from the
DGLA (and DGA) $\Dp$ to the DGLA (and DGA) $H^{\bul}(\OmD,D)$\,.
\end{pred}
{\bf Proof.} Since both the map $\tau$ (\ref{lift})
and the map $\nu$ (\ref{eq:nu})
respect the cup-product (\ref{cup})
it suffices to prove the compatibility
with the DGLA structures.
I will prove that inverse map
$$
\la_D^{-1}=
\nu \circ \si : \ker D \cap \G(M, \cDp^{\bul})
\mapsto \Dp\,,
$$
\begin{equation}
\label{la-D-1}
\la_D^{-1} (\mP)(a_0, \dots, a_k) =
\Big( \mP (\tau(a_0), \dots, \tau(a_k))\Big)\Big|_{y=0}\,,
\end{equation}
$$
\mP \in \ker D \cap \G(M, \cDp^k)\,,\qquad
a_i \in \OM
$$
respects the Gerstenhaber
bracket (\ref{Gerst}) and the Hochschild
differential (\ref{pa}).

To prove the compatibility with
the bracket I observe that
applying $\tau$ (\ref{lift}) to both
sides of (\ref{la-D-1}) and using
(\ref{tau-inv}) one gets
\begin{equation}
\label{D-mP-a}
\tau \Big(\la_D^{-1} (\mP)(a_0, \dots, a_k)\Big) =
\mP (\tau(a_0), \dots, \tau(a_k))\,,
\end{equation}
$$
\forall \qquad
\mP\in \ker D \cap \G(M, \cDp^k)\,, \qquad
a_i\in \OM\,.
$$
Using this equation I conclude that
for any $\mP_1\in  \ker D \cap \G(M, \cDp^{k_1})$
and $\mP_2\in  \ker D \cap \G(M, \cDp^{k_2})$
$$
\la_D^{-1}( \mP_1 )
(a_0, \dots, \la_D^{-1}( \mP_2 )(a_i,  \dots, a_{i+k_2}),
\dots, a_{k_1+k_2})=
$$
$$
\mP_1
(\tau(a_0), \dots,  \mP_2 ( \tau(a_i),  \dots, \tau(a_{i+k_2}) ),
\dots, \tau(a_{k_1+k_2}))\,.
$$
Therefore, for any $\mP_1\in  \ker D \cap \G(M, \cDp^{k_1})$
and $\mP_2\in \ker D \cap \G(M, \cDp^{k_2})$
\begin{equation}
\label{la-bullet}
\la_D^{-1}(\mP_1)\bul
\la_D^{-1}(\mP_2) = \la_D^{-1}(\mP_1 \bul \mP_2)\,,
\end{equation}
where the operation $\bul$ is defined in
(\ref{bullet}).

Thus $\la_D^{-1}$ is compatible with the
Gerstenhaber bracket (\ref{Gerst}).

To prove that $\la_D^{-1}$ respects
the differentials (\ref{pa}) in
$\Dp$ and $H^{\bul}(\OmD, D)$ it suffices to show
that the multiplication $\mu\in \G(M, \cDp^1)$
in $\G(M,\SM)$ is sent to the multiplication
$\mu_0\in D^1_{poly}(M)$ in $\OM$\,.
This is immediate from the definition of
$\la_D^{-1}$ (\ref{la-D-1}). Thus the
proposition is proved. $\Box$

It is obvious that the restriction of the
map $\nu$ (\ref{eq:nu}) to $\G(M, \cTp)$ gives a map
\begin{equation}
\label{eq:nu-T}
\nu : \FT \mapsto \Tp\,.
\end{equation}
By the abuse of notation I denote this
map by the same letter.

It is easy to see that due to
equation (\ref{deriv}) the map (\ref{eq:nu-T}) is also
an isomorphism of graded vector spaces.
Furthermore,
\begin{pred}
\label{Pr-la-T}
The composition
\begin{equation}
\label{la-T}
\la_T =  \tau \circ \nu^{-1} : \Tp \mapsto \ker D \cap \G(M, \cTp^{\bul})
\end{equation}
induces an isomorphism from the
graded Lie algebra $\Tp$ to the
graded Lie algebra $H^{\bul}(\OmT,D)$\,.
\end{pred}
{\bf Proof.} To show that $\la_T$ is compatible
with Lie brackets I observe that the following
diagram
\begin{equation}
\begin{array}{ccc}
\Tp &\stackrel{\la_T}{\,\rightarrow\,} &   \G(M, \cTp^{\bul})  \\[0.3cm]
\downarrow^{\cV}   & ~  &     \downarrow^{\cV^{fib}} \\[0.3cm]
\Dp   &\stackrel{\la_D}{\,\rightarrow\,} &  \G(M, \cDp^{\bul}),
\end{array}
\label{eq:diag}
\end{equation}
commutes. Here $\cV$ is the map of Vey (\ref{U-1}) and
$\cV^{fib}$ denotes its fiberwise analogue.

Thus for any pair $\ga_1, \ga_2\in \Tp$ I have
$$
\cVf(\la_T([\ga_1, \ga_2]_{SN}) -
[\la_T(\ga_1), \la_T(\ga_2)]_{SN})=
\la_D \cV([\ga_1, \ga_2]_{SN}) -
[\cVf \la_T(\ga_1), \cVf \la_T(\ga_2)]_{SN}
$$
modulo $\pa$-exact terms in $\G(M, \cDp^{\bul})$\,.
Continuing this line of equations and using
proposition \ref{Pr-la-D} I conclude
that
$$
\cVf (\la_T([\ga_1, \ga_2]_{SN}) -
[\la_T(\ga_1), \la_T(\ga_2)]_{SN}) \in \pa (\G(M, \cDp^{\bul}))\,.
$$
Therefore, since $\cVf$ is a quasi-isomorphism
of complexes $(\G(M, \cTp^{\bul}),0)$ and
$(\G(M, \cDp^{\bul}), \pa )$
$$
\la_T([\ga_1, \ga_2]_{SN}) -
[\la_T(\ga_1), \la_T(\ga_2)]_{SN}=0
$$
and the proposition follows. $\Box$

Since the Fedosov differential (\ref{DDD}) is compatible
with the DGLA module structures on $\OmE$ and $\OmC$
the cohomology groups $H^{\bul}(\OmE, D)$ and
$H^{\bul}(\OmC, D)$ acquire the DGLA module structures
over $H^{\bul}(\OmT, D)$ and  $H^{\bul}(\OmD, D)$,
respectively. Due to theorem \ref{teo} and
propositions \ref{Pr-la-D}, \ref{Pr-la-T}
$H^{\bul}(\OmT, D)\cong \Tp$ and
$H^{\bul}(\OmD, D)\cong \Dp$ as DGLAs.
My next task is to show that $H(\OmE,D)\cong \AM$
and $H(\OmC,D)\cong \Cp$ as modules over
the corresponding DGLAs.

The desired statement about chains follows from
proposition \ref{jet-chain} and
\begin{pred}
\label{vot-ono}
For any $q > 0$
\begin{equation}
\label{H0D-C}
H^{q} (\OmC ,D) = 0\,.
\end{equation}
The map
\begin{equation}
\label{map-vr}
\vr : \G(M, \cCp_{\bul}) \to \JM\,,
\qquad
\vr(a)(P)= (\la_D (P))(a)\Big|_{y^i=0}\,,
\end{equation}
$$
a\in \G(M, \cCp_k)\,, \quad P\in D^k_{poly}(M)
$$
is an isomorphism of the DG modules over the DGLA
$\Dp \cong \ker D \cap \G(M, \cDp^{\bul})$.
Moreover, this isomorphism sends the Fedosov
connection (\ref{DDD}) on $\cCp$ to
the Gro\-then\-dieck connection
(\ref{eq:gro}) on $J_{\bul}$\,.
\end{pred}
{\bf Proof.}
The first statement (\ref{H0D-C}) follows easily from the spectral sequence
argument. Indeed, using the zeroth collection of the fiber coordinates
$y_0^i$ (\ref{Om-f-chain}) I introduce the decreasing filtration
on the sheaf $\OmC$
$$
\dots \subset F^p (\OmC) \subset
F^{p-1}(\OmC)
\subset \dots \subset F^0(\OmC)=\OmC\,,
$$
where the components of the forms (\ref{Om-f-chain})
in $F^p(\OmC)$ have degree in $y_0^i$
$\ge p$.

Since $D (F^p(\OmC)) \subset F^{p-1}(\OmC)$ the corresponding
spectral sequence starts with
$$
E_{-1}^{p,q} = F^p (\Om^{p+q}(M, \cCp))\,.
$$
Next, I observe that
$$
d_{-1} = dy^i \frac{\pa}{\pa y_0^i}\,,
$$
and hence, due to the Poincar\'e lemma for the formal disk
I have
$$
E_0^{p,q} = E_1^{p,q} = \dots = E_{\infty}^{p,q} =0
$$
whenever $p+q> 0$.
Thus, the first statement (\ref{H0D-C}) of the
proposition follows.

Thanks to observation (\ref{deriv}) the map
(\ref{map-vr}) is indeed an isomorphism of
graded vector spaces.

The compatibility with the action (\ref{cochain-act})
of the sheaf of DGLAs $\cDp^{\bul}$ on
the sheaf $\cCp_{\bul}$ and with the action
(\ref{ono}) of the sheaf of DGLAs $D^{\bul}_{poly}$
on the sheaf $J_{\bul}$
$$
\hR_P (\vr(a)) = \vr (R_{\la_D(P)}(a) )
$$
follows from the compatibility of $\la_D$
with the operation $\bul$ (\ref{bullet})
(see (\ref{la-bullet})), with the
cyclic permutations, and
with the cup products (\ref{cup}) in $\Dp$
and $\cDp^{\bul}(M)$ (see proposition \ref{Pr-la-D}).

It remains to prove that the map (\ref{map-vr})
sends the Fedosov connection (\ref{DDD}) to the Grothendieck
connection (\ref{eq:gro}). This statement is proved
by the following line of equations:
\begin{align*}
\vr(D_u a)(P)& = (\la_D (P))(D_u a)\Big|_{y^i=0} =
(D_u [\la_D (P)(a)]) \Big|_{y^i=0}\cr
&= u [\la_D (P)(a)]\Big|_{y^i=0}
- ( i_u \de \, \bul \, [\la(P)(a)] )\Big|_{y^i=0}\cr
&= u [\la_D (P)(a)]\Big|_{y^i=0}
- (\la_D (u)\, \bul \, \la_D(P))(a)\Big|_{y^i=0}
\cr
&= u [\la_D (P)(a)]\Big|_{y^i=0}
- \la_D (u \, \bul \,  P)(a)\Big|_{y^i=0}
\cr
&= u (\vr(a))(P) - (\vr(a))(u \, \bul \,  P) =
(\n^G_u \vr(a))(P)\,,
\end{align*}
where $u\in \G(M,TM)$, $a \in \G(M, \cCp_{k})$,
$P\in D^k_{poly}(M)$\,, $i$ denotes the
contraction of a vector field with
differential forms, $\bul$ is as in (\ref{bullet}),
and $u$ is viewed
both as a vector field and a differential
operator of the first order. $\Box$

Let me conclude this chapter with
\begin{pred}
\label{now}
The map (\ref{lift})
\begin{equation}
\label{nowE}
\tau\,:\, \AM \mapsto \OmE
\end{equation}
induces an isomorphism of DG modules
$\AM$ and $H^{\bul}(\OmE, D)$ over
the DGLA $H^{\bul}(\OmT, D)\cong \Tp$\,.
\end{pred}
{\bf Proof.} I have to prove that for any
exterior form $a=a_{i_1 \dots i_q}(x)dx^{i_1} \dots dx^{i_q}$
and any polyvector field
$\ga=\ga^{i_0 \dots i_k}(x) \pa_{x^{i_0}} \wedge \dots \wedge
\pa_{x^{i_k}}$
\begin{equation}
\label{nado}
\tau (L_{\ga}(a))= L_{\tau\circ\nu^{-1}(\ga)} (\tau (a))\,.
\end{equation}
Since Fedosov differential $D$ is compatible
with the fiberwise Lie derivative $L$
the form $L_{\tau\circ\inu (\ga)} (\tau (a))$ is $D$-closed.
Therefore by (\ref{tau-inv}) it suffices to the show that
\begin{equation}
\label{och-nado}
L_{\tau\circ\nu^{-1}(\ga)} (\tau (a))\Big|_{y=0}=L_{\ga}(a)\,.
\end{equation}
To prove (\ref{och-nado}) I need the expressions
for $\tau(\inu(\ga))$ and $\tau(a)$ only up to the
second order terms in $y$. They are
\begin{equation}
\label{tau-ga}
\tau (\inu(\ga)) = \inu(\ga) + y^i \frac{\pa \inu(\ga)}{\pa x^i} -
y^i [\,\G_i(x), \inu(\ga)\,]_{SN} \quad
mod \quad (y)^2\,,
\end{equation}
\begin{equation}
\label{tau-form}
\tau (a) = a  + y^i \frac{\pa a}{\pa x^i}  -
y^i L_{\G_i(x)}(a)\quad
mod \quad (y)^2\,,
\end{equation}
$$
\G_i=\G^k_{ij}(x)y^j\pa_{y^k}\,,
$$
where $\G^k_{ij}(x)$ are Christoffel symbols and
$$
\inu(\ga)=\ga^{i_0 \dots i_k}(x) \pa_{y^{i_0}} \wedge \dots \wedge
\pa_{y^{i_k}}\,.
$$
Using symmetry of indices for the Christoffel symbols
$\G^k_{ij}= \G^k_{ji}$ I can rewrite (\ref{tau-ga}) and
(\ref{tau-form}) in the form
\begin{equation}
\label{tau-ga1}
\tau (\inu(\ga)) = \inu(\ga) + y^i \frac{\pa \inu(\ga)}{\pa x^i} -
[\,\tG(x), \inu(\ga)\,]_{SN} \quad
mod \quad (y)^2\,,
\end{equation}
\begin{equation}
\label{tau-form1}
\tau (a) = a  + y^i \frac{\pa a}{\pa x^i}  -
L_{\tG(x)}(a)\quad
mod \quad (y)^2\,,
\end{equation}
where $\displaystyle
\tG=\frac12 \G^k_{ij}y^i y^j\frac{\pa}{\pa y^k}\,.$
Using these formulas it is not hard to show that
equation (\ref{och-nado}) is equivalent to
$$
L_{\inu(\ga)}L_{\tG} (a) + L_{[\tG,\inu(\ga)]_{SN}} (a)=0\,,
$$
which obviously holds because $L_{\inu(\ga)} (a)=0$\,. $\Box$

%%
%% CHAPTER 4 ENDS
%%
%% CHAPTER 5  1 mile
%%
%%

\chapter{Formality theorems for $(\Dp$, $\Cp)$ and their applications}
\section{Proof of the formality theorem for $\Cp$}

The results of the previous chapter can be represented in the
form of the following commutative diagrams of DGLAs,
their modules, and quasi-isomorphisms
given by honest (not $\Linf$) morphisms
\begin{equation}
\begin{array}{ccc}
\Tp &\stackrel{\la_T}{\,\longrightarrow\,} &(\OmT, D, [,]_{SN})\\[0.3cm]
\downarrow^{L}_{\,mod}  & ~  &     \downarrow^{L}_{\,mod} \\[0.3cm]
\AM   &\stackrel{\la_{\cA}}{\,\longrightarrow\,} & (\OmE, D),\\[1cm]
(\OmD, D+\pa, [,]_{G}) &\stackrel{\,\la_D}{\,\longleftarrow\,} & \Dp\\[0.3cm]
\downarrow^{R}_{\,mod}  & ~  &     \downarrow^{R}_{\,mod} \\[0.3cm]
(\OmC, D+\mb) &\stackrel{\,\la_C}{\,\longleftarrow\,} &   \Cp,
\end{array}
\label{diag-T-D}
\end{equation}
where $\displaystyle\la_T= \tau\circ \inu \Big|_{\Tp}$,
$\displaystyle \la_{\cA}=\tau\Big|_{\AM}$,
$\la_D = \tau\circ \inu$,
$\la_C = \vr^{-1} \circ\chi^{-1}$, the map
$\chi$ is defined in (\ref{chi}), the map
$\tau$ is defined in (\ref{iter_a}), and
the map $\vr$ is defined in (\ref{map-vr}).

Next, due to properties $1$ and $2$ in theorem \ref{aux} I have
a fiberwise quasi-iso\-mor\-phism (which I denote
by the same letter $\cK$)
\begin{equation}
\label{cal-K}
\cK :  (\OmT,0,[,]_{SN}) \brarrow (\OmD,\pa,[,]_G)\,.
\end{equation}
from the DGLA $(\OmT,0,[,]_{SN})$ to the DGLA $(\OmD,\pa,[,]_G)$\,.

Due to properties $1$ and $2$ in theorem \ref{aux1} I have
a fiberwise quasi-iso\-mor\-phism (which I denote
by the same letter $\cS$)
\begin{equation}
\label{cal-S}
\cS :  (\OmC,\mb) \bbrarrow (\OmE,0,L)
\end{equation}
from the $\Linf$-module $\OmC$ to the DGLA module
$\OmE$ over $\OmT$\,.

Thus I get the following commutative diagram
\begin{equation}
\begin{array}{ccc}
(\OmT, 0, [,]_{SN}) &\stackrel{\cK}{\brarrow} &(\OmD, \pa, [,]_{G})\\[0.3cm]
\downarrow^{L}_{\,mod}  & ~  &     \downarrow^{R}_{\,mod} \\[0.3cm]
(\OmE, 0)  &\stackrel{\cS}{\bblarrow} & (\OmC, \mb),
\end{array}
\label{diag-K-Sh}
\end{equation}
where by commutativity I
mean that $\cS$ is a morphism of the $\Linf$-modules
$(\OmC, \mb)$ and $(\OmE, 0)$ over the DGLA $(\OmT, 0, [,]_{SN})$,
and the $\Linf$-module structure on $(\OmC, \mb)$
over $(\OmT, 0, [,]_{SN})$ is obtained by composing
the quasi-isomorphism $\cK$ with the action $R$ of
$(\OmD, \pa, [,]_{G})$ on $(\OmC, \mb)$\,.

Having the complete decreasing
filtration on $\OmT$, $\OmD$, $\OmE$, and $\OmC$
induced by the exterior degree I can
now apply the technique developed in
section \ref{section-MC}.
To do this I first restrict myself to
an open coordinate subset
$$
V\subset M\,.
$$
On $V$ it makes sense to speak about the ordinary De
Rham differential\footnote{Let me recall that $\Om(\cdot)$
stands for $d y$-forms.} in the DGLA modules
$(\Omb(V,\cTp), \Omb(V, \cE))$, and
$(\Omb(V, \cDp), \Omb(V, \cCp))$
\begin{equation}
\label{de-Rham}
d=dy^i \pa_{x^i} : \Om^{\bul}(V,\cB) \mapsto  \Om^{\bul+1}(V,\cB)\,,
\end{equation}
where $\cB$ is either $\cTp$ or $\cDp$, $\cE$, or
$\cCp$\,.

Since the quasi-isomorphisms (\ref{cal-K}) and (\ref{cal-S})
are fiberwise I can add to all the differentials
in diagram (\ref{diag-K-Sh}) the
De Rham differential (\ref{de-Rham}),
and thus, get the new commutative diagram
\begin{equation}
\begin{array}{ccc}
(\Omb(V,\cT_{poly}), d, [,]_{SN}) &\stackrel{\cK}{\brarrow} &
(\Omb(V, \cD_{poly}), d+\pa, [,]_{G})\\[0.3cm]
\downarrow^{L}_{\,mod}  & ~  &     \downarrow^{R}_{\,mod} \\[0.3cm]
(\Omb(V, \cA), d)  &\stackrel{\cS}{\bblarrow} & (\Omb(V, \cCp), d + \mb)\,.
\end{array}
\label{diag-V}
\end{equation}
I claim that
\begin{pred}
The $\Linf$-morphism $\cK$
and the morphism of $\Linf$-modules $\cS$
in (\ref{diag-V}) are quasi-isomorphisms.
\end{pred}
{\bf Proof.}
This statement follows easily from
the standard argument of the spectral sequence.
Indeed, the $\Linf$-morphism $\cK$
(resp. the morphism of $\Linf$-modules $\cS$)
is compatible the descending filtration
induced by the exterior degree
\begin{equation}
\label{ext-filtr}
\cF^p(\Omb(V,\cB)) = \bigoplus_{k\ge p} \Om^k(V,\cB)\,,
\end{equation}
where $\cB$ is either $\cTp$ or $\cDp$ (resp.
$\cE$ or $\cCp$).

The corresponding versions of Vey's \cite{Vey} and
Hoch\-schild-Kos\-tant-Ro\-sen\-berg-Con\-nes-Te\-le\-man
 \cite{Connes}, \cite{HKR}, \cite{Tel} theorems for
$\bbRf$ imply that the first
structure map $\cK_1$ (resp. the zeroth structure
map $\cS_0$)
induces a quasi-isomorphism on the level of
$E_0$. Therefore, $\cK_1$ (resp. $\cS_0$) induces a
quasi-isomorphism on the level of $E_{\infty}$.
The standard snake lemma argument of homological
algebra implies that $\cK_1$ (resp. $\cS_0$) is
a quasi-isomorphism. Hence, so is $\cK$
(resp. $\cS$). $\Box$

On the subset $V$ I can represent the Fedosov
differential (\ref{DDD}) in the following form
$$
D=d+ B \, \cdot \,,
$$
where
\begin{equation}
\label{d+B}
B=\sum^{\infty}_{p=0} dy^i B^k_{i\, ; \, j_1 \dots j_p}(x) y^{j_1} \dots
y^{j_p} \frac{\pa}{\pa y^k} \in \Om^1(V,\cT^0_{poly})\,,
\end{equation}
and $\cdot$ denotes the action of the
sheaf $\cT^0_{poly}$. (See proposition \ref{ono1}.)

The nilpotency condition $D^2=0$ says that $B$ is a
Maurer-Cartan element of the DGLA
$(\Omb(V,\cT_{poly}),d,[,]_{SN})$ with the
filtration (\ref{ext-filtr}).
Thus using the terminology of section \ref{section-MC}
one can say that the DGLA $(\Omb(V,\cT_{poly}),D,[,]_{SN})$
is obtained from $(\Omb(V,\cT_{poly}),d,[,]_{SN})$ by
twisting via $B$.

Due to property $3$ in theorem \ref{aux} the Maurer-Cartan
element in $(\Om(V,\cD_{poly}),d+\pa,[,]_{G})$
$$
B_D=\sum_{m=1}^{\infty}\frac{1}{m!} \cK_m(B, \dots, B)
$$
corresponding to the Maurer-Cartan element $B$
in $(\Om(V,\cT_{poly}),d,[,]_{SN})$ coincides with $B$
viewed as an element of $\Om^1(V, \cD_{poly})$\,.
Thus twisting of the quasi-isomorphism $\cK$
via the Maurer-Cartan element $B$ I get
the $\Linf$-morphism
\begin{equation}
\label{K-tw}
\cK^{tw}\,:\,(\Omb(V,\cT_{poly}),D,[,]_{SN}) \brarrow
(\Omb(V,\cD_{poly}),D+\pa,[,]_{G})\,,
\end{equation}
which is a quasi-isomorphism due to
claim $5$ of proposition \ref{twist-morph}.

Next, using (\ref{twist-vf-str}) it is not hard to see
that the graded module structures on $\Omb(V,\cE)$ and  $\Omb(V, \cC^{poly})$
over $(\Omb(V$,$\cT_{poly})$,$D$,$[,]_{SN})$ and
$(\Omb(V$,$\cD_{poly})$,$D+\pa$,$[,]_{G})$, respectively,
remain unchanged under the twisting procedures, while
the differentials get shifted. Namely, $d$ on $\Omb(V,\cE)$
gets replaced by $D$ and $d+\mb$ on $\Omb(V, \cC^{poly})$
gets replaced by $D+ \mb$\,.

Hence, by virtue of propositions \ref{twist-mod} and
\ref{functor} twisting procedure turns diagram
(\ref{diag-V}) into the commutative diagram
\begin{equation}
\begin{array}{ccc}
(\Omb(V,\cT_{poly}), D, [,]_{SN}) &\stackrel{\cK^{tw}}{\brarrow} &
(\Omb(V, \cD_{poly}), D+\pa, [,]_{G})\\[0.3cm]
\downarrow^{L}_{\,mod}  & ~  &     \downarrow^{R}_{\,mod} \\[0.3cm]
(\Omb(V,\cE), D)  &\stackrel{\cS^{tw}}{\bblarrow} & (\Omb(V,\cC^{poly}), D + \mb),
\end{array}
\label{diag-V1}
\end{equation}
where $\cS^{tw}$ is morphism of
$\Linf$-modules obtained from $\cS$
by twisting via the Maurer-Cartan element
$B\in \Om^1(V,\cTp)$\,. Due to
claim $5$ of proposition \ref{twist-mod}
$\cS^{tw}$ is a quasi-isomorphism.

Surprisingly, due to property $4$ in theorem \ref{aux}
and proposition \ref{John} the ``morphisms''
$\cK^{tw}$ and $\cS^{tw}$ are defined globally. Indeed,
using (\ref{twist-F-str}) and (\ref{twist-ka-str}) I get
the structure maps of $\cK^{tw}$ and $\cS^{tw}$
\begin{equation}
\label{cK-tw}
\cK^{tw}_n(\ga_1, \dots,\,\ga_n)=
\sum_{k=0}^{\infty} \frac1{k!}
\cK_{n+k} (B,\dots,\, B, \ga_1, \dots, \, \ga_n)\,,
\end{equation}
\begin{equation}
\label{cS-tw}
\cS^{tw}_n(\ga_1, \dots,\,\ga_n, a)=
\sum_{k=0}^{\infty} \frac1{k!}
\cS_{n+k} (B,\dots,\, B, \ga_1, \dots, \, \ga_n, a)\,,
\end{equation}
$$
\ga_i \in \Om(V, \cT_{poly})\,,\qquad  a\in \Om(V, \cC^{poly})
$$
in terms of the structure maps of $\cK$ and $\cS$\,.
But the only term in $B$ that does not transform as
a tensor is
\begin{equation}
\label{Gamma}
\G=- dy^i\G^k_{ij}y^j \frac{\pa}{\pa y^k}\,.
\end{equation}
This term contributes neither to $\cK^{tw}_n$
nor to $\cS^{tw}_n$ since it is linear in $y$'s.

Thus the quasi-isomorphisms $\cK^{tw}$ and $\cS^{tw}$ are defined
globally and I arrive at the following commutative diagram
\begin{equation}
\begin{array}{ccc}
(\Omb(M,\cT_{poly}), D, [,]_{SN}) &\stackrel{\cK^{tw}}{\brarrow} &
(\Omb(M, \cD_{poly}), D+\pa, [,]_{G})\\[0.3cm]
\downarrow^{L}_{\,mod}  & ~  &     \downarrow^{R}_{\,mod} \\[0.3cm]
(\OmE, D)  &\stackrel{\cS^{tw}}{\bblarrow} & (\OmC, D + \mb).
\end{array}
\label{diag-M}
\end{equation}
Assembling (\ref{diag-M}) with (\ref{diag-T-D}) I
get the desired commutative diagram
\begin{equation}
\begin{array}{ccccccc}
\Tp &  \stackrel{\la_T}{\lrarrow} & \OmT &
\stackrel{\cK^{tw}}{\brarrow} & \OmD &
\stackrel{\,\la_D}{\llarrow} & \Dp\\[0.3cm]
\downarrow^{L}_{\,mod} & ~ & \downarrow^{L}_{\,mod} &  ~  &
\downarrow^R_{\,mod}& ~ & \downarrow^R_{\,mod}\\[0.3cm]
\AM & \stackrel{\la_{\cA}}{\lrarrow} &
   \Omb(M,\cE) & \stackrel{\cS^{tw}}{\bblarrow} &
\OmC &
\stackrel{\,\la_C}{\llarrow} & \Cp\,,
\end{array}
\label{diag}
\end{equation}
where the DGLAs $\OmT$ and $\OmD$ are taken with the
differentials $D$ and $D+\pa$, respectively, where
as the DGLA modules $\OmE$ and $\OmC$ are
taken with the differentials $D$ and $D+\mb$,
respectively.

Let $f$ be a diffeomorphism of the
pairs $(M,\n)$, $(\tM, \tn)$
$$
f: (M,\n) \mapsto (\tM, \tn)\,,
$$
where $M$ and $\tM$ are $d$-dimensional
manifolds and $\n$, $\tn$ are torsion free
connections on $M$ and $\tM$, respectively.

It is obvious that the corresponding
isomorphisms between the DGLA modules
$$
f_* : (T^{\bul}_{poly}(M), \cA^{\bul}(M)) \mapsto
(T^{\bul}_{poly}(\tM), \cA^{\bul}(\tM))\,,
$$
$$
f_* : (D^{\bul}_{poly}(M), C^{poly}_{\bul}(M)) \mapsto
(D^{\bul}_{poly}(\tM), C^{poly}_{\bul}(\tM))\,,
$$
$$
f_* : (\OmT, \OmE) \mapsto
(\Omb(\tM, \cTp), \Omb(\tM, \cE))\,,
$$
and
$$
f_* : (\OmD, \OmC) \mapsto
(\Omb(\tM, \cDp), \Omb(\tM, \cCp))
$$
commute with the maps in the diagrams (\ref{diag-T-D})
for $M$ and $\tM$.

Furthermore, since the term (\ref{Gamma})
of the Fedosov connection form (\ref{d+B}) does not enter
the definition of the $\Linf$-morphism $\cK^{tw}$
(\ref{cK-tw}) (resp. the morphism of $\Linf$-modules $\cS^{tw}$
(\ref{cS-tw})) the isomorphism $f_*$
commutes with $\cK^{tw}$, and $\cS^{tw}$ as well.
Thus the terms and
the morphisms of diagram (\ref{diag}) are functorial for
diffeomorphisms
of pairs $(M,\n)$, where $\n$ is a torsion
free connection on $M$.

Theorem \ref{thm-chain} is proved. $\Box$

%
%%
%%
%% Kontsevich's story in 1/2 MILE
%%
%%
%

\section{Kontsevich's formality theorem revisited}
\label{revisited}
In this section I prove the existence of a
quasi-isomorphism from $\Tp$ to $\Dp$ which
is functorial for diffeomorphisms of
pairs $(M,\n)$, where $\n$ is a torsion free
connection on $M$\,. Although a proof
of this statement is outlined in
Appendix $3$ of \cite{K1}
some people \cite{C1} think that my proof
is more thorough and refer to my paper \cite{CEFT}
instead of \cite{K1}.

First, I observe that composing the quasi-isomorphisms
$\la_T$ and $\cK^{tw}$ one can shorten the
upper row in diagram (\ref{diag}) to
\begin{equation}
\label{upper}
\begin{array}{ccccc}
\Tp & \stackrel{\cU}{\brarrow} & (\OmD, D+\pa, [,]_G) &
\stackrel{\,\la_D}{\llarrow} & \Dp\,,
\end{array}
\end{equation}
in which $\cU$ is a quasi-isomorphism of DGLAs.

On the other hand due to proposition \ref{Pr-la-D}
the DGLA $\Dp$ is isomorphic to
the DG Lie subalgebra
\begin{equation}
\label{sub-DGLA}
\ker D \cap \G(M, \cDp) \subset \OmD\,.
\end{equation}
This observation raises the question
of whether one can contract the
quasi-isomorphism $\cU$ in (\ref{upper})
to a quasi-isomorphism
\begin{equation}
\label{cU-c}
\cU^c : \Tp \brarrow \ker D \cap \G(M, \cDp^{\bul})
\end{equation}
in a functorial way with respect to
diffeomorphisms of the pair $(M,\n)$\,.
The following theorem gives a positive
answer to this question:
\begin{teo}[M. Konstevich, \cite{K1}, construction 4]
\label{Konets}
For smooth real ma\-ni\-folds $M$
there exists a construction of
 DGLA quasi-isomorphisms
\begin{equation}
\label{cU-K}
\cU^K : \Tp \brarrow \Dp
\end{equation}
which is functorial for diffeomorphisms
of pairs $(M,\n)$, where $\n$ is a
(torsion free) connection on $M$\,.
\end{teo}
{\bf Proof.} First, I construct
a collection of quasi-isomorphisms $(n\ge 0)$
\begin{equation}
\label{upper1}
\cU^{(n)} \,:\, \Tp \brarrow  (\OmD, D+\pa, [,]_G)
\end{equation}
satisfying the following properties:
\begin{equation}
\label{p}
\cU^{(n)}_m (\wedge^m \Tp) \subset
\ker D \cap \G(M, \cDp^{\bul})\,,
\qquad \forall ~ m \le n\,,
\end{equation}
\begin{equation}
\label{p1}
\cU^{(n-1)}_m = \cU^{(n)}_m \,,
\qquad \forall ~ m < n\,,
\end{equation}
where $\cU^{(n)}_m$ denote the
structure maps of $\cU^{(n)}$\,.

I start with $\cU^{(0)}=\cU$ and observe that
due to (\ref{q-iso})
\begin{equation}
\label{D-cU}
(D + \pa )\, \cU^{(0)}_1(\ga)=0\,, \qquad \forall~ \ga\in \Tp\,.
\end{equation}

Since the map $\Phi$ (\ref{Phi}), (\ref{nado-mne})
satisfies equation (\ref{eq:Phi}) I conclude that
for any $\ga \in \Tp$ the combination
$$
\cU^{(0)}_1(\ga) - (D+\pa)\,\Phi (\cU^{(0)}_1(\ga))
$$
does not have the top exterior degree
component.
Thus, applying lemma \ref{styag}
from section \ref{homotopy} for $n=1$
I get a quasi-isomorphism
$$
\tcU \,:\, \Tp \brarrow  (\OmD, D+\pa, [,]_G)\,,
$$
such that for any $\ga\in \Tp$
$$
\tcU_1(\ga)\in \bigoplus_{k=1}^{d-1} \Om^{k}(M,\cDp^{\bul})\,,
$$
where $d=\dim M$\,.

Proceeding in this way I construct
a quasi-isomorphism of DGLAs
$$
\cU^{(1)} \,:\, \Tp \brarrow (\OmD, D+\pa, [,]_G)\,,
$$
such that for any $\ga\in \Tp$
$$
\cU^{(1)}_1(\ga)\in \G(M,\cDp^{\bul})\,.
$$
On the other hand due to
equation (\ref{q-iso})
$$
(D + \pa ) \, \cU^{(1)}_1(\ga)=0
$$
and hence, $\cU^{(1)}_1(\ga)$ belongs to the kernel of
the Fedosov differential $D$. Thus
$\cU^{(1)}$ satisfies (\ref{p}).

Suppose that I have already constructed
$\cU^{(k)}$ up to $k=n-1$ satisfying
(\ref{p}) and (\ref{p1})\,.
Due to (\ref{q-iso})
$$
(D+\pa)\, \cU^{(n-1)}_n(\ga_1, \ga_2, \ldots, \ga_n) =
$$
\begin{equation}
=\frac12 \sum_{p,q\ge 1}^{p+q=n}
 \sum_{\ve\in Sh(p,q)}
\pm [\cU^{(n-1)}_p (\ga_{\ve_1}, \ldots, \ga_{\ve_p}),
\cU^{(n-1)}_q (\ga_{\ve_{p+1}}, \ldots,
\ga_{\ve_{n}})]-
\label{q-iso-cU}
\end{equation}
$$
-\sum_{i\neq j}
\pm \cU^{(n-1)}_{n-1}
([\ga_i,\ga_j]_{SN}, \ga_1, \ldots, \hat{\ga_i}, \ldots, \hat{\ga_j}, \ldots \ga_n),
\qquad \ga_i \in \Tp\,.
$$
By the assumption (\ref{p}) of induction
the right hand side of equation (\ref{q-iso-cU}) is of exterior
degree zero. Hence, using the map $\Phi$ (\ref{Phi}),
(\ref{eq:Phi}), (\ref{nado-mne}) once again I conclude that
for any collection $\ga_i \in \Tp$ the
combination
$$
\cU^{(n-1)}_n (\ga_1, \dots, \ga_n) -
(D+\pa)\, \Phi (\cU^{(n-1)}_n (\ga_1, \dots, \ga_n))
$$
does not have the top exterior degree
component.

Thus applying lemma \ref{styag} enough times
I get a quasi-isomorphism of DGLAs
$$
\cU^{(n)} \,:\, \Tp \brarrow (\OmD, D+\pa, [,]_G)\,,
$$
such that for any $\ga\in \Tp$
$$
\cU^{(n)}_{n}(\ga)\in \G(M,\cDp^{\bul})\,,
$$
and for any $m< n$
$$
\cU^{(n)}_m = \cU^{(n-1)}_m\,.
$$
Due to the corresponding version of
(\ref{q-iso-cU}) $\cU^{(n)}_{n}(\ga)$
is also annihilated by the Fedosov
differential $D$.

Thus, I have constructed the desired
collection (\ref{upper1}). The projective
limit of this inverse system gives me
a quasi-isomorphism $\cU^c$ (\ref{cU-c}).
Composing it with $\la^{-1}_D$ I get
the desired quasi-isomorphism $\cU^K$
(\ref{cU-K})\,.

The construction of $\cK^{tw}$  (\ref{K-tw})
is functorial for diffeomorphisms of pairs
$(M,\n)$ since the term (\ref{Gamma}) of
the Fedosov connection form (\ref{d+B})
does not
contribute to (\ref{cK-tw}).
Thus, the construction of (\ref{cU-K}) is
functorial for diffeomorphisms of pairs
$(M,\n)$ since so are the constructions of
$\tau$ (\ref{lift}), $\Phi$ (\ref{Phi}),
$\la_D$ (\ref{la-D}), and $\la_T$ (\ref{la-T}). $\Box$

%
%%
%%
%% Applications in 1/2 MILE
%%
%%
%

\section{Applications}
The first obvious applications of the formality theorem for
$\Cp$ are related to computation of
Hochschild homology for the quantum algebra of
functions on a Poisson manifold and to description
of traces on this algebra. These applications
were suggested in Tsygan's paper \cite{Tsygan}
(see the first
part of corollary $4.0.3$ and corollary $4.0.5$) as
immediate corollaries of the conjectural formality
theorem (conjecture $3.3.1$ in \cite{Tsygan}).

Let $M$ be, as above, a smooth manifold. Recall that
\begin{defi}[\cite{Bayen, Ber}]
\label{D-Q}
A deformation quantization  of $M$ is a Maurer-Cartan
element (\ref{MC}) $\Pi \in \h D^1_{poly}(M)[[\h]]$
of the DGLA $\Dp[[\h]]$\,, where $\h$ is
an auxiliary variable which plays the role
of the deformation parameter.
Furthermore, two deformation quantizations
$\Pi$ and $\tPi$ are called equivalent
if they are connected by the action
(\ref{act-onMC}) of an element $U$ in
the prounipotent group
$$
\mG_D = \{\, I+ \h D^0_{poly}(M)[[\h]] \, \}
$$
corresponding to the
Lie algebra $\h D^0_{poly}(M) [[\h]]$
\end{defi}
Notice that, since the DGLA $\Dp[[\h]]$ is
endowed with the complete filtration given
by degrees in $\h$ the above definition
makes sense.

In plain English, the Maurer-Cartan element
$\Pi$ in the above definition gives rise to
an associative product $*$ (the so-called
star-product) in the algebra $\OM[[\h]]$
\begin{equation}
\label{Pi-star}
a * b = a \cdot b + \Pi(a,b)\,, \qquad a,b\in \OM[[\h]]\,,
\end{equation}
which deforms the ordinary commutative
multiplication in $\OM[[\h]]$\,.
Moreover, two deformation quantizations
$\Pi$ and $\tPi$ corresponding to
star-products $*$ and $\widetilde{*}$
are equivalent if there
exists a formal series of differential
operators (the element $U$ in $\mG_D$)
$$
U = I + \h U_1 + \h^2 U_2 + \dots  \in \{I+ \h D^0_{poly}(M) [[\h]]\}\,,
$$
which establishes an isomorphism between
the algebras $(\OM[[\h]], *)$ and
$(\OM[[\h]]$, $\widetilde{*})$
\begin{equation}
\label{equiv}
U(a * b) = U(a)\, \widetilde{*} \, U(b)\,,
\qquad a,b \in \OM[[\h]]\,.
\end{equation}
{\bf Remark.} Sometimes it is required
that the Maurer-Cartan element $\Pi$ belongs
to the subalgebra of normalized
polydifferential operator. This requirement
corresponds to the compatibility condition
with the unit function:
$$
a * 1 = 1 * a =0\,.
$$
However, since the subcomplex of normalized
Hochschild chains is quasi-isomorphic
to the total complex it is
very easy to switch from one definition
to another using the action of
the group $\{I+ \h D^0_{poly}(M) [[\h]]\}$\,.

Thanks to quasi-isomorphisms of the
upper row in the diagram (\ref{diag})
and proposition \ref{Fukaya} I have a
bijective correspondence between the moduli
space of Maurer-Cartan elements of
the DGLA $\Tp[[\h]]$ of polyvector fields
and the moduli space of Maurer-Cartan elements
of the DGLA $\Dp[[\h]]$ of polydifferential
operators (tensored with $\bbR[[\h]]$).
In other words, if
we consider the cone
\begin{equation}
\label{cone}
\begin{array}{c}
\al  = \h \al_1 + \h^2 \al_2 +
\h^3 \al_3 + \dots \in \h\, T^1_{poly}(M)[[\h]]\,, \\[0.3cm]
[\al , \al ]_{SN}=0
\end{array}
\end{equation}
acted upon by the Lie
algebra $\h\, \G(M, TM)[[\h]]$
\begin{equation}
\label{action}
\al \to [u, \al]_{SN}\,, \qquad
u \in \h\, \G(M, TM)[[\h]]
\end{equation}
then
\begin{cor}[M. Kontsevich, \cite{K}, theorem 1.1]
The deformation quantizations (\ref{Pi-star}) of $M$
modulo the equivalence relation (\ref{equiv})
are in a bijective correspondence
with the points of the cone (\ref{cone})
modulo the action (\ref{action}) of
the prounipotent group corresponding to the
Lie algebra $\h\, \G(M, TM)[[\h]]$\,. $\Box$
\end{cor}
An orbit $[ \al ]$ on the cone (\ref{cone})
corresponding to a deformation $\Pi$
is called {\it Kontsevich's class of the
deformation} $\Pi$ and
any point $\al$ of this orbit is called
{\it a representative} of this class.

Theorem \ref{thm-chain} allows me to
describe Hochschild homology of
the algebra $(\OM[[\h]], *)$ for
any deformation quantization $\Pi$
(\ref{Pi-star}) of $M$\,.
Namely\footnote{see the first part of
corollary $4.0.3$ in \cite{Tsygan}}
\begin{cor}
If $\Pi$ is a deformation quantization
and $\al\in \h \,T^1_{poly}(M)[[\h]]$
represents Kontsevich's class $[\al]$
of $\Pi$
then the complex of Hochschild homology
\begin{equation}
\label{Hoch-hom}
(\Cp[[\h]], \mb + R_{\Pi})
\end{equation}
is quasi-isomorphic to the complex of
exterior forms
\begin{equation}
\label{Forms}
(\AM[[\h]], L_{\al})
\end{equation}
with the differential $L_{\al}$.
\end{cor}
Here, as above, $R$ denotes the action
(\ref{cochain-act}) of Hochschild cochains
on Hochschild chains and $L$
stands for the Lie derivative
(\ref{Tpoly-act})\,.\\
{\bf Remark.}
In the symplectic case the above corollary
reduces to the well-known theorem
of R. Nest and B. Tsygan (theorem A$2.1$ in \cite{NT})
which is proved for the quantum algebra
of compactly supported functions of
a smooth symplectic manifold.
An equivariant version of this
result in the symplectic case
is discussed in paper \cite{Hoch}
(see proposition $4$) and
paper \cite{Tang1} (see theorem $5.2$).

~\\
{\bf Proof.}
Since $\al$ represents Kontsevich's class of the
deformation quantization $\Pi$ the Maurer-Cartan elements
$\la_D (\Pi)$ and
\begin{equation}
\label{mP}
\mP=\sum_{m=1}^{\infty}\frac1{m!}
\cK^{tw}_m(\la_T(\al), \dots, \la_T(\al))
\end{equation}
are connected by the action (\ref{act-onMC})
of an element $\mU$ of the prounipotent
group $\mH$ corresponding to
the Lie algebra
$$
\mh=(\Om^0(M,\cD^0_{poly})\oplus
\Om^1(M,\cD^{-1}_{poly}))\otimes \h \bbR[[\h]]\,.
$$
Therefore $\mU$ provides me with
a quasi-isomorphism
(actually isomorphism)
\begin{equation}
\label{qis-mU}
R_{\mU} : (\OmC[[\h]], D +\mb + R_{\la_D(\Pi)})
\mapsto
(\OmC[[\h]], D + \mb + R_{\mP})
\end{equation}
from the complex
$(\OmC[[\h]], D + \mb + R_{\la_D(\Pi)})$
to the complex
$(\OmC[[\h]], D+ \mb + R_{\mP})$\,.

Twisting the terms in the second diagram in (\ref{diag-T-D})
by the Maurer-Cartan element $\Pi$
I get the new commutative diagram
\begin{equation}
\begin{array}{ccc}
(\OmD[[\h]], D+\pa+[\la_D(\Pi),\cdot \,]_G) &
\stackrel{\la_D}{\llarrow} &
(\Dp[[\h]],\pa+[\Pi,\cdot \,]_G)\\[0.3cm]
\downarrow^{R}_{\,mod}  & ~  &     \downarrow^{R}_{\,mod} \\[0.3cm]
(\OmC[[\h]], D+\mb+R_{\la_D(\Pi)}) &\stackrel{\la_C}{\llarrow} &
(\Cp[[\h]], \mb + R_{\Pi})\,,
\end{array}
\label{diag-Pi}
\end{equation}
in which the DGLAs $\OmD[[\h]]$ and $\Dp[[\h]]$
carry the initial Lie bracket $[,]_G$
(\ref{Gerst}).

Due to claim $5$ of proposition \ref{twist-mod}
the map $\la_C$ in the above diagram is
a quasi-isomorphism of complexes.

On the other hand twisting the terms in the left part
of diagram (\ref{diag}) by the Maurer-Cartan element
$\al\in \Tp[[\h]]$ I get
the new commutative diagram
\begin{equation}
\begin{array}{ccc}
(\Tp[[\h]], [\al,\cdot ]_{SN}) &  \stackrel{\la_T}{\lrarrow}
& (\OmT[[\h]], D + [\la_T(\al),\cdot ]_{SN} ) \\[0.3cm]
\downarrow^{L}_{\,mod} & ~ & \downarrow^{L}_{\,mod} \\[0.3cm]
(\AM[[\h]], L_{\al}) & \stackrel{\la_{\cA}}{\lrarrow} &
(\Omb(M,\cE)[[\h]], L_{\la_T(\al)} )\,,
\end{array}
\label{diag-al-tw}
\end{equation}
in which the DGLAs $\OmT[[\h]]$ and $\Tp[[\h]]$
carry the initial Lie bracket $[,]_{SN}$
(\ref{eq:Leib}).

Due to claim $5$ of proposition \ref{twist-mod}
the map $\la_{\cA}$ in diagram (\ref{diag-al-tw}) is
a quasi-isomorphism of complexes.

Similarly, twisting the terms in the middle part
of diagram (\ref{diag}) by the Maurer-Cartan element
$\la_T(\al)\in \OmT[[\h]]$ I get
\begin{equation}
\begin{array}{ccc}
 (\Omb(\cTp)[[\h]], D + [\la_T(\al),\cdot ]_{SN} )&
\stackrel{\cK^{\al}}{\brarrow} &
(\Omb(\cDp)[[\h]], D +\pa + [\mP,\cdot ]_{G} )\\[0.3cm]
\downarrow^{L}_{\,mod} &  ~  &
\downarrow^R_{\,mod} \\[0.3cm]
(\Omb(\cE)[[\h]], L_{\la_T(\al)} ) & \stackrel{\cS^{\al}}{\bblarrow} &
(\Omb(\cCp)[[\h]], D+\mb+R_{\mP} )\,,
\end{array}
\label{diag-al-tw1}
\end{equation}
where $\cK^{\al}$ and $\cS^{\al}$ are obtained
from $\cK^{tw}$ and $\cS^{tw}$, respectively,
by twisting via $\la_T(\al)\in \OmT[[\h]]$\,,
$\mP$ is defined in (\ref{mP}), $\OmT[[\h]]$
goes with the initial bracket $[,]_{SN}$, and
$\OmD[[\h]]$ goes with the initial
bracket $[,]_G$\,.

Again, due to claim $5$ of proposition \ref{twist-mod}
the morphism of $\Linf$-modules
$\cS^{\al}$ in diagram (\ref{diag-al-tw1}) is
a quasi-isomorphism.

The desired statement is proved
since the complexes (\ref{Hoch-hom})
and (\ref{Forms}) are connected by
a chain of quasi-isomorphisms. $\Box$

Another application of theorem $1$ is related to description
of traces on the algebra $(C^{\infty}_c(M)[[\h]], *)$, where
by $C^{\infty}_c(M)$ I denote the vector space of smooth functions
with a compact support.

By definition {\it a  trace} is a continuous $\bbR[[\h]]$-linear
$\bbR[[\h]]$-valued functional $tr$ on $C^{\infty}_c(M)[[\h]]$
vanishing on commutators
$$
tr(R_{\Pi}(a))=0\,,
$$
where $a=a(x_0,x_1)$ is a function in $C^{\infty}(M\times M)$
with a compact support in its first argument, and
$R$ is as in (\ref{cochain-act}).

One can easily verify that my constructions still make
sense if I replace the first version (\ref{H-chains1}) of $\Cp$
 by
$$
C^{poly-com}(M)= \bigoplus_{n\ge 0} C^{\infty}_{com}(M^{n+1})\,,
$$
and the vector space of exterior forms
$\AM$ by the vector space $\cA^{\bul}_c(M)$ of
exterior forms with a compact support.
Here by $ C^{\infty}_{com}(M^{n+1})$ I denote the vector
space of smooth functions on $M^{n+1}$ with a compact
support in the first argument.

Then the corresponding version of the above
corollary implies that
\begin{cor}[\cite{Tsygan}, Corollary 4.0.5]
\label{trace}
If $\Pi$ is deformation quantization
(\ref{Pi-star}) and $\al$ represents
Kontsevich's class of $\Pi$ then
the vector space of traces on the algebra
$(C^{\infty}_c(M)[[\h]],*)$ is isomorphic to
the vector space of continuous $\bbR[[\h]]$-linear
$\bbR[[\h]]$-valued functionals on $C^{\infty}_c(M)[[\h]]$
vanishing on all Poisson brackets $\al(a,b)$
for $a,b\in C^{\infty}_c(M)[[\h]]$\,. $\Box$
\end{cor}
For a symplectic manifold this statement
has been proved in \cite{CFS, Fedosov1, NT}.

~\\
{\bf Remark.} Corollary \ref{trace} still
holds if one replaces real valued functions (resp. traces)
by smooth complex valued functions (resp. complex valued
traces), as well as the ring $\bbR[[\h]]$ by
the field $\bbC[[\h, \h^{-1}]$.

I would like to mention that the functoriality
of the chain of quasi-isomor\-phisms (\ref{diag-thm})
in theorem \ref{thm-chain} implies the
following interesting results
\begin{cor}
Let $M$ be a smooth manifold
equipped with a smooth action
of a group $G$. If one can construct
a $G$-invariant connection $\n$ on $M$
then there exists a chain of $G$-equivariant
quasi-isomorphisms between the
DGLA modules $(\Tp, \AM)$ and
$(\Dp, \Cp)$\,.  $\Box$
\end{cor}
In particular,
\begin{cor}
If $M$ is equipped with a smooth action
of a finite or compact group $G$
then the DGLA modules $\Big( (\Tp)^G$, $(\AM)^G\Big)$
and $\Big( (\Dp)^G$, $(\Cp)^G \Big)$ of
$G$-invariants are quasi-isomorphic. $\Box$
\end{cor}

Due to the functoriality of the quasi-isomorphism
(\ref{cU-K}) in theorem \ref{Konets} I have
the following result:
\begin{cor}
If $M$ is equipped with a smooth action
of a finite or compact group $G$
then there exists a quasi-isomorphism
from the DGLA $(\Tp)^G$ of $G$-invariant
polyvector fields to the DGLA $(\Dp)^G$
$G$-invariant polydifferential
operators on $M$. $\Box$
\end{cor}

Using this corollary and proposition
\ref{Fukaya} I get a solution of a deformation quantization
problem for an arbitrary Poisson orbifold. Namely,
\begin{cor} Given a smooth action of a finite group
$G$ on a ma\-ni\-fold $M$ and a $G$-invariant Poisson
structure $\al\in (\wedge^2 T M)^G$ one can always
construct a $G$-invariant star-product $*$\,,
corresponding to $\al$\,.
Furthermore, $G$-invariant star-products on $M$ corresponding
to the Poisson bracket $\al$ are classified up to
equivalence by non-trivial $G$-invariant deformations
of $\al$\,. $\Box$
\end{cor}

%%
%% CHAPTER 5 ENDS
%%
%% Conclusion in  1 MILE
%%
%%

\chapter{Conclusion}
I am glad that the results of my thesis
have been already applied to two interesting
problems. In his PhD thesis
\cite{Tang} X. Tang used theorem \ref{thm-chain}
to compute Hochschild homology of
formal symplectic deformations of proper
\'etale Lie
groupoids\footnote{See paper \cite{Tang1}
in which cyclic and Hochschild
(co)homology groups of formal symplectic
deformations of proper \'etale Lie groupoids
were computed without making use of
formality theorems.} and
in \cite{C1} D. Calaque used
theorem \ref{Konets} in order to
quantize a class of formal classical
dynamical $r$-matrices in the reductive case.

I would like to mention paper \cite{C} in which
 D. Calaque generalized theorem \ref{Konets} to the case when
the tangent bundle of $M$ is replaced by
an arbitrary smooth Lie algebroid.
In our joint paper \cite{CDH} we generalized
the result of \cite{C} to Hochschild chains
and extended our constructions to
the holomorphic setting. In this way we
proved a version of
Tsygan's formality conjecture
for Hochschild chains for
any complex manifold.

I would like to mention parallel
results of the MIT alumnus A. Yekutieli.
In his papers \cite{Ye, Ye1, Ye11} he
proved that for any smooth algebraic
variety $X$ over a field
$\bbK$ ($\bbR \subset \bbK$) the sheaf
of polyvector fields and
the sheaf of polydifferential
operators are quasi-isomorphic as
sheaves of DGLAs.

In \cite{Tomsk} S.L. Lyakhovich and A.A. Sharapov
suggested that a generalization of theorem \ref{Konets}
for super-manifolds can be applied to
quantum reduction. In this paper they proposed
the most general setting of a reduction
which leads to a Poisson manifold and showed
that under certain cohomological
conditions the ``super''-version of
Kontsevich's formality theorem would lead to
deformation quantization of the reduced
manifold.
In \cite{Tomsk} the authors also discussed
a possible path integral approach
\cite{BLN,CF} to the ``super''-version of
Kontsevich's formality theorem.

Two relative versions of Kontsevich's
formality theorem were suggested
simultaneously in papers \cite{Gilles}
and \cite{Rel-FT1}. In both of these papers
the authors considered a smooth submanifold
$C$ of a smooth manifold $M$. In paper
\cite{Gilles} it is conjectured that
the DGLA (and more generally
$G_{\infty}$-algebra) of polydifferential
operators compatible with the
ideal $I\subset \OM$ of functions vanishing
on $C$ is formal. In \cite{Gilles} the authors
proved this conjecture for the case $M=\bbR^d$ and
$C = \bbR^{d-k}$ if $k\ge 2$, and
applied this result to the construction of
representations of the star-product
algebras. In paper \cite{Rel-FT1} the authors
proved the formality theorem for the
DGLA of polydifferential operators
acting on the exterior algebra of the
conormal bundle of $C$ in $M$ and applied
this result to the quantum reduction.
I would like to mention that the
question of globalization is not
properly addressed in \cite{Rel-FT1}.
However, I do not think that
this question is very difficult since
the authors reduced their problem to
a formal neighborhood of $C$ in $M$\,.

There are still many interesting
open questions in this subject.
For example, it would be very interesting to develop
the applications \cite{Rel-FT1},
\cite{Tomsk} of the ``super''-version
of Kontsevich's formality theorem
to the quantum reduction and find out how
the characteristic classes of deformations
fit into the reduction procedure \cite{Cherviv}.
It is also interesting to further
investigate the relative versions
\cite{Gilles, Rel-FT1} of Kontsevich's
formality theorem and apply them
to the conjectural correspondence
\cite{Rel-FT} between the category of
Poisson manifolds with dual pairs as morphisms and
the category of deformation
quantization algebras with bimodules
as morphisms. Finally, the cyclic formality conjecture
\cite{Tsygan} as well as the most general
version of the algebraic
index theorem \cite{TT} still remain open
questions.

\appendix
\chapter{Figures}

\clearpage
\newpage

\clearpage
\newpage

\end{document}